\documentclass[10pt]{article}
\usepackage{amssymb,amsmath,amsfonts,bbm,pifont,upgreek,bbold,accents,ulem} 
\usepackage[colorlinks=true]{hyperref}
\hypersetup{urlcolor=blue, citecolor=red}

\font\teneufm=eufm10
\font\seveneufm=eufm7
\font\fiveeufm=eufm5
\newfam\eufmfam
\textfont\eufmfam=\teneufm
\scriptfont\eufmfam=\seveneufm
\scriptscriptfont\eufmfam=\fiveeufm
\def\eufm#1{{\fam\eufmfam\relax#1}}

\newcommand\beq[1]{ \begin{equation}\label{#1} }
\newcommand{\eeq}{ \end{equation} }
\newcommand{\beqno}{ \[ }
\newcommand{\eeqno}{ \] }
\newcommand\beqa[1]{ \begin{eqnarray} \label{#1}}
\newcommand{\eeqa}{ \end{eqnarray} }
\newcommand{\beqano}{ \begin{eqnarray*} }
\newcommand{\eeqano}{ \end{eqnarray*} }
\newcommand\arr[1]{\left\{\begin{array}{l}#1\end{array}\right.}
\renewcommand{\theequation}{\arabic{section}.\arabic{equation}}
\newtheorem{conjecture}{Conjecture}[section]
\newtheorem{theorem}{Theorem}[section]
\newtheorem{definition}{Definition}[section]
\newtheorem{proposition}{Proposition}[section]
\newtheorem{lemma}{Lemma}[section]
\newtheorem{sublemma}{Sublemma}[section]
\newtheorem{remark}{Remark}[section]
\newtheorem{notationalremark}{Notational Remark}[section]
\newtheorem{corollary}{Corollary}[section]
\newtheorem{assumption}{Assumption}[section]
\newtheorem{claim}{Claim}[section]
\newtheorem{tools}{$\negsp\negsp$}[subsection]

\newcommand\thm[1]{ \begin{theorem}\label{#1}}
\newcommand\thmtwo[2]{ \begin{theorem}[#1]\label{#2}}
\newcommand\ethm{ \end{theorem} }
\newcommand\dfn[1]{ \begin{definition}\label{#1} \rm}
\newcommand\dfntwo[2]{ \begin{definition}[#1]\label{#2} \rm}
\newcommand\edfn{ \end{definition} }
\newcommand\pro[1]{ \begin{proposition}\label{#1}}
\newcommand\protwo[2]{ \begin{proposition}[#1]\label{#2}}
\newcommand\epro{ \end{proposition} }
\newcommand\lem[1]{ \begin{lemma}\label{#1}}
\newcommand\lemtwo[2]{ \begin{lemma}[#1]\label{#2}}
\newcommand\elem{ \end{lemma} }
\newcommand\sublem[1]{ \begin{sublemma}\label{#1}}
\newcommand\sublemtwo[2]{ \begin{sublemma}[#1]\label{#2}}
\newcommand\esublem{ \end{sublemma} }
\newcommand\rem[1]{ \begin{remark}\label{#1} \rm}
\newcommand\erem{ \end{remark} }
\newcommand\notrem[1]{ \begin{notationalremark}\label{#1} \rm}
\newcommand\enotrem{ \end{notationalremark} }
\newcommand\cor[1]{ \begin{corollary}\label{#1}}
\newcommand\cortwo[2]{ \begin{corollary}[#1]\label{#2}}
\newcommand\ecor{ \end{corollary} }
\newcommand\asmp[1]{ \begin{assumption}\label{#1}}
\newcommand\asmptwo[2]{ \begin{assumption}[#1]\label{#2}}
\newcommand\easmp{ \end{assumption} }
\newcommand\clm[1]{ \begin{claim}\label{#1}}
\newcommand\eclm{ \end{claim} }
\newcommand{\proof}{\par\medskip\noindent{\bf Proof\ }}
\newcommand\equ[1]{{\rm (\ref{#1})}}
\expandafter\chardef\csname pre amssym.def
at\endcsname=\the\catcode`\@
\catcode`\@=11
\def\undefine#1{\let#1\undefined}
\def\newsymbol#1#2#3#4#5{\let\next@\relax
 \ifnum#2=\@ne\let\next@\msafam@\else
 \ifnum#2=\tw@\let\next@\msbfam@\fi\fi
 \mathchardef#1="#3\next@#4#5}
\def\mathhexbox@#1#2#3{\relax
 \ifmmode\mathpalette{}{\m@th\mathchar"#1#2#3}%
 \else\leavevmode\hbox{$\m@th\mathchar"#1#2#3$}\fi}
\def\hexnumber@#1{\ifcase#1 0\or 1\or 2\or 3\or 4\or 5\or 6\or 7\or
8\or
 9\or A\or B\or C\or D\or E\or F\fi}
\ifcase\@ptsize
 \font\tenmsb=msbm10
 \font\sevenmsb=msbm7
 \font\fivemsb=msbm5
\or
 \font\tenmsb=msbm10 scaled \magstephalf
 \font\sevenmsb=msbm7 scaled \magstephalf
 \font\fivemsb=msbm5  scaled \magstephalf
\or
 \font\tenmsb=msbm10 scaled \magstep1
 \font\sevenmsb=msbm7 scaled \magstep1
 \font\fivemsb=msbm5 scaled \magstep1
\fi
\newfam\msbfam
\textfont\msbfam=\tenmsb
\scriptfont\msbfam=\sevenmsb
\scriptscriptfont\msbfam=\fivemsb
\edef\msbfam@{\hexnumber@\msbfam}
\def\Bbb#1{\fam\msbfam\relax#1}
\def\widehat#1{\setboxz@h{$\m@th#1$}%
 \ifdim\wdz@>\tw@ em\mathaccent"0\msbfam@5B{#1}%
 \else\mathaccent"0362{#1}\fi}
\def\widetilde#1{\setboxz@h{$\m@th#1$}%
 \ifdim\wdz@>\tw@ em\mathaccent"0\msbfam@5D{#1}%
 \else\mathaccent"0365{#1}\fi}

\def\RIfM@{\relax\ifmmode}
\def\nonmatherr@#1{\errmessage{\string#1\space allowed only in math mode}}
\def\Bbb{\RIfM@\expandafter\Bbb@\else
 \expandafter\nonmatherr@\expandafter\Bbb\fi}
\def\Bbb@#1{{\Bbb@@{#1}}}
\def\Bbb@@#1{\fam\msbfam\relax#1}
\def\setboxz@h{\setbox\z@\hbox}
\def\wdz@{\wd\z@}
\catcode`\@=\csname pre amssym.def at\endcsname
\newcommand{\ie}{{\rm i.e.\,}}
\newcommand{\eg}{{\it e.g.\,}}

\newcommand{\Giu}{{\bigskip\noindent}}
\newcommand{\nl}{{\smallskip\noindent}}
\newcommand{\noi}{{\noindent}}

\newcommand{\qed}{\hskip.5truecm
\vrule width 1.7truemm height 3.5truemm depth 0.truemm
\par\Giu}
\newcommand{\qedeq}{\hskip.5truecm
\vrule width 1.7truemm height 3.5truemm depth 0.truemm}

\newcommand{\negsp}{\hspace{-.09truecm}}  

\newcommand{\dst}{\displaystyle}
\newcommand\ovl[1]{ \overline {#1} }

\newcommand\su[1]{ \frac{1}{ {#1}} }

\newcommand{\torus}{ {\Bbb T}   }
\renewcommand{\natural}{ {\Bbb N}   }
\newcommand{\real}{ {\Bbb R}   }
\newcommand{\integer}{ {\Bbb Z}   }

\renewcommand{\a }{ {\alpha}   }
\renewcommand{\b}{ {\beta}   }
\newcommand{\g}{ {\gamma}   }
\newcommand{\G}{ {\Gamma}   }
\renewcommand{\d}{ {\delta}   }

\newcommand{\e }{ {\epsilon}   }

\renewcommand{\th }{ {\theta}   }

\renewcommand{\k}{ {\kappa}   }
\renewcommand{\l}{ {\lambda}   }
\renewcommand{\L}{ {\Lambda}   }
\newcommand{\m}{ {\mu}   }
\newcommand{\n}{ {\nu}   }

\newcommand{\p}{ {\pi}   }
\renewcommand{\P}{ {\Pi}   }
\renewcommand{\r}{ {\rho}   }
\newcommand{\s}{ {\sigma}   }

\renewcommand{\t}{ {\tau}   }
\newcommand{\f}{ {\varphi}   }

\renewcommand{\o}{ {\omega}   }
\renewcommand{\O}{ {\Omega}   }

\newcommand{\const}{{\, \rm const\, }}

\newcommand{\cA}{ {\cal A} }

\newcommand{\cR}{ {\cal R} }
\newcommand{\cH}{ {\cal H} }
\newcommand{\cK}{ {\cal K} }

\newcommand{\cD}{ {\cal D} }

\newcommand{\cM}{ {\cal M} }

\newcommand{\cP}{ {\cal P} }

\newcommand{\cS}{ {\cal S} }

\newcommand\ppu{{ (1) }}
\newcommand\ppd{{ (2) }}
\newcommand\ppt{{ (3) }}

\newcommand\ppj{{ (j) }}

\newcommand\ppn{{ (n) }}

\newcommand\ppi{{ (i) }}

\newcommand\ppo{{ (0) }}

\newcommand\ul{{\uplambda}}
\newcommand\ux{{\upxi}}
\newcommand\uh{{\upeta}}
\newcommand\up{{\rm p}}
\newcommand\uq{{\rm q}}
\newcommand\uz{{\rm z}}

\newcommand\meas{{\, \rm meas\,}}

\newcommand\fd{{n_1}}
\newcommand\sd{{n_2}}
\newcommand\fg{{\g_1}}
\newcommand\sg{{\bar\g_2}}
\newcommand\pertnorm{{E}}
\newcommand\KAM{{\hat E}}

\begin{document}

\title{
Aspects of the planetary Birkhoff normal form\footnote{Research Supported by ``Prin 2009 project ÒCritical Point Theory and Perturbative Methods for Nonlinear Differential Equations''.}
}

\author{  
Gabriella Pinzari  \\ 
\vspace{-.2truecm}
{\scriptsize Dipartimento di Matematica}
\\{\scriptsize  Universit\`a  ``Roma Tre''}
\vspace{-.2truecm}
\\{\scriptsize Largo S.L. Murialdo 1, I-00146 Roma (Italy) }
\vspace{-.2truecm}
\\{\scriptsize  
pinzari@mat.uniroma3.it}
}
\date{}
\maketitle
\vskip.1in
\noi
\begin{flushright}
{\sl This paper is dedicated to  Professor Alain Chenciner on his 70${}^{\rm th}$ birthday}
\end{flushright}
\vskip.5in
\noi

\begin{abstract}\footnotesize The discovery in \cite{pinzari-th09}, 
 \cite{chierchiaPi11b}  of the Birkhoff normal form for the planetary many--body problem opened new insights and hopes for the comprehension of the dynamics of this problem. Remarkably, it allowed to give a {\sl direct} proof of the celebrated Arnold's Theorem \cite{arnold63} on the stability of planetary motions. In this paper, using a ``ad hoc'' set of symplectic variables, we develop an  asymptotic formula for this normal form that may turn to be useful in applications. As an example, we provide two very simple applications to  the three--body problem: we  prove a conjecture by V. I. Arnold \cite{arnold63} on the {\sl Kolmogorov set} of this problem and, using Nehoro{\v{s}}ev Theory \cite{nehorosev77}, we prove, in the planar case,
stability of {\sl all} planetary actions over  exponentially--long times, provided mean--motion resonances are excluded. We also briefly discuss perspectives and problems for full  generalization of the results in the paper.
 \end{abstract}

\nl
{\footnotesize{\bf Keywords: } Averaging Theory, Birkhoff normal form, Nehoro{\v{s}}ev Theory, Planetary many--body problem, {Arnold's Theorem on the stability of planetary motions}, Properly--degenerate {\sc kam} Theory, steepness.

\nl
{{\bf MSC2000 numbers:}
34D10,  34C20, 70E55,  70F10, 70F15, 70F07,  37J10, 37J15, 37J25, 37J35, 37J40, 70K45

}}

\maketitle
\tableofcontents

\section{
Introduction and results}
\setcounter{equation}{0}
\label{intro}
\renewcommand{\theequation}{\arabic{equation}}


\vskip.2in
\noi
{\bf\large 1.1} 
The planetary many--body problem
consists in determining the dynamics of $(1+n)$ masses
undergoing Newtonian attraction. The term ``planetary'' is reserved to the case when one mass, the ``sun'', or ``star'', denoted with $\bar m_0$, is taken to be much greater than the others, $\m\bar m_1$, $\cdots$, $\m\bar m_n$, which are called ``planets''. Here $\m\ll1$ is a small number.
After the ``heliocentric\footnote{See, \eg, \cite{robutel95}.} reduction'' of  invariance by translations, this dynamical system is governed by the $3n$ degrees of freedom Hamiltonian
\beq{planetary}{\rm H}_{\rm plt}=\sum_{i=1}^n(\frac{|y^\ppi|^2}{2m_i}-\frac{m_iM_i}{|x^\ppi|})+\m\sum_{1\le i<j\le n}(\frac{y^\ppi\cdot y^\ppj}{\bar m_0}-\frac{\bar m_i \bar m_j}{|x^\ppi-x^\ppj|})\eeq
on the phase space
\beqno
(y,x)=(y^\ppu,\cdots, y^\ppn,x^\ppu,\cdots, x^\ppn)\in (\real^3)^{2n}:\quad x^\ppi\ne 0\ ,\quad x^\ppi\ne x^\ppj\eeqno
endowed with the standard 2-- form
$$\O:=dy\wedge dx:=\sum_{i=1}^n \sum_{j=1}^3dy_j^\ppi\wedge dx^\ppi_j$$
where $y^\ppi=(y^\ppi_1,y^\ppi_2,y^\ppi_3)$, $x^\ppi=(x^\ppi_1,x^\ppi_2,x^\ppi_3)$. Here, $m_i$, $M_i$ are suitable auxiliary masses related to $\bar m_i$ and $\m$ via
$$M_i=\bar m_0+\m\bar m_i\qquad m_i=\frac{\bar m_0\bar m_i}{\bar m_0+\m\bar m_i}\ .$$  
A procedure commonly followed in the past \cite{arnold63}, \cite{fejoz04}, \cite{herman09}, \cite{nehorosev77} to regard the system as a ``close to integrable'', was to use a symplectic set of variables, usually called ``Poincar\'e variables''. These   variables, that we denote 
\beqno(\L_i,\ul_i,\uh_i,\ux_i,\up_i,\uq_i)\qquad 1\le i\le n\ ,\eeqno
are  ``six per planet''. They  were introduced by H. Poincar\'e by modifying another set of ``action--angle''  variables $(\L_i,\G_i, \Theta_i, \ell_i, {\rm g}_i,\theta_i)\in \real^3\times \torus^3$  (where $\torus:=\real/(2\p\integer$), having the $\L_i$'s in common, called ``Delaunay variables''.  Delaunay variables are ``natural'',  ``action--angle''  variables  related to the ``{Cartesian variables}'' $(y^\ppi, x^\ppi)$ in \equ{planetary} via the integration of each of the ``two--body'' Hamiltonians
\beqno\frac{|y^\ppi|^2}{2m_i}-\frac{m_iM_i}{|x^\ppi|}\ .\eeqno
 The Poincar\'e variables  are in part ``action--angle'' (\ie, $(\L_i,\ul_i)\in \real\times\torus$), in part ``rectangular'' (\ie, $(\uh_i,\ux_i,\up_i,\uq_i)\in \real^4$). The definition of Delaunay
 and Poincar\'e variables may be found, \eg, in \cite{chierchiaPi11c}.
 In  Delaunay--Poincar\'e variables, any of the two--body Hamiltonian above takes the  ``Kepler form'' 
$$h^\ppi_{\rm Kep}(\L_i)=-\frac{M_i^2m_i^3}{2\L_i^2}\ .$$
It is ``properly degenerate'': two degrees of freedom disappear, as it is well known. This proper degeneracy  naturally reflects on the system \equ{planetary}, which in fact takes  the form
\beq{prop deg}\cH_{\rm P}(\L,\ul,\uz)=h_{\rm Kep}(\L)+\m f_{\rm P}(\L,\ul,\uz)\eeq
where   $h_{\rm Kep}(\L)$ is the $n$ degrees of freedom ``unperturbed'' part $-\sum_{i=1}^n\frac{M_i^2m_i^3}{2\L_i^2}$, while  $f_{\rm P}(\L,\ul,\uz)$ is 
the   $3n$ degrees of freedom  ``perturbation'' 
 \beq{pert**}\sum_{1\le i<j\le n}(\frac{y^\ppi\cdot y^\ppj}{\bar m_0}-\frac{\bar m_i \bar m_j}{|x^\ppi-x^\ppj|})\eeq in \equ{planetary}, expressed in Poincar\'e variables. Here, we have denoted as $(\L,\ul,\uz)$ the $3n$--dimensional collection of
\beq{PV} \L=(\L_1,\cdots,\L_n)\ ,\quad \ul=(\ul_1,\cdots,\ul_n)\ ,\quad \uz=(\uz_1,\cdots,\uz_n)\eeq
with $\uz_i:=(\uh_i,\ux_i,\up_i,\uq_i)$.
 \vskip.2in
\noi

  \nl
 A long outstanding problem lasted about fifty years concerned the existence of a {\sl Birkhoff normal form} for the system \equ{prop deg}.

 \nl
Namely, if it were possible to conjugate the Hamiltonian \equ{prop deg}  to an analogue one, 
 $$\cH_{\rm bnf}=h_{\rm Kep}+\m f_{\rm bnf}\ ,$$
whose average (``secular'') perturbing function
$$(f_{\rm bnf})_{\rm av}:=\frac{1}{(2\p)^n}\int_{\torus^n} f_{\rm bnf}$$ were in Birkhoff normal form of some order (see \cite{hoferZ94}, \cite{arnold89}, \cite{williamson36} for information on Birkhoff theory). The claim is perfectly natural, since in fact the average  $(f_{\rm P})_{\rm av}$ of $f_{\rm P}$ in \equ{prop deg} turns to have an elliptic equilibrium point in $\{\uz=0\}$, for any any choice of $\L$. We recall that, physically,  $\{\uz=0\}$, the ``secular origin'', corresponds to circular and co--inclined unperturbed motions, a  configuration with relevant physical meaning, being commonly observed in nature in many--body systems.

\nl
The problem was settled by V. I. Arnold, who, in the 60's announced (1962's International Congress for Mathematicians; Stockolm, \cite{icm62}) and next  (1963) published his more than celebrated  ``theorem on the stability of planetary motions''; or ``The Planetary Theorem'', for short.
 \begin{theorem}[{V. I. Arnold,  \cite[p. 127]
 {arnold63}}]\label{Planetary Theorem}
In the $n$--body problem there exists a set of initial conditions having positive Lebesgue measure and such that, if the initial positions and velocities belong to this set, the distances of the bodies from each other will remain perpetually bounded.
\end{theorem}
Arnold gave the details of the proof of the Planetary Theorem  the case of  three bodies constrained on a plane: the ``first'' non trivial case. He was aware that, to extend the result to the general problem, some extra--difficulty
related to the 
 ``rotation invariance'' of the system \equ{planetary} was to be overcome. Namely, the invariance by the two--parameter  group of (non--commuting) transformations
\beq{rotation invariance}(y^\ppi, x^\ppi)\to (\cR y^\ppi, \cR x^\ppi)\ ,\quad \cR\in {\rm SO}(3)\ .\eeq
From a dynamical point of view, rotation invariance is caused by the conservation, along the ${\rm H}_{\rm plt}$--trajectories, of the three components, ${\rm C}_1$, ${\rm C}_2$ and ${\rm C}_3$, of the ``angular momentum''
\beq{C}{\rm C}=\sum_{i=1}^n x^\ppi\times y^\ppi\ ,\eeq
where ``$\times$'' denotes skew--product. 

\nl
To prove his Planetary Theorem, Arnold proved an abstract theorem (that he called ``The Fundamental Theorem'', see Appendix  \ref{ArnoldKAM}; Theorem \ref{fundamental theorem})  on the conservation of quasi--periodic motions precisely suited for properly--degenerate systems. For such systems, indeed, ``standard'' non--degeneracy assumptions as the ones appeared in \cite{kolmogorov54}, \cite{moser1962} or  \cite{arnold63c} are strongly violated.
The non--degeneracy condition of the Fundamental Theorem is a ``strong'' non--linearity condition   we  shall refer to as  ``full torsion''. It requires (besides the non--degeneracy of the unperturbed part) the existence of the Birkhoff normal form and the invertibility of the matrix of the coefficients of the second--order term (second--order ``Birkhoff invariants''): see conditions {\rm (ii)}  and {\rm (iii)} in Theorem \ref{fundamental theorem}.

\nl
In the case of the problem in the space, the integral \equ{C} causes another strong degeneracy in the perturbation: one of the first order Birkhoff invariants associated to $(f_{\rm P})_{\rm av}$, $\O_{2n}(\L)$, vanishes identically. This ``resonance'' is apparently a problem for the construction of the Birkhoff normal form. See, \eg, \cite{hoferZ94}.  It is worth to remark that, moreover, another resonance is to be taken in account, which, even though not mentioned in \cite{arnold63}, was later pointed out by M. Robert Herman: the sum of the remaining first invariants $\O_1(\L)$, $\cdots$, $\O_{2n-1}(\L)$, vanishes identically (see \cite{abdullahA01} for a study on Herman resonance). Such two resonances,  
$$\O_{2n}(\L)\equiv0\ ,\quad \sum_{i=1}^{2n-1}\O_i(\L)\equiv 0\ ,$$
are usually referred to, respectively, as ``rotational'', ``Herman'' resonance or, jointly,  ``secular resonances''.

\nl
To overcome the problem of the secular resonances (or, at least, of the rotational one), Arnold proposed, in \cite{arnold63}, a sketchy program of which he did not give the complete details. Such details revealed to be not trivial at all. The  spatial three--body case later was proved in the PhD dissertation by P. Robutel   \cite{robutel95} (see also \cite{laskarR95}), on the basis of a rigorous development of the ideas in \cite{arnold63}.

\nl
The first complete proof of Arnold's Planetary Theorem in the general case appeared in \cite{fejoz04}, including  efforts by M. Robert Herman. This important and beautiful result was reached with a different {\sc kam} technique, avoiding Birkhoff normal form: only the properties of the first order invariants are exploited in  \cite{fejoz04}. 
The underlying elegant,  {\sc kam} Theory  in \cite{fejoz04} (for ``smooth'' systems) is different from the one in \cite{arnold63}; it goes back to \cite{Russmann:2001} (analytic) and exploits non--degeneracy conditions previously studied in the 80's by Arnold, Piartly, Parasyuk, Sprinzuk and others; see   \cite{fejoz04} and references therein, for more information.
Moreover, the problem of secular  resonances is solved in  \cite{fejoz04} via arguments of  abstract reductions, \cite{arnoldKS06}.

\nl
The complete achievement of Arnold's program for $n\ge 3$  was  reached  in the PhD thesis \cite{pinzari-th09}  (next, published in \cite{chierchiaPi10}, \cite{chierchiaPi11a}, \cite{chierchiaPi11b}).

\nl
Switching from $n=2$ (spatial)  to $n\ge 3$ (spatial)  required new  ideas. 
Indeed,   as   Arnold pointed out in \cite[Chapter III, \S 5, 4--5]{arnold63}, while in the spatial three--body case   a classical tool for reducing the integral \equ{C}, the so--called ``Jacobi reduction of the nodes'' \cite{jacobi1842} was available and in fact used in \cite{robutel95}, this tool was instead lacking for the general spatial problem with more than two planets. By this reason,  Arnold suggested  a qualitatively different strategy to handle this latter case. He conjectured 
\cite[Chapter III, \S 5, 5]{arnold63}
it were possible  to reduce only two (out of three) non--commuting components of ${\rm C}$ (or functions of them) and, simultaneously, keep the ({\sl regular}) structure of the Hamiltonian \equ{prop deg} in Poincar\'e variables.
He believed this should
let the system free of the vanishing eigenvalue.
 Note that this seems in contrast with the strategy \cite[Chapter III, \S 5, 4]{arnold63} for three bodies (based on Jacobi reduction), which reduces {\sl all} the integrals and the reduction is {\sl singular} for co--planar motions. 
It turns out that such two apparently different programs  are  both realizable general and, besides, intimately related. Indeed, they have been realized in  \cite{pinzari-th09}. 

\nl
The starting point in \cite{pinzari-th09} was the construction of a set 
of  ``action--angle''  variables, \ie, taking values in $\real^{3n}\times \torus^{3n}$ and  denoted as
$(\L,\G,\Psi,\l,\g,\psi)$, that should extend to the case of  $n\ge 3$ planets Jacobi's reduction of the nodes. Such variables actually already existed: they  had been  considered, in a slightly different form, in the 80's  by F. Boigey \cite{boigey82} for $n=3$ and  A. Deprit \cite{deprit83} for $n\ge 4$. 
Next, Boigey--Deprit variables were rediscovered  in ``planetary''  form\footnote{The   rediscovered variables $(\L,\G,\Psi,\l,\g,\psi)$  are different from the ones in \cite{boigey82}--\cite{deprit83}. They  correspond to be the ``planetary version'' of Deprit's variables.   They are defined only for negative unperturbed energies, so are less general, but turn to be better fitted to the planetary problem, since involve the elliptic elements of the planets.  Also the proof of their symplectic character is different from \cite{boigey82}--\cite{deprit83}: for $n=2$ the $(\L,\G,\Psi,\l,\g,\psi)$  were obtained in
\cite{pinzari-th09} constructively (via generating function, starting with Delaunay variables). This proof was however never published, after realizing the partial coincidence with the variables of \cite{deprit83}. For $n\ge 2$, the proof in \cite{pinzari-th09} is by induction. Part of this proof  was later published in  \cite{chierchiaPi11a} as a ``new proof of Deprit variables''. In \cite{chierchiaPi11a} the relation between the two sets  is also clarified.} by the author (who was strongly motivated by the present application to the Planetary Theorem) during her PhD,  in the first months of 2008. Incidentally, the author would be grateful to anyone who let her know of applications of Boigey--Deprit variables to physical systems with $n\ge 3$ particles (the first not known case, after Jacobi), before \cite{pinzari-th09}. 

\nl
Next, a new set of {\sl regular} variables, named  ``Regular'', ``Planetary'' and ``Symplectic'' -- {\sc rps} -- and denoted with analogue symbols  as  Poincar\'e variables, 
 \beqa{z}
\L&=&(\L_1,\cdots,\L_n)\ ,\quad \l=(\l_1,\cdots,\l_n)\nonumber\\
 z&=&(\eta_1,\cdots,\eta_n, \xi_1,\cdots,\xi_n, p_1,\cdots,p_n, q_1,\cdots,q_n)\ ,\eeqa
 having all the properties conjectured by Arnold for the many--body case was determined, in \cite{pinzari-th09}. Such variables were not discussed by Boigey and Deprit.
They were obtained 
 by applying to the $(\L,\G,\Psi,\l,\g,\psi)$ 's a regularization similar to Poincar\'e's regularization of Delaunay variables. Though being qualitatively 
 similar to Poincar\'e variables,  at contrast with them, {\sc rps} variables  are not ``six per planet'' (the coordinates the $i^{\rm th}$ planet are determined by the variables $\l_i$ and $(\L_j,\eta_j,\xi_j,p_j,q_j)$ with $i\le j\le n$, because of a certain hierarchical structure in their definition, actually inherited by the $(\L,\G,\Psi,\l,\g,\psi)$ 's). Moreover, {\sc rps} variables  are better fitted to  rotation invariance of the problem, since they exhibit a cyclic couple $(p_n,q_n)$ of conjugated variables (integrals of motion). The disappearing of this latter couple of variables from the Hamiltonian implies that the number of degrees of freedom is reduced of one unit (it is $(3n-1)$, one over the minimum) {\sl and}, moreover, the system is let free of two (out of three\footnote{For this reason, following \cite{maligeRL02}, the reduction performed by the {\sc rps} variables is sometimes called ``partial reduction'', at contrast with the ``full reduction'', also discussed in \cite{pinzari-th09}, that reduces the system to  the minimum number $(3n-2)$, of degrees of freedom. Pay attention not to confuse, however, the {\sl regular} ``partial reduction'' performed by {\sc rps} variables with the elementary (but {\sl singular}) reduction that can be obtained reducing the integral ${\rm C}_3$ in Poincar\'e variables. This latter one does not exhibit a cyclic couple and has nothing to do with the aforementioned Arnold's claim in \cite[Ch. 3, \S 5, 5]{arnold63}.}) non--commuting integrals, as Arnold claimed.

 \nl
In place of the  ``Poincar\'e Hamiltonian'' \equ{prop deg}, we consider the  ``{\sc rps} Hamiltonian''
\beq{rps ham}\cH_{\rm rps}=h_{\rm Kep}(\L)+\m f_{\rm rps}(\L,\l,\bar z)\eeq
with \beq{bar z}\bar z=(\eta,\xi,\bar p,\bar q)\ ,\eeq where $\eta=(\eta_1,\cdots,\eta_n)$, $\bar p=(p_1,\cdots,p_{n-1})$ and so on.

\nl
Fixing the value of $(p_n,q_n)$ corresponds to fix one of the $\infty^2$ invariant manifolds that foliate the phase space; letting the other $2(3n-1)$ vary gives a symplectic chart on any of such manifolds. On any of such invariant manifolds, the Birkhoff normal form has been proved to exist (with the properties described at the beginning of the paragraph, but with $(3n-1)$ degrees of freedom, instead of $3n$). Moreover, this normal form satisfies the non--degeneracy condition required by the Fundamental Theorem and the direct proof of Arnold's Planetary Theorem follows.

\nl
It has also  been proved \cite{chierchiaPi11c} that this construction is necessary. Namely that the unreduced system in Poincar\'e variables \equ{prop deg}  would admit a Birkhoff normal form (we remark, despite of the secular resonances), but this normal form would be degenerate {\sl at any order}: the lowest order of it corresponding to the rotational resonance. At the fourth order,  the system would exhibit an identically vanishing torsion (given by the torsion of the partially reduced system, bordered with a row and a column of zeroes) and so on. In particular, no  {\sc kam} theory might   be {\sl directly} applied to the unreduced system \equ{prop deg}. 

\nl We refer to \cite{chierchiaFP13} for
more information on this topic. Other reviews appeared in \cite{fejoz13}, \cite{chierchia13}.

\vskip.2in
\noi
{\large \bf 1.2} 
This paper is concerned with a more detailed study of the normal form constructed in \cite{pinzari-th09}, \cite{chierchiaPi11b}. Before describing it, we anticipate two  applications.


\vskip.1in
{\bf a) 
A ``uniform''  theorem on quasi--periodic motions} The former result of this paper is an improvement of the statements of the Planetary Theorem found in \cite{robutel95} and \cite{pinzari-th09}--\cite{chierchiaPi11b}, in the case of he spatial three--body problem. 
In such papers, a positive measure set of quasi--periodic motions has been obtained, provided eccentricities and the mutual inclination among the planets are suitably small. Moreover, the {\sl Kolmogorov set} (the union of quasi--periodic motions) {\sl depends} strictly on eccentricities and the inclination, in the sense that its density  tends to one as eccentricities and the inclination go to zero.
In fact, the proofs in such papers  are based on the application of the Fundamental Theorem (or improved formulations of it, \cite{chierchiaPi10}), where  this assumption is essential: compare the first inequality in \equ{A5} and the measure of $\cK_{\m,\e}$ below.

\nl
In the case of the {\sl planar} three--body problem this assumption can be relaxed. In Arnold's words:

\nl
{\bf \cite[p. 128]{arnold63}}
{ \sl ``In the case of three bodies {\rm [on a plane]}  we can obtain stronger results (...). It turns out that it is not necessary to require the eccentricities to be small; all that is necessary is that they
should be small enough to exclude the possibility of collision.''}

\nl
And in fact, he stated (we refer to Appendix \ref{ArnoldKAM} for notations)
\begin{theorem}[V. I. Arnold, {\cite[p. 128]{arnold63}}]\label{simpler planar} In the case of the planar three--body problem,   it is possible to find $\m_*>0$, $a_*>0$ such that if
\beq{cond on mu}|\m|<\m_*\eeq 
an invariant set $\cK_{\m}\subset
  \cP_{\e_0}$, with $$ {\meas{\cK_\m}}\ge(1-\m^{a_*}){\meas \cP_{\e_0}}$$
formed by the union of invariant four--dimensional tori, on which the motion is analytically conjugated to  linear Diophantine quasi--periodic motions.
\end{theorem}
He then {\sl conjectured}  the same should hold also for the {\sl spatial} problem:
 \begin{conjecture}[V. I. Arnold, {\cite[p. 129]{arnold63}}]\label{conjecture}
An analogous {\rm[}{\rm to} {\rm Theorem}  {\rm \ref{simpler planar}}{\rm]} theorem is valid for the space three-body problem. In this case, one has to add to condition  \equ{cond on mu} a smallness condition for inclinations.
 \end{conjecture}

 \nl
In \cite{arnold63}, Arnold gave some hints to prove Conjecture \ref{conjecture}. In the 90's M. Robert Herman pointed out a serious gap in such indications. Since then, this stronger case of the Planetary Theorem  remained unproved.

\nl
We shall prove the following
\vskip.1in
 \noi
{\bf Theorem A} {\it
In the spatial three--body problem, there exist numbers $\a_*$, $\m_*$, $\e_*$, $c_*<C_*$ and $\b_*$ such that, if {the numbers $\a$ and  $\m$ (where $\m$ is the masses ratio) verify
$$0<\m<\m_*\ ,\quad 0<\a<\a_*\ ,\quad \m<c_*\log(\a^{-1})^{-4\b_*}$$
in the domain $\cD_\a$ where semi--axes $a_1$, $a_2$, eccentricities $e_1$, $e_2$ and mutual inclination $\iota$ verify 
$$\cD_\a:\quad a_-\le a_1<\a\, a_2\ ,\quad |(e_1,e_2,\iota)|<\e_*$$
a
set $\cK_{\m,\a}\subset  \cD_\a$
may be found,  formed by the union of invariant $5$--dimensional tori, on which the motion is analytically conjugated to  linear Diophantine quasi--periodic motions. The  set $\cK_{\m,\a}$  is of positive Liouville--Lebesgue measure and satisfies, uniformly in $\e$,
$$
\meas \cK_{\m,\a}> \Big(1- C_*({\sqrt[4]\m (\log \a^{-1})^{\b_*}}+\sqrt\a) \Big) \meas \cD_{\a}\ .
$$
The same assertion holds for the planar $(1+n)$--body problem.}}
 \vskip.1in
 \noi
 Note that the thesis of Theorem {\rm A} is a bit weaker  than the one of Theorem \ref{simpler planar}, since, in Theorem {\rm A}  this density is not uniform with respect to the semi--major axes ratio.


\vskip.1in
 {\bf b) 
A ``full'' Nehoro{\v{s}}ev stability theorem} The latter result of the paper is concerned with the  stability for the planetary system. To introduce it, we recall the following fundamental result by N. N. Nehoro{\v{s}}ev
\footnote{A more technical statement of Theorem \ref{Nek thm simplified} is given in Appendix \ref{Nekhorossev}: Compare Theorem \ref{Nek thm}. Recall that other improved statements of Theorem \ref{Nek thm simplified}  have later been found in particular cases: see, for example, \cite{poschel93}, \cite{niedermanB12} and references therein.}, 
mainly motivated by its application to the  Hamiltonian \equ{prop deg}.
 \begin{theorem}[N. N. Nehoro{\v{s}}ev, 1977, {\cite{nehorosev77}, \cite{nehorosev79}}]\label{Nek thm simplified} Let 
 $$H(I,\varphi,p,q)=H_0(I)+\m P(I, \varphi, p, q)\ ,\quad (I, \varphi, p, q)\in\cP
 \subset \real^{n_1}\times \torus^{n_1}\times \real^{2n_2}$$ 
  be of the form of \equ{prop deg}, real--analytic. 
Assume that $H_0(I)$ is ``steep''.
Then, one can find $a$, $b>0$, $C$ and $\m_0$ such that, if
$\m<\m_0$,
any trajectory $t\to \g(t)=(I(t), \varphi(t), p(t), q(t))$ solution of $H$ such that \beq{inclusion simple}(p(t), q(t))\in \P_{(p,q)}\cP\ ,
\quad \forall\ 0\le t\le T_0:= \frac{1}{C\m}e^{\frac{1}{C\m^a}}\eeq
verifies
\beqno|I(t)-I(0)|\le r_0:=\frac{C}{2}\m^b\qquad \forall\ 0\le t\le T_0\ .\eeqno
\end{theorem}

\nl
As for the definition of ``steepness'', we refer to the papers 
 \cite{nehorosev77}, \cite{nehorosev79} and \cite{nehorosev73}. See also \cite{niederman06} for an equivalent definition. We aim to point out that,  despite of the almost 150--pages length of the proof of Theorem \ref{Nek thm simplified} and the complication of notion of steepness, in \cite{nehorosev77} Nehoro{\v{s}}ev 
easily\footnote{The only delicate point in the application of Theorem  \ref{Nek thm simplified} to $\cH_{\rm P}$ consisted in checking assumption \equ{inclusion simple}, that Nehoro{\v{s}}ev accomplished using the conservation of the third component ${\rm C}_3$ of the total angular momentum \equ{C} along the $\cH_{\rm P}$--trajectories. Note that, in the non--degenerate case, \ie, when the variables $(p, q)$ do not appear, this assumption is void.}
applied  Theorem \ref{Nek thm simplified} to the planetary Hamiltonian $\cH_{\rm P}$ in \equ{prop deg} (with $I=\L$, the actions related to the semi--axes, and $(p,q)=\uz$ in \equ{PV}, the secular variables related to eccentricities and inclinations), since the unperturbed term $H_0=h_{\rm Kep}$ is {\sl concave}, a special case of steepness. Nehoro{\v{s}}ev then obtained a spectacular result of stability   for the planetary semi--axes (hence, absence of collisions) over exponentially--long times  {\sl for all initial data} in phase space (see also \cite{niederman96} for a different approach and improved estimates).
Up no now, Nehoro{\v{s}}ev's result is the only  rigorous, global (\ie, valid on the whole phase space, or, possibly, on a very large open subset of it) stability result for the planetary problem. Indeed, there do exist in literature results involving also strong numerical efforts for physical systems (see,\eg,  \cite{sansotteraLG13}, \cite{giorgilliLS09} and references therein) true on Cantor sets (in general, they are obtained via {\sc kam} techniques).

\nl
A physically relevant and widely studied  open problem
is related to the study of the stability of the whole system; \ie, the study of the secular variation  of eccentricities and inclinations of the planets' instantaneous orbits, besides the ones of semi--axes. See, for example, 
 \cite{laskar12} and references therein. Partial rigorous results  in this direction have been obtained in \cite{chierchiaPi11c}, where it has been proved that, if eccentricities and inclinations are initially suitably small, they remain confined with respect to their initial values over {\sl polynomially} long times, up to exclude the so--called\footnote{I. e., resonances of the Keplerial frequencies $\o_{\rm Kep}:=\partial h_{\rm Kep}$.} ``mean--motion resonances''.
More precisely, the following result has been proved.
\begin{theorem}[{\cite{chierchiaPi11c}}]\label{CP11c}
Whatever is the number of planets, for any arbitrarily fixed $s\in \natural$, with $s\ge5$, one can find positive numbers ${\rm C}$, $\underline a_j$, $\ovl a_j$, $\underline\e$, $\overline\e$ with $\underline a_j<\ovl a_j<\underline a_{j+1}$ and $\underline\e<\overline\e$ such that for any $\k>0$,  in the domain where semi--major exes $a_i$, eccentricities $a_i$ and mutual inclinations $\iota_j$ verify
$$\hat\cD_{s,\e}:\quad \underline a_j\le a_j\le \ovl a_j\, \quad \underline\e<\max_{i,j}\{e_i, \iota_j\}<\e<\ovl\e$$
under suitable relations between $\m$ and $\e$, one can find an open set $\hat{\hat\cD}_{s,\m,\e}$  such that, for all the motions starting in $\hat{\hat\cD}_{s,\m,\e}$, the displacement of eccentricities and inclinations with respect to their initial values is bounded by $\k \underline\e$, for all
$$|t|\le \frac{{\rm C}\k}{\m \underline\e^{s}}\ .$$ 
\end{theorem}
 The proof of Theorem \ref{CP11c} again relies with the Birkhoff normal form of the system: the time of stability is related in fact to the remainder of this normal form. No analysis of resonance zones, trapping arguments... is used  for its proof. An undesirable aspect of  Theorem \ref{CP11c}, is that the size of $\hat{\hat\cD}_{\s,\e}$ decreases with with the time of stability.

\nl
In this paper, we prove a stronger result, at least for the planar three--body problem.
 \vskip.1in
 \noi
{\bf Theorem B} {\it
In the planar three--body problem, there exist numbers $\bar a_-$, $\bar\a$, $\bar\e$, $\bar a$, $\bar b$, $\bar c$,  $\bar d$ such that, 
in the domain  
$$\bar\cD_{\e}:\quad \bar a_-\le a_1<\bar\a\, a_2\ ,\quad \underline\e<|(e_1,e_2)|<\e<\bar\e$$ under suitable relations between $\m$ and $\e$, one can find an open set $\bar{\bar\cD}_{\m,\e}\subset \bar\cD_{\e}$, defined by absence of mean--motion resonances up to a suitable order,  such that,  for all the motions with initial datum in $\bar{\bar\cD}_{\m,\e}$, one has
$$|a_i(t)-a_i(0)|,\ |e_i(t)-e_i(0)|\le \bar r:=\max\{\d^{\bar b},\m^{1/12},\ \e\}\quad \forall\ 0\le t\le \bar T=\frac{e^{\frac{1}{\bar\d^{\bar a}}}}{\bar\d}$$
where $\bar\d:=\frac{\m^{\bar d}\e}{\bar c}$.
}

\nl

\vskip.1in
\noi


\vskip.2in
\noi
{\bf\large 1.3} 
Let us  sketch the proofs of Theorems A and B and make some comment.

\nl
The proof of  Theorem {\rm A} is  a remake of an idea by V. I. Arnold in \cite{arnold63}. His proof of Theorem \ref{simpler planar}   relies on the observation  that the planar three--body system\footnote{In ``planar'' Poincar\'e variables $(\L_i,\ul_i,\uh_i,\ux_i)$, $i=1$, $2$.}
$\cH_{\rm pl3b}=h_{\rm Kep}(\L)+\m f_{\rm pl3b}(\L,\ul,\uh,\ux)$ 
enjoys the strong property that secular perturbation $(f_{\rm pl3b})_{\rm av}$ is {\sl integrable}. It has two degrees of freedom (related to the secular variables $(\uh_1,\ux_1)$ and $(\uh_2,\ux_2)$) and two commuting integrals: the third component of the angular momentum \equ{C} (the only one non to vanish, since the problem is planar) and itself. Then the Birkhoff  series of $(f_{\rm pl3b})_{\rm av}$ converges and Arnold can use a {\sc kam} theory (recalled in Appendix \ref{ArnoldKAM}, Theorem \ref{simplifiedFT}) that is less general than the Fundamental Theorem but better fitted to this case. In the proof of Theorem A we use a similar idea. Let us denote as $f_{\rm 3b}$ the function  $f_{\rm rps}$ for the three--body case; $(f_{\rm 3b})_{\rm av}$, its the averaged value. We shall see below that a suitable approximation  $(f_{\rm 3b})_{\rm av}^\ppd$ defined in Eq. \equ{general a-expansion} below, is {\sl integrable}. This fact has been already used, in different settings, in \cite{lidovZ76},  \cite{Zhao-th13} and \cite{palacianSY12}. 
Moreover, the same property of integrability is proved  to hold for the {\sl planar} many--body problem; see below for more details on this assertion.
Then, we apply Arnold's argument,   but working on $(f_{\rm 3b})_{\rm av}^\ppd$, $(f_{\rm pl})_{\rm av}^\ppd$, respectively,  simply suitably modifying Theorem \ref{simplifiedFT}: see  Theorem \ref{thm: ArnoldConjnew}.

\nl
 The proof of Theorem {\rm B} is an application of the  Nehoro{\v{s}}ev's Theorem in the non--degenerate case.
Essentially, it relies on checking ``steepness'' of some integrable truncation of the ``Birkhoff--normalized'' system 
  $${\rm H}_0:=h_{\rm Kep}+\m (f_{\rm bnf})_{\rm av}$$
 in all of its degrees of freedom. Here the difficulty is that, at contrast with the application in \cite{nehorosev77} (where only the concavity of $h_{\rm Kep}$ is exploited), the ``full torsion'' of the system, given by the Hessian of $h_{\rm Kep}$ and the matrix $\b$ of the second--order Birkhoff invariants,  {\sl is not convex}, nor quasi--convex. Its eigenvalues are alternating in sign. 
  Therefore, it is necessary to consider  higher orders of Birkhoff normal form and apply more refined conditions for steepness. It is not clear (and actually an open question) what is the right order of the Birkhoff 
series to be involved for general $n$ and, especially, how steepness can be checked for systems with many degrees of freedom (see 
\cite{schirinziG13} for progresses in this direction).
 For three--degrees of freedom systems Nehoro{\v{s}}ev proved that the ``three--jet condition'' (recalled in Appendix \ref{Nekhorossev})
is ``generic''. But the {\sl planar} three--body problem,  after reducing completely rotations, has three degrees of freedom, so it is not surprising that this problem satisfies three--jet. We do this check in \S \ref{steepness}.

\nl
Before passing to describe technical aspects, we provide a few comments.

\begin{itemize}
\item[--] Theorem B is stated for the planar three--body problem. As previously outlined, the secular problem associated to it is {\sl integrable}: its Birkhoff normal form converges. And in fact this circumstance allowed Arnold to obtain refined results for this case (see \S 1.2): the independence of the {\sl Kolmogorov set} on the eccentricities. One might ask if such independence holds also in the statement of Theorem B. I. e., if the set $\bar{\bar\cD}_{\m,\e}$ may be chosen to be independent of $\e$.  However, with our proof we are not able\footnote{The dependence of $\bar{\bar\cD}_{\m,\e}$ on $\e$ may be read in inequality just before \equ{Kappa} and by the formula \equ{Kappa}, that define this set.} to refine the result in that direction. The reason  is technical:  instead of the (integrable) secular system $\cH_{\rm pl3b}:=h_{\rm Kep}+\m(f_{\rm pl3b})_{\rm av}$ that would be more natural, during the proof we consider a {\sl non integrable} system close\footnote{Compare the system $h_{\rm Kep}+\m(\hat N+\hat N_*)$ in \equ{reduced averaged phsp}.} to it, by performing not only one but {\sl many} steps of averaging with respect to fast (mean motion) frequencies. Therefore, we need to {\sl truncate} the Birkhoff series associated to this closely to integrable system and this is the reason we have the dependence of $\e$. In turn, the exigency of many\footnote{Compare Lemma \ref{aver3b}.}  steps comes succeeding in applying the  theory  developed in \cite{nehorosev77}.
\item[--] In \S \ref{proof of stability} we do more  than we need for Theorem B. We compute the Birkhoff normal form of the {\sl spatial} three--body problem, which is\footnote{In particular,  truncating  this formula to the fourth order we recover the formulae found in \cite{pinzari-th09}--\cite{chierchiaPi11b}.}, which is  
\beqa{normal form of order 6} 
(f_{\rm bnf})_{\rm av}&=&-\frac{\bar m_1\bar m_2}{a_2}-\bar m_1\bar m_2\frac{a_1^2}{4a_2^3}\Big(\big(1+3\frac{t_1}{\L_1}+3\frac{t_2}{\L_2}-3(\frac{1}{\L_1}+\frac{1}{\L_2}) t_3\Big)\nonumber\\
&-&\bar m_1\bar m_2\frac{a_1^2}{4a_2^3}\Big(-\frac{3}{2}\frac{t_1^2}{\L_1^2}+6\frac{t_2^2}{\L_2^2}+\frac{3}{2}\frac{t_3^2}{\L_1^2}+9\frac{t_1t_2}{\L_1\L_2}-12 \frac{t_1t_3}{\L_1^2}-9\frac{t_2t_3}{\L_1\L_2}\nonumber\\
&+&10\frac{t_2^3}{\L_2^3}-\frac{3}{2}\frac{t_3^3}{\L_1^2\L_2}-\frac{9}{2}\frac{t_1^2t_2}{\L_1^2\L_2}-\frac{105}{4}\frac{t_1^2t_3}{\L_1^3}-18\frac{t_2^2t_3}{\L_1\L_2^2}+18\frac{t_1t_2^2}{\L_1\L_2^2}\nonumber\\
&+&\frac{105}{4}\frac{t_1t_3^2}{\L_1^3}+\frac{9}{2}\frac{t_2t_3^2}{\L_1^2\L_2}-36\frac{t_1t_2t_3}{\L_1^2\L_2}\big)\big(1+{\rm O}(\frac{\L_1}{\L_2})\big)+\frac{a_1^2}{a_2^3}{\rm O}(|t|^{7/2})+{\rm O}(\frac{a_1^3}{a_2^4})\Big)
\eeqa
and then we reduce to the planar case setting $t_3=0$. 
However, we are not able to extend Theorem B to the spatial case, since we are not able to check steepness for this case. The three--jet condition might fail at least on manifolds of co--dimension one: see Remark \ref{steepness in space}.
\item[--] Besides the previous case, a possible extension of Theorem B to the general {\sl planar} problem might be helped by the fact that, for this case {\sl we know a good approximation of $(f_{\rm bnf})_{\rm av}$, at any order.} This  result is a corollary of the analysis of \S \ref{An asymptotic formula for the secular perturbation}. See also \S 1.4 below.
\item[--] In our strategy of proofs, the planetary Birkhoff normal form (hence, the system \equ{rps ham} in {\sc rps} variables) plays a central r\^ole. The author is not aware (and would be interesting to know) what kind of results could be obtained (and what would be the relative  difficulty) via Herman--F\'ejoz's   normal form  \cite{fejoz04}. 
\end{itemize}
\vskip.2in
\noi
{\bf\large 1.4} 
The main novelty of this paper (with respect to our previous ones on this subject) is a technical lemma of geometrical nature
that helps in the analysis of the secular perturbing function of  the system \equ{rps ham}. This reflects on the computation of the Birkhoff invariants at higher orders.

\nl
Let us remark, at this respect that, in general, computing the Birkhoff invariants of the planetary problem is a huge work. See, for example the computations of the torsion in \cite{arnold63} ($n=2$, planar), \cite{robutel95} ($n=2$, spatial), \cite{herman09} ($n\ge 2$, planar), \cite{pinzari-th09}--\cite{chierchiaPi11b} ($n\ge 2$, spatial). So, our main progress relies on an improvement of the technique of computation of such invariants, which is particularly desirable if one wants to extend Theorem B to the general problem.

\nl
Let us introduce it briefly, referring to the following section for details.

\nl
Consider the system \equ{rps ham} and, in particular, its secular perturbing function $(f_{\rm rps})_{\rm av}$. Since the indirect\footnote{The former term in \equ{pert**} is of often referred to as ``indirect part''; the latter as ``direct part''. As far as the author knows, this terminology has been introduced by the French school. The vanishing of the average of the indirect part, known Poincar\'e variables, holds also in {\sc rps} variables.} part has zero $\l$--average, $(f_{\rm rps})_{\rm av}$ is given by
 $$(f_{\rm rps})_{\rm av}=-\sum_{1\le i<j\le n}\frac{\bar m_i\bar m_j}{(2\p)^2}\int_{\torus^2}\frac{d\l_id\l_j}{|x^\ppi(\L,\l_i,\bar z)-x^\ppj(\L,\l_j,\bar z)|}\ .$$
De--homogeneizating with respect to $a_j$, we expand each of the terms
 $$(f_{\rm rps}^{(ij)})_{\rm av}:=-\frac{\bar m_i\bar m_j}{(2\p)^2}\int_{\torus^2}\frac{d\l_id\l_j}{|x^\ppi(\L,\l_i,\bar z)-x^\ppj(\L,\l_j,\bar z)|}\ .$$
in powers of the ratio $\frac{a_i}{a_j}$, with $a_j$ fixed:
\beq{aij exp}(f_{\rm rps}^{(ij)})_{\rm av}=(f_{\rm rps}^{(ij)})_{\rm av}^\ppo+(f_{\rm rps}^{(ij)})_{\rm av}^\ppu+(f_{\rm rps}^{(ij)})^\ppd_{\rm av}+\cdots\ .\eeq
Clearly, to this expansion there corresponds an analogue expansion of
\beq{general a-expansion}(f_{\rm rps})_{\rm av}=(f_{\rm rps})_{\rm av}^\ppo+(f_{\rm rps})^\ppu+(f_{\rm rps})^\ppd_{\rm av}+\cdots\ .\eeq
Analogously to what happens for the Poincar\'e Hamiltonian \equ{prop deg}, one has that, in these expansions, the zeroth order terms $(f_{\rm rps}^{(ij)})_{\rm av}^\ppo$ are independent\footnote{They are given by given by$-\frac{\bar m_i\bar m_j}{a_j}$.} of $\bar z$  by well known properties of the two--body potential and  that 
the linear terms $(f_{\rm rps}^{(ij)})_{\rm av}^\ppu$ vanish by Fubini's and Newton equation\footnote{I. e., by the vanishing of $$\frac{1}{2\p}\int_{\torus}\frac{x^\ppj(\L,\l_j,\bar z)}{|x^\ppj(\L,\l_j,\bar z)|^3}d\l_j=\frac{1}{{\rm T}_j}\int_{0}^{{\rm T}_j}\frac{d}{dt}y^\ppj(\L,\o_j t,\bar z)dt$$
with some ${\rm T}_j$ and $\o_j=\frac{2\p}{{\rm T}_j}$.
}. The lowest order information on $(f_{\rm rps})_{\rm av}$ is then given by the second--order terms $(f_{\rm rps})^\ppd_{\rm av}$.

\nl
By \cite{pinzari-th09}--\cite{chierchiaPi11b} $(f_{\rm rps})_{\rm av}^\ppd$ may be splitted into a sum
 \beq{sum}(f_{\rm rps})_{\rm av}^\ppd=(f_{\rm pl})_{\rm av}^\ppd+(f_{\rm vert})_{\rm av}^\ppd\eeq
 of a ``planar''  and\footnote{We follow the terminology in \cite{fejoz04}.} a ``vertical'' part, where $(f_{\rm pl})^\ppd_{\rm av}$  corresponds to the term that we would have for the problem in the plane, while $(f_{\rm vert})^\ppd_{\rm av}$  vanishes for $(\bar p,\bar q)=0$ and is even in $(\bar p,\bar q)$. In \S \ref{An asymptotic formula for the secular perturbation} we prove that $(f_{\rm pl})^\ppd_{\rm av}$, $(f_{\rm vert})^\ppd_{\rm av}$  are given by, respectively,
 \beqa{averaged planar}
&&(f_{\rm pl})^\ppd_{\rm av}=-\frac{1}{4}\sum_{1\le i<j\le n}\bar m_i\bar m_j\frac{a_i^2}{a_j^3}\frac{\frac{1}{2\p}\int_\torus\frac{d\zeta}{1-e_j\cos\zeta}}{(1-\frac{\eta_j^2+\xi_j^2}{2\L_j})^2}(1+\frac{3}{2} e_i^2)\nonumber\\
&&(f_{\rm vert})^\ppd_{\rm av}=+\frac{3}{4}\sum_{1\le i<j\le n}\bar m_i\bar m_j\frac{a_i^2}{a_j^3}\frac{\frac{1}{2\p}\int_\torus\frac{d\zeta}{1-e_j\cos\zeta}}{(1-\frac{\eta_j^2+\xi_j^2}{2\L_j})^2}\frac{1}{2\p}\int_{\torus}(\hat x^\ppi\cdot \hat{\rm C}^\ppj)^2d\l_i\ ,\eeqa
where $e_i$'s are the eccentricities,  expressed in terms of  $\L_i$ and $\frac{\eta_i^2+\xi_i^2}{2}$; $\hat{\rm C}^\ppj$ are the planets' normalized angular momenta $\frac{{\rm C}^\ppj}{|{\rm C}^\ppj|}$ and $\hat x^\ppi:=\frac{x^\ppi(\L,\l_i, z)}{a_i}$. 

\nl
The author is not aware if the formulae  \equ{averaged planar} had been already noticed before (they hold also in the case of the Poincar\'e system \equ{prop deg}). Such formulae are the thesis of Proposition  \ref{asymptotic formula}, that we prove using a new set of symplectic variables, defined in \equ{def}, and tools of normal form theory. The variables \equ{def} in a sense resemble the  well known Adoyer--Deprit variables of the rigid body, with the difference that have six degrees of freedom instead of three. Also the thesis of Proposition  \ref{asymptotic formula} resembles certain formulae for the rigid body, as outlined in Remark \ref{Andoyer-Deprit}.
 
 \nl
In particular,  inspecting  \equ{averaged planar}, it is to be remarked that $(f_{\rm pl})_{\rm av}^\ppd$ not only is {\sl integrable}, but {\sl is in Birkhoff normal form}. This fact implies the validity of Theorem A for the planar general problem and, especially, is of great help in the computation of its Birkhoff invariants {\sl at any order}.

\nl
Secondly, formulae \equ{averaged planar} imply that, in the three--body case ($n=2$), $(f_{\rm 3b})_{\rm av}^\ppd:=(f_{\rm rps})_{\rm av}^\ppd|_{n=2}$ is independent of  the argument of $(\eta_2,\xi_2)$, therefeore, it is {\sl integrable} (compare \cite{lidovZ76}   for an analogue assertion  in a different setting and  \cite{Zhao-th13} and \cite{palacianSY12} for applications). More in general, for $n\ge 2$, $(f_{\rm rps})_{\rm av}^\ppd$ is independent on the argument of $(\eta_n,\xi_n)$.
But while, for this general case, the expression of  $(f_{\rm vert})_{\rm av}^\ppd$ in terms of {\sc rps} variables is complicated, due to the factors $(\hat x^\ppi\cdot\hat{\rm C}^\ppj)^2$, it is not so for three bodies, where there is only one of such factors ($i=1$, $j=2$). The  aspect of the corresponding vertical term is nice
\beq{averaged spatial3b}(f_{\rm 3bvert})_{\rm av}^\ppd=\frac{3}{4}\bar m_1\bar m_2\frac{a_1^2}{a_2^3}\frac{\frac{1}{2\p}\int_\torus\frac{d\zeta}{1-e_2\cos\zeta}}{(1-\frac{\eta_2^2+\xi_2^2}{2\L_2})^2}
\Big((1+\frac{3}{2} e_1^2)({\rm i} v v^\star)+\frac{5}{2}\big(( u_1^\star)^2 v^2+(v^\star)^2  u_1^2\big){ \bar e_1^2}\Big)\bar{\eufm s}^2\eeq
where $u_i$, $u_i^*$ are the Birkhoff variables associated to $(\eta_i,\xi_i)$; $(v, v^*)$ to $(p_1,q_1)$,
$\bar e_1$
and $\bar{\eufm s}$ are suitable functions in normal form. Since the first non--normal terms in this formula appear from the fourth order on, the computation of the sixth orders Birkhoff invariants for the three--body case is quickly done: it takes less than two pages (see \S \ref{Some dynamical application}) and gives \equ{normal form of order 6}.

\vskip.2in
\noi
{\bf Acknowledgments} I wish to express my regards and deep admiration to Professor Alain Chenciner in occasion of his seventieth birthday and  thank warmly Alain Albouy and Alexey Borisov for inviting me to present this paper in that important circumstance.

\nl
Moreover, I would like to thank Massimiliano Berti, Jacques F\'ejoz, Boris Khesin, Sergei Kuksin, Vittorio Coti Zelati and Edi  Zehnder for honoring me, since my PhD, of their comments on the planetary problem and their encouragement. 

\nl
Doubly thanks to Jacques F\'ejoz, refer\'ee of my PhD thesis,  for mentioning in \cite{fejoz13}  my contribution to the proof of Arnold's Theorem, especially, my rediscovery of Deprit's reduction. I remember with much pleasure the long, relaxing and enlightening discussions on the planetary problem since we met for the first time in early Spring of 2008. I feel deeply indebted with him, since then.

\nl
Thanks to Alessandra Celletti, Giancarlo Benettin, Luca Biasco, Francesco Fass\`o, Massimiliano Guzzo and  Lei Zhao, who helped me with precious bibliographic advices.

\nl
To the anonymous Referee, for his  thoughtful advices, without which this paper would not have this form. In particular, for advertising me on the possibility of more global results than I stated in Theorem B, in the case of the planar three--body problem.
  
 \nl
To my early adviser Luigi Chierchia, without whom nothing of what I did since my  PhD would exist.

 \nl To my husband, without whom nothing would have meaning. 
\section{An asymptotic formula for the secular perturbation}\label{An asymptotic formula for the secular perturbation}
Let, for fixed $1\le i<j\le n$,
\beqno
f_{ij}(\L,\bar z):=\frac{1}{(2\p)^2}\int_{\torus^2}\frac{d\l_id\l_j}{|x^\ppi(\L,\l_i,\bar z)-x^\ppj(\L,\l_j,\bar z)|}
\eeqno
so as to write
 \beqno
(f_{\rm rps})_{\rm av}(\L,\bar z)=-\sum_{1\le i<j\le n} \bar m_i\bar m_j f_{ij}(\L,\bar z)\ .
\eeqno
Here\footnote{\label{foot24}Actually, the map \equ{rps map} depends on $z$, rather than $\bar z$. However, by the independence of the Hamiltonian \equ{rps ham} of $(p_n,q_n)$, we may  arbitrarily fix such couple of variables to some value, \eg, $(0,0)$. Abusively, {\sl just in \equ{double average new}} and similar formulae below, we denote again as $(\L,\l,\bar z)\to (y(\L,\l,\bar z), x(\L,\l,\bar z))$ the map  $\phi_{\rm rps}^{-1}|_{(p_n,q_n)=(0,0)}$. } $(\L,\l_i,\bar z)\to x^\ppi(\L,\l_i,\bar z)$ denotes the $x^\ppi$--projection of the map 
\beq{rps map}
\phi_{\rm rps}^{-1}:\quad (\L,\l,z)\to (y, x)\in \real^{3n}\times \real^{3n}
 \eeq
Consider the formal expansions
\beq{alpha expansion}f_{ij}=f_{ij}^\ppo+f_{ij}^\ppd+\cdots\eeq
  in powers of the semi--major axes ratio $\a
_{ij}
:=a_i/a_j$, with $a_j$ fixed.
Here,
$$f_{ij}^{(k)}:=\frac{1}{k!}\frac{d^k}{d\varepsilon^k}\Big[\frac{1}{(2\p)^2}\int_{\torus^2}\frac{d\l_id\l_j}{|\varepsilon x^\ppi(\L,\l_i,\bar z)- x^\ppj(\L,\l_j,\bar z)|}\Big]
_{\varepsilon=0}\ .$$

\nl
In particular, we focus on the  second--order term of this expansion,  given by
\beqa{double average new}
f_{ij}^{(2)}&&=\frac{1}{(2\p)^2}\int_{\torus^2}d\l_id\l_j\nonumber\\
&&\quad\frac{3(x^\ppi(\L,\l_i,\bar z)\cdot x^\ppj(\L,\l_j,\bar z))^2-|x^\ppi(\L,\l_i,\bar z)|^2|x^\ppj(\L,\l_j,\bar z)|^2}{2|x^\ppj(\L,\l_j,\bar z)|^5}\ .\nonumber\\
\eeqa
Note that $(f_{\rm rps})^\ppd_{\rm av}$ in \equ{general a-expansion} corresponds to \beq{a-exp neew}(f_{\rm rps})^\ppd_{\rm av}=-\sum_{1\le i<j\le n}\bar m_i\bar m_j f^\ppd_{ij}\ .\eeq
Let ${\rm C}^\ppi(\L,\bar z):=x^\ppi(\L,\l_i,\bar z)\times y^\ppi(\L,\l_i,\bar z)$ (by definition of the map \equ{rps map}, ${\rm C}^\ppi(\L,\bar z)$ is independent of $\l_i$). We have the following identity
\begin{proposition}\label{asymptotic formula}
\beq{2aver**}f_{ij}^{(2)}=-\frac{M_jm_j^2}{4}\frac{\dst\frac{1}{2\p}\int_{\torus}\big(3({\rm C}^\ppj\cdot x^\ppi)^2-|x^\ppi|^2|{\rm C}^\ppj|^2\big)d\l_i}{\dst|{\rm C}^\ppj|^4}(\frac{1}{2\p}\int_{\torus}\frac{d\l_j}{\dst |x^\ppj|^2})\eeq
\end{proposition}
Note that Eqs. \equ{a-exp neew}, \equ{2aver**} and the  formulae of $|{\rm C^\ppj}|$, $|x^\ppj|$ in terms of {\sc rps} variables (see \cite{pinzari-th09}, \cite{chierchiaPi11b} and eventually Appendix \ref{formula for Theta}) imply \equ{sum}--\equ{averaged planar}.

\vskip.1in
\noi
We first discuss
\subsection{The three--body case}
Let 
\beqano
{\rm P}^\ppd&&:=\frac{1}{(2\p)^2}\int_{\torus^2}d\ul_1d\ul_2\nonumber\\
&&\quad\frac{3(x^\ppu(\L,\ul_1,\uz)\cdot x^\ppd(\L,\ul_2,\uz))^2-|x^\ppu(\L,\ul_1,\uz)|^2|x^\ppj(\L,\ul_2,\uz)|^2}{2|x^\ppd(\L,\ul_2,\uz)|^5}\ .
\eeqano
where, for $i=1$, $2$, $$(\L_1,\L_2, \ul_i,\uz)\in \cA^2\times\torus^1\times B^8\to(y^\ppi(\L_1,\L_2, \ul_i,\uz), x^\ppi(\L_1,\L_2, \ul_i,\uz))\in \real^3\times \real^3$$ are two mappings  such that  
 \begin{itemize}
\item[{\bf (A)}] The map $(\L_1,\L_2, \ul_2,\uz)\to(y^\ppd(\L_1,\L_2, \ul_2,\uz), x^\ppd(\L_1,\L_2, \ul_2,\uz))$  solves the two--body problem {\rm ODE}
\beq{Newton law}\partial_{t}y^\ppd(\L_2, \o^\ppd_{\rm Kep}t,\uz)=-m_2M_2\frac{x^\ppd(\L_2, \o^\ppd_{\rm Kep}t, \uz)}{|x^\ppd(\L_2, \o^\ppd_{\rm Kep}t, \uz)|^3}\eeq
where $\o^\ppd_{\rm Kep}=\o^\ppd_{\rm Kep}(\L_2)=\frac{M_2^2m_2^3}{\L_2^3}$;
\item[{\bf (B)}]  The map
\beq{map**}\bar\phi:\qquad (\L_1,\L_2, \ul_1,\ul_2,\uz)\to (y^\ppu, y^\ppd, x^\ppu, x^\ppd)\eeq
is  symplectomorphism of $\cA^2\times\torus^2\times B^8$ into $\real^{12}$ (where $\cA^2\subset\real^2$, $B^8\subset\real^8$ sre open and connected).
\end{itemize}
\begin{proposition}\label{NF3BP}
Under assumptions  {\bf (A)} and {\bf (B)}, the following identity holds 
\beq{Pav}
{\rm P}^\ppd=-\frac{M_2m_2^2}{4}\frac{\frac{1}{2\p}\int_{\torus}\big(3({\rm C}^\ppd\cdot x^\ppu)^2-|x^\ppu|^2|{\rm C}^\ppd|^2\big)d\ul_1}{|{\rm C}^\ppd|^4}(\frac{1}{2\p}\int_{\torus}\frac{d\ul_2}{|x^\ppd|^2})\eeq
where ${\rm C}^\ppd(\L_1,\L_2, \uz):=x^\ppd(\L_1,\L_2,\ul_2,\uz)\times y^\ppd(\L_1,\L_2,\ul_2,\uz)$.\end{proposition}
\begin{remark}\label{Andoyer-Deprit}\rm

\item[--]
Note that, in the case $n=2$, the map \equ{rps map} satisfies assumptions (A) and (B), hence Proposition \ref{NF3BP} is just Proposition \ref{asymptotic formula} in this particular case.
\item[--] We shall prove more than \equ{Pav}: letting ${\rm P}^\ppu(\L,\ul_1,\uz)$ as in \equ{P1} below, then ${\rm P}^\ppu$ satisfies an analogue identity as in \equ{Pav}, but neglecting the first average $\frac{1}{2\p}\int_{\torus}d\ul_1$.
\item[--] The formula \equ{Pav} resembles  the expression of the averaged quartic term in the spin--orbit problem, using Andoyer--Deprit coordinates: see \cite[Eq. (24)]{benettinGM08}, in turn based on the expansions in \cite[\S 12]{chierchiaG94}.
\end{remark}

\nl
In the next sections, we  prove Proposition \ref{NF3BP}. Next (in \S \ref{proof of 2.1}), we discuss the general case.
\subsection{A six--degrees of freedom set of symplectic variables}\label{sympl var}
The proof of Proposition \ref{NF3BP} is based on the use of a ``ad hoc'' variables for the three--body problem. Let us introduce them.

\vskip.1in
\noi
 Let $(k^\ppu, k^\ppd, k^\ppt)$ be a prefixed orthonormal frame in $\real^3$ and let
 $$(y^\ppu, y^\ppd, x^\ppu, x^\ppd)\in (\real^3)^4\ ,\quad (y^\ppi, x^\ppi)=(y^\ppi_1, y^\ppi_2, y^\ppi_3, x^\ppi_1, x^\ppi_2, x^\ppi_3)$$
be a system of ``Cartesian coordinates'' in the configuration space $\real^3$, with respect to $(k^\ppu, k^\ppd, k^\ppt)$.

 \nl
Denote as
 $${\rm C}^\ppi:=x^\ppi\times y^\ppi$$
 (with ``$\times$'' denoting skew product) the $i^{\rm th}$ angular momentum, and let ${\rm C}:={\rm C}^\ppu+{\rm C}^\ppd$ the total angular momentum. For $u,v\in\real^3$ lying in the plane orthogonal to a vector $w$,  let $\a_w(u,v)$ denote the positively oriented angle (mod $2\p$) between $u$ and $v$  (orientation follows  the ``right hand rule''). Define the  ``nodes''
$$\n_1:=k^\ppt\times {\rm C}\ ,\quad \n_2:={\rm C}\times x^\ppu\ ,\quad \n_3:=x^\ppu\times {\rm C}^\ppd\ .$$
Let $\cP^{12}_\star$ denote the subset of $(\real^3)^4$ where ${\rm C}$, ${\rm C}_2$, $x^\ppu$, $x^\ppd$, $\n_1$, $\n_2$ and $\n_3$ simultaneously do not vanish. On $\cP^{12}_\star$ define a map 
 \beqno
 \phi^{-1}:\quad  (y^\ppu, y^\ppd, x^\ppu, x^\ppd)\to  ({\rm C}_3, {\rm G}, {\rm R}_1, \Theta, {\rm R}_2, \Phi_2,\zeta, {\eufm g}, {\rm r}_1, \vartheta, {\rm r}_2, \varphi_2)\eeqno
 via the following formulae
\beqa{def}
\phi^{-1}:\qquad \arr{
 {\rm C}_3:={\rm C}\cdot k^\ppt\\
 {\rm G}:=|{\rm C}|\\
  {\rm R}_1:=\frac{y^\ppu\cdot x^\ppu}{|x^\ppu|}\\
 \Theta:=\frac{{\rm C}^\ppd\cdot x^\ppu}{|x^\ppu|}\\
   {\rm R}_2:=\frac{y^\ppd\cdot x^\ppd}{|x^\ppd|}\\
\Phi_2:=|{\rm C}^\ppd|
 }\qquad \arr{
 \dst \zeta:=\a_{k^\ppt}(k^\ppu, \n_1)\\
 \dst{\eufm g}:=\a_{\rm C}(\n_1, \n_2)\\
 \dst {\rm r}_1:=|x^\ppu|\\
 \dst\vartheta:=\a_{x^\ppu}(\n_2, \n_3)\\
 \dst {\rm r}_2:=|x^\ppd|\\
\varphi_2:=\a_{{\rm C}_2}(\n_3, x^\ppd) }
\eeqa
 \nl
\begin{proposition}\label{good variables}     The map $\phi^{-1}$   in \equ{def} is invertible on  $\cP^{12}_*$ and  preserves the standard Liouville 1--form $\l=\sum_{i=1}^6{\rm P}_id{\rm Q}_i$.
\end{proposition}
\vskip.1in
\noi
 We denote as 
\beqano\cR_{1}(i)=\left(
\begin{array}{cccc} 
1&0&0\\
0&\cos i&-\sin i\\
0&\sin i&\cos i
\end{array}
 \right)\ ,\qquad \cR_{3}(\th)=\left(
\begin{array}{ccc} 
\cos\th&-\sin\th&0\\
\sin\th&\cos\th&0\\
0&0&1
\end{array}
 \right)\eeqano
 The invertibility is proven by exhibiting the inverse $\phi$.  Indeed, the definitions in \equ{def} and
elementary geometric considerations easily imply the following
 \begin{lemma}\label{lem: good change}
On $\phi^{-1}(\cP^{12}_*)$, the  inverse map of $\phi^{-1}$ in in \equ{def},
has the following analytical expression:
\beqa{map}
\phi:\ \arr{\dst x^\ppu={\cal R}_3(\zeta){\cal R}_1(i){\cal R}_3({\eufm g}){\cal R}_1(i_1)\left(
\begin{array}{lrrr}
0\\
0\\
{\rm r}_1\\
\end{array}
\right)\\ \\
\dst y^\ppu:=\frac{{\rm R}_1}{{\rm r}_1}x^\ppu+\frac{1}{{\rm r}_1^2}{\rm C}^\ppu\times x^\ppu
\\ \\
\dst x^\ppd=
{\cal R}_3(\zeta){\cal R}_1(i){\cal R}_3({\eufm g}){\cal R}_1(i_1)
{\cal R}_3(\vartheta){\cal R}_1(i_2)
\left(
\begin{array}{ccc}
{\rm r}_2\cos\varphi_2\\
{\rm r}_2\sin\varphi_2\\
0
\end{array}
\right)\\ \\
\dst y^\ppd={\cal R}_3(\zeta){\cal R}_1(i){\cal R}_3({\eufm g}){\cal R}_1(i_1){\cal R}_3(\vartheta){\cal R}_1(i_2)\left(
\begin{array}{ccc}
{\rm R}_2\cos\varphi_2-\frac{\Phi_2}{{\rm r}_2}\sin\varphi_2\\
{\rm R}_2\sin\varphi_2+\frac{\Phi_2}{{\rm r}_2}\cos\varphi_2\\
0
\end{array}
\right)
}
\eeqa
where, if $i$, $i_1$, $i_2\in (0,\pi)$  are  defined by
\beq{incli}\cos i=\frac{{\rm C}_{3}}{\rm G}\ ,\quad \cos i_{1}=\frac{\Theta}{{\rm G}}\ ,\quad \cos i_{2}=\frac{\Theta}{\Phi_2}
\eeq
and $ {\rm C}$, 
 $ {\rm C}^\ppd$ by
\beqa{C2}
 {\rm C}&:=&{\cal R}_3(\zeta){\cal R}_1(i)\left(
\begin{array}{ccc}
0\\
0\\
{\rm G}
\end{array}
\right)\nonumber\\
 {\rm C}^\ppd&:=&{\cal R}_3(\zeta){\cal R}_1(i){\cal R}_3({\eufm g}){\cal R}_1(i_1){\cal R}_3(\vartheta){\cal R}_1(i_2)\left(
\begin{array}{ccc}
0\\
0\\
\Phi_2
\end{array}
\right)
\eeqa
then
\beq{C1}{\rm C}^\ppu:={\rm C}-{\rm C}^\ppd\ .\eeq
\end{lemma}
To prove symplecticity we shall use the following easy
\begin{lemma}[\cite{chierchiaPi11a}]\label{lemma on reduction}
Let
$$x=\cR_3(\theta)\cR_1(i)\bar x\ ,\quad y=\cR_3(\theta)\cR_1(i)\bar y\ ,\quad {\rm C}:=x\times y\ ,\quad \bar{\rm C}:=\bar x\times \bar y\ ,$$
with $x,\bar x, y, \bar y\in\real^3$ . Then,
$$y\cdot dx={\rm C}\cdot k^{(3)}d\theta+\bar{\rm C}\cdot k^{(1)} di+\bar y \cdot d\bar x\ .$$
\end{lemma}
\proof {\bf of Proposition} \ref{good variables}.
Let us preliminarly verify that, if ${\rm C}^\ppi$ are as in \equ{C2}--\equ{C1}, and $y^\ppi$, $x^\ppi$ as in \equ{map}, then as expected,
\beq{CC2}x^\ppi\times y^\ppi={\rm C}^\ppi\ .\eeq
Indeed, for $i=2$, this identity is follows trivially from the definitions. To check that it holds also for $i=1$, one can do as follows: firstly, to check that $x^\ppu\cdot {\rm C}^\ppu=0$. This  is an elementary consequence of  \equ{map} and, in particular, of  \equ{incli}. Next, using  the rule of the double skew product, one has
\beqano
x^\ppu\times y^\ppu&=&x^\ppu\times\big(\frac{{\rm R}_1}{{\rm r}_1} x^\ppu+\frac{1}{{\rm r}_1^2}{\rm C}^\ppu\times x^\ppu\big)\nonumber\\
&=&0+\frac{1}{{\rm r}_1^2}\big({\rm r}_1^2\,{\rm C}^\ppu-(x^\ppu\cdot {\rm C}^\ppu)\, x^\ppu\big)={\rm C}^\ppu\ .
\eeqano
 Define now
\beqano
\bar{\rm C}^\ppu&:=&{\cal R}_1(-i){\cal R}_3(-\zeta){\rm C}^\ppu
\nonumber\\
\bar{\rm C}^\ppd&:=&{\cal R}_3({\eufm g}){\cal R}_1(i_1)
{\cal R}_3(\vartheta){\cal R}_1(i_2)\left(
\begin{array}{ccc}
0\\
0\\
\Phi_2
\end{array}
\right)\nonumber\\
\bar{\bar{\rm C}}^\ppu&:=&{\cal R}_1(-i_1){\cal R}_3(-{\eufm g})\left(
\begin{array}{ccc}
0\\
0\\
{\rm G}
\end{array}
\right)-{\cal R}_3(\vartheta){\cal R}_1(i_2)\left(
\begin{array}{ccc}
0\\
0\\
\Phi_2
\end{array}
\right)\nonumber\\
\bar{\bar{\rm C}}^\ppd&:=&
{\cal R}_3(\vartheta){\cal R}_1(i_2)\left(
\begin{array}{ccc}
0\\
0\\
\Phi_2
\end{array}
\right)\nonumber\\
\bar{\bar{\bar{\rm C}}}^\ppd&:=&
\left(
\begin{array}{ccc}
0\\
0\\
\Phi_2
\end{array}
\right)\nonumber\\
\eeqano
and
\beqano
\bar{\bar y}^\ppu&:=&\left(
\begin{array}{ccc}
0\\
0\\
{\rm R}_1
\end{array}
\right)+\frac{1}{{\rm r}_1^2}\bar{\bar{\rm C}}_1\times\left(
\begin{array}{ccc}
0\\
0\\
{\rm r}_1
\end{array}
\right)\ ,\quad 
\bar{\bar x}^\ppu:=\left(
\begin{array}{ccc}
0\\
0\\
{\rm r}_1
\end{array}
\right)\nonumber\\
\bar{\bar{\bar x}}^\ppd&:=&\left(
\begin{array}{ccc}
{\rm r}_2\cos\varphi_2\\
{\rm r}_2\sin\varphi_2\\
0
\end{array}
\right)\ ,\qquad \bar{\bar{\bar y}}^\ppd:=\left(
\begin{array}{ccc}
{\rm R}_2\cos\varphi_2-\frac{\Phi_2}{{\rm r}_2}\sin\varphi_2\\
{\rm R}_2\sin\varphi_2+\frac{\Phi_2}{{\rm r}_2}\cos\varphi_2\\
0
\end{array}
\right)
\eeqano
so as to write
\beqano
&&y^\ppu={\cal R}_3(\zeta){\cal R}_1(i){\cal R}_3({\eufm g}){\cal R}_1(i_1)\bar{\bar y}^\ppu\nonumber\\
&& x^\ppu={\cal R}_3(\zeta){\cal R}_1(i){\cal R}_3({\eufm g}){\cal R}_1(i_1)\bar{\bar x}^\ppu\ .\nonumber\\
&&x^\ppd=
{\cal R}_3(\zeta){\cal R}_1(i){\cal R}_3({\eufm g}){\cal R}_1(i_1)
{\cal R}_3(\vartheta){\cal R}_1(i_2)
\bar{\bar {\bar x}}^\ppd\nonumber\\
&&y^\ppd={\cal R}_3(\zeta){\cal R}_1(i){\cal R}_3({\eufm g}){\cal R}_1(i_1)
{\cal R}_3(\vartheta){\cal R}_1(i_2)\bar{\bar {\bar y}}^\ppd\nonumber\\
\eeqano
Applying repeatedly Lemma \ref{lemma on reduction}, Eq. \equ{CC2} and the rule
$${\cal R}x\times {\cal R}y={\cal R}(x\times y)\qquad {\rm for\ all}\quad  {\cal R}\in {\rm SO(3)},\ x, \ y\in \real^3$$
gives
\beqano
y^\ppu\cdot dx^\ppu&=&{\rm C}^\ppu\cdot k^\ppt d\zeta+\bar{\rm C}^\ppu\cdot k^\ppu d i+\bar{\rm C}^\ppu\cdot k^\ppt d{\eufm g}+\bar{\bar{\rm C}}^\ppu\cdot k^\ppu di_1+{\rm R}_1 d{\rm r}_1\nonumber\\
y^\ppd\cdot dx^\ppd&=&{\rm C}^\ppd\cdot k^\ppt d\zeta+\bar{\rm C}^\ppd\cdot k^\ppu d i+\bar{\rm C}^\ppd\cdot k^\ppt d{\eufm g}+\bar{\bar{\rm C}}^\ppd\cdot k^\ppu di_1+\bar{\bar{\rm C}}^\ppd\cdot k^\ppt d{\vartheta}\nonumber\\
&+&\bar{\bar{\bar{\rm C}}}^\ppd\cdot k^\ppu di_2+{\rm R}_2 d{\rm r}_2+\Phi_2 d\varphi_2
\eeqano
Taking the sum of the two equations and recognizing that, if
$$e^\ppi:={\cal R}_3(\zeta){\cal R}_1(i)k^\ppi\ ,\quad f^\ppi:={\cal R}_3(\zeta){\cal R}_1(i) {\cal R}_3({\eufm g}){\cal R}_1(i_1)k^\ppi$$
then
\beqano
&&({\rm C}^\ppu+{\rm C}^\ppd)\cdot k^\ppt={\rm C}\cdot k^\ppt={\rm G}\cos i={\rm C}_3\nonumber\\
&&(\bar{\rm C}^\ppu+\bar{\rm C}^\ppd)\cdot k^\ppu={\rm C}\cdot e^\ppu=0\nonumber\\
&&(\bar{\rm C}^\ppu+\bar{\rm C}^\ppd)\cdot k^\ppt={\rm C}\cdot e^\ppt={\rm G}\nonumber\\
&&(\bar{\bar{\rm C}}^\ppu+\bar{\bar{\rm C}}^\ppd)\cdot k^\ppu={\rm C}\cdot f^\ppu=({\rm G}k^\ppt)\cdot({\cR}_3({\eufm g})k^\ppu)=0\nonumber\\
&&\bar{\bar{\rm C}}^\ppd\cdot k^\ppt=\Phi_2\cos i_2=\Theta\nonumber\\
&&\bar{\bar{\bar{\rm C}}}^\ppd\cdot k^\ppu=0
\eeqano
we have the thesis:
$$y^\ppu\cdot dx^\ppu+y^\ppd\cdot dx^\ppd={\rm C}_3 d\zeta+{\rm G}d{\eufm g}+\Theta d\vartheta+{\rm R}_1 d{\rm r}_1+{\rm R}_2 d{\rm r}_2+\Phi_2 d\varphi_2\ .\qedeq$$

\subsection{Two--steps averaging for properly--degenerate systems}\label{twost}
In this section we discuss a unicity argument for normal forms of degenerate systems.
\vskip.1in
\noi
Consider a real--analytic and properly--degenerate Hamiltonian
$${\rm H}(I,\varphi, u,v)={\rm H}_0(I)+\a{\rm P}(I,\varphi, u,v)\ ,\qquad 0<\a< 1$$
defined on some phase $(n+m)$--dimensional phase space of the form $V\times\torus^{n_1}\times B^{2n_2}$, where $V$ is an open, connected set of $\real^{n_1}$.
Perturbation theory (\eg, \cite{arnold63}, \cite{nehorosev77}, \cite{poschel93}, \cite{biascoCV03}, \cite{chierchiaPi10}) tells us that, under suitable assumptions of non resonance of the unperturbed frequency map $\o:=\partial_I {\rm H}_0$ and of smallness of the perturbation $\a{\rm P}$, the system may be conjugated, at least formally, to a new system 
\beq{Hp}{\rm H}_p(I,\varphi, u,v)={\rm H}_0(I)+(\a\bar{\rm P}_1(I, u,v)+\cdots+\a^p\bar{\rm P}_p)+\a^{p+1}{\rm P}_{p+1}\ ,\qquad ({\rm P}_{1}\equiv{\rm P})\eeq
where the term inside parentheses (``p--step normal form'') is of degree $p$ and is {\sl independent of $\varphi$}. Quantitative versions of this fact are well known in the literature since \cite{arnold63} and have been more and more refining themselves  (depending on needs) both in the non--degenerate \cite{poschel93}, \cite{chierchiaG94} and degenerate case \cite{arnold63}, \cite{biascoCV03}, \cite{nehorosev77}, \cite{niederman96}. Moreover, we know that, when the system in non--degenerate, \ie, the variables $(u, v)$ do not appear, the {\sl p--step normal form} is uniquely determined (though the change of variables realizing it may be not). In general, when the system is degenerate, uniqueness does not hold.
However, the following lemma is easily proved.
 \begin{lemma}\label{NFT}
Let\footnote{We assume $n_1=1$ to avoid complications due to resonances of the frequency--map. This is enough for the purposes of the paper. Analogue statements for the case $n_1\ge 1$ may be available. } $n_1=1$ and ${\rm H}$
be a properly--degenerate system, such that
\beq{zero average}{\rm P}_{\rm av}:=\frac{1}{2\p}\int_{\torus^n}{\rm P}(I,\varphi;u,v)d\varphi\equiv0\ .\eeq
Then, the two--step normal form
\beqno
\tilde {\rm H}(\tilde I,\tilde\varphi;\tilde u,\tilde v)={\rm H}_0(\tilde I)+(\a\bar{\rm P}_1(\tilde I;\tilde u,\tilde v)+\a^2\bar{\rm P}_2(\tilde I;\tilde u,\tilde v))+{\rm O}(\a^3)\eeqno
is uniquely determined, up to real--analytic and symplectic changes $(\tilde I,\tilde\varphi;\tilde u,\tilde v)\in \tilde{V}\times\torus^{n_1}\times\tilde B^{2n_2}\to (I,\varphi;u,v)\in V^n\times\torus^n\times B^{2n_2}$, $\a$--close to the identity. 
\end{lemma}
\proof 
Let $p\ge 0$. Assuming to have reached the  form in \equ{Hp} (with the term inside parentheses identically vanishing for $p=0$), the $(p+1)^{\rm th}$ Hamiltonian ${\rm H}_{p+1}$ is obtained applying to ${\rm H}_p$ any transformation in the class of  {\sl infinitesimal transformations} having as $\a^{p+1}$ germ the time--one flow of $\a^{p+1}\psi_{p+1}$, where
\beqno
\psi_{p+1}:=\sum_{k\ne 0}\frac{{\rm P}^{(p+1)}_k(I;u,v)}{{\rm i}k\cdot\o(I)}e^{{\rm i}k\cdot\varphi}+\bar\psi_p\eeqno
if ${\rm P}_{p+1}$ has the Fourier expansion
$${\rm P}_{p+1}=\sum_{k\ne 0}{\rm P}^{(p+1)}_k(I;u,v)e^{{\rm i}k\cdot\varphi}$$
and $\bar\psi_p$ is any function independent of $\varphi$. Moreover, as it is known, $\bar{\rm P}_{j}$'s and ${\rm P}_{j}$'s
are related by
$$\bar{\rm P}_{p+1}=({\rm P}_{p+1})_{\rm av}=\frac{1}{(2\p)^n}\int_{\torus^n}{\rm P}_{p+1}d\varphi\ .$$
Therefore, if we perform two steps of the procedure, \ie, with $p=0$, $1$, we find the {\sl two--step normal form}  is defined  by $\bar{\rm P}_1={\rm P}_{\rm av}=0$ and
\beqno
\bar{\rm P}_2=\frac{1}{2}\frac{1}{(2\p)^n}\int_{\torus^n}\{\psi_1, {\rm P}\}d\varphi\ ,\eeqno
where $\{\cdot\ ,\ \cdot\}$ denotes Poisson parentheses with respect to all the variables. (The relative  transformation
  will be given by $\phi_1\circ\phi_2$, where $\phi_j$ is generated by $\a^j\psi_j$.) 
Therefore,  to prove uniqueness, all we have to do is to check that, if we change $\psi_1\to \psi_1+\tilde\psi_1$, where $\tilde\psi_1$ is independent of $\varphi$, the function $\bar{\rm P}_2$ does not change. And in fact  this term changes   by adding
$$\frac{1}{2}\frac{1}{(2\p)^n}\int_{\torus^n}\{\tilde\psi_1, {\rm P}\}d\varphi\ .$$
Since $\tilde\psi_1$ is independent of $\varphi$, Poisson parentheses and the integral may be exchanged and we see that this term vanishes
$$\frac{1}{2}\frac{1}{(2\p)^n}\int_{\torus^n}\{\tilde\psi_1, {\rm P}\}d\varphi=\frac{1}{2}\frac{1}{(2\p)^n}\Big\{\tilde\psi_1,\int_{\torus^n}\ {\rm P}d\varphi\Big\}=\Big\{\tilde\psi_1,\frac{1}{2}{\rm P}_{\rm av}\Big\}=0$$
because of \equ{zero average}.  \qed
\subsection{Proof of Proposition \ref{NF3BP}}\label{proof of 3.1}
To prove Proposition \ref{NF3BP}, we write ${\rm P}^\ppd$
as
\beq{P1}{\rm P}^\ppd=\frac{1}{2\p}\int_{\torus}{\rm P}^\ppu(\L,\ul_1,\uz) d\ul_1\eeq 
where
\beqa{tildeP}
{\rm P}^\ppu(\L,\ul_1,\uz)&:=&\frac{1}{2\p}\int_{\torus}d\ul_2\nonumber\\
&&\frac{3( x^\ppu(\L_1,\ul_1,\uz)\cdot x^\ppd( \L_2,\ul_2, \uz))^2-| x^\ppu(\L_1,\ul_1,\uz)|^2| x^\ppd(\L_1,\ul_1,\uz)|^2}{2|x^\ppd( \L_2,\ul_2, \uz)|^5}\nonumber\\
\eeqa
Then we consider the auxiliary Hamiltonian
\beqno
{\rm H}_{\rm Dip}(y^\ppu, x^\ppu, y^\ppd, x^\ppd):=\frac{|y^\ppd|^2}{2m_2}-\frac{m_2M_2}{|x^\ppd|}-\a m_2M_2\frac{x^\ppu\cdot x^\ppd}{|x^\ppd|^3}\eeqno
on the phase space
$$\{(y^\ppu, x^\ppu, y^\ppd, x^\ppd)\in (\real^3)^4:\ x^\ppd\ne 0\}$$
endowed with the standard symplectic form
$$\o:=dy^\ppu\wedge dx^\ppu+dy^\ppd\wedge dx^\ppd$$
and  $\a\ll1 $ a small positive parameter.
 \vskip.1in
\noi
For $\a=0$, ${\rm H}_{\rm Dip}$ reduces to the two--body Hamiltonian
$${\rm H}_{\rm 2b}=\frac{|y^\ppd|^2}{2m_2}-\frac{m_2M_2}{|x^\ppd|}\ .$$ 
Therefore, letting $$(\L_i,\ul_i,\uz)\to (y^\ppi(\L_i,\ul_i,\uz), x^\ppi(\L_i,\ul_i,\uz))$$ the projection over the $i^{\rm th}$ planet
of the  map \equ{map**}, in such variables, ${\rm H}_{\rm Dip}$ takes the  form
\beqa{initial}
{\rm H}(\L_1,\L_2,\ul_1,\ul_2, \uz)&=&{\rm H}_{\rm Kep}(\L_2)+\a{\rm P}(\L_1,\L_2,\ul_1,\ul_2, \uz)\nonumber\\
&=&-\frac{M_2^2 m_2^3}{2\L_2^3}-\a M_2 m_2\frac{x^\ppu(\L_1, \ul_1,\uz)\cdot x^\ppd(\L_2, \ul_2,\uz)}{|x^\ppd(\L_2, \ul_2, \uz)|^3}\ .\eeqa

\begin{lemma}\label{normal form lemma}
Under assumptions of Proposition \ref{NF3BP},
the Hamiltonian  in \equ{initial}, endowed with the symplectic form
$$\sum_{i=1}^2d\L_i\wedge d\l_i+\sum_{i=1}^4du_i\wedge dv_i\ ,\qquad \uz=(u,v)$$
 verifies the assumptions of Lemma \ref{NFT}, with to the ``variables'' $(I,\varphi):=(\L_2, \ul_2)$ and the ``parameters'' $(\L_1,\ul_1, \uz)$. Its  (unique) two--step  normal form is 
 \beqno
 \tilde {\rm H}(\L_1,\L_2,\ul_1, \uz)={\rm H}_{\rm Kep}(\L_2)+\a^2M_2m_2 {\rm P}^\ppu(\L_1,\L_2, \ul_1, \uz)+{\rm O}(\a^3)\eeqno
with ${\rm P}^\ppu$ as in \equ{tildeP}.
\end{lemma}
\proof We   apply Lemma \ref{NFT} to the Hamiltonian ${\rm H}$ in \equ{initial}.
Indeed,
the assumption \equ{Newton law}  
implies that the zero--averaging (with respect to $\ul_2$) assumption for ${\rm P}$ is satisfied:
\beqano
\int_{\torus}{\rm P}d\ul_2&=&\int_{\torus}(- m_2M_2\frac{x^\ppu\cdot x^\ppd}{|x^\ppd|^3})d\ul_2=- m_2M_2x^\ppu\cdot \int_{\torus}\frac{x^\ppd}{|x^\ppd|^3}d\ul_2\nonumber\\
&=& m_2M_2\o^\ppd_{\rm Kep}x^\ppu\cdot \int_{\torus}\partial_{\ul_2}y^\ppd=0
\eeqano
Denote as
\beqno\tilde {\rm H}(\L_1,\L_2,\ul_1, \uz)={\rm H}_{\rm Kep}(\L_2)+\a^2M_2m_2 {\rm P}^\ppu(\L_1,\L_2; \ul_1, \uz)+{\rm O}(\a^3)\eeqno
with $$\dst{\rm H}_{\rm Kep}(\L_2)=-\frac{M_2^2m_2^3}{2\L_2^2}$$
the two--step normal form which is achieved via Lemma \ref{NFT}. Let $\psi$ denote the symplectic, $\a$--close--to--the identity transformation realizing this normal form.
Consider  the auxiliary Hamiltonian
\beq{Hstar}{\rm H}^\star={\rm H}+\a^2 M_2 m_2{\rm Q}\ ,\eeq
where
$${\rm Q}:=-\frac{3( x^\ppu(\L_1,\ul_1,\uz)\cdot x^\ppd( \L_2,\ul_2, \uz))^2-| x^\ppu(\L_1,\ul_1,\uz)|^2| x^\ppd(\L_1,\ul_1,\uz)|^2}{2|x^\ppd( \L_2,\ul_2, \uz)|^5}\ .$$
Being $\a$--close to the identity,  $\psi$ transforms ${\rm H}^\star$ into
\beqno\tilde {\rm H}^\star=\tilde{\rm H}+\a^2 M_2 m_2{\rm Q}+{\rm O}(\a^3)\ .\eeqno
Hence, at expenses of  a further $\ul_2$--averaging ($\a^2$--close to the identity), $\widetilde {\rm H}$ can be let into
\beqano\widetilde {\rm H}(\L_1,\L_2,\ul_1, \uz)&=&\tilde{\rm H}+\a^2 M_2 m_2{\rm Q}^\ppu(\L_1,\L_2; \ul_1, \uz)+{\rm O}(\a^3)\nonumber\\
&=&{\rm H}_{\rm Kep}(\L_2)+\a^2M_2m_2 {\rm P}^\ppu(\L_1,\L_2; \ul_1, \uz)+\a^2 M_2 m_2{\rm Q}^\ppu(\L_1,\L_2; \ul_1, \uz)\nonumber\\
&+&{\rm O}(\a^3)\ ,\eeqano
with
$${\rm Q}^\ppu:=\frac{1}{2\p}\int_{\torus}{\rm Q}d\ul_2$$
On the other hand, 
one one immediately sees that ${\rm H}^\star$ in \equ{Hstar} may be written as
$${\rm H}^\star={\rm H}^\star_{\textrm{\sc 2b}}+{\rm O}(\a^3)$$
${\rm H}^\star_{\textrm{\sc 2b}}$ is the well ``familiar'' one
$${\rm H}^\star_{\textrm{\sc 2b}}:=\frac{|y^\ppd(\L_2,\ul_2,\uz)|^2}{2m_2}-\frac{m_2M_2}{|x^\ppd(\L_2,\ul_2,\uz)-\a x^\ppu(\L_1,\ul_1,\uz)|}\ . $$
But  (using, note, assumption (B)) ${\rm H}^\star_{\textrm{\sc 2b}}$ may be symplectically conjugated, via an $\a$--close to the identity map $\psi'$, to $\dst{\rm H}_{\rm Kep}=-\frac{M_2^2m_2^3}{2\L_2^2}$. This implies that $\psi'$ lets ${\rm H}^\star$ into
$$\widetilde{\rm H}'=-\frac{M_2^2m_2^3}{2\L_2^2}+{\rm O}(\a^3)$$
Uniqueness (claimed by Lemma \ref{NFT}) implies $\widetilde {\rm H}\equiv \widetilde{\rm H}'+{\rm O}(\a^3)$, namely,
$${\rm P}^\ppu(\L_1,\L_2; \ul_1, \uz)=-{\rm Q}^\ppu(\L_1,\L_2; \ul_1, \uz)$$
which is the thesis.
 \qed
 We are now ready for the
 \proof {\bf of Proposition \ref{NF3BP}}
For the purposes of this proof, if $f:x\in\torus\to f(x)\in\real$ is continuous, we denote as $\langle f\rangle_{x}:=\frac{1}{2\p}\int_{\torus}f(x)dx$.

\nl
Consider  the Hamiltonian ${\rm H}$ in \equ{initial}; let $\bar\phi$ as in \equ{map**} and $\phi$ as in \equ{map}. Denote as ${\rm H}_{\rm red}:={\rm H}\circ \bar\phi^{-1}\circ\phi$ the expression of ${\rm H}$ in the variables \equ{def}.
This is 
\beq{reduction}
{\rm H}_{\rm red}={\rm H}\circ \bar\phi^{-1}\circ\phi=\frac{{\rm R}_2^2}{2m_2}-\frac{M_2m_2}{{\rm r}_2}+\frac{\Phi_2^2}{2m_2{\rm r}_2^2}-M_2 m_2\a\frac{{\rm r}_1}{{\rm r}_2^2}\sqrt{1-(\frac{\Theta}{\Phi_2})^2}\sin\varphi_2\ .\eeq
Let us split ${\rm H}_{\rm red}$ into two parts, a ``radial'' and a ``tangential'' one: 
\beqno
{\rm H}_{\rm rad}:=\frac{{\rm R}_2^2}{2m_2}-\frac{M_2m_2}{{\rm r}_2}\eeqno
and
\beqno
{\rm H}_{\rm tan}:=\frac{\Phi_2^2}{2m_2{\rm r}_2^2}-M_2 m_2\a\frac{{\rm r}_1}{{\rm r}_2^2}\sqrt{1-(\frac{\Theta}{\Phi_2})^2}\sin\varphi_2\eeqno
and focus on ${\rm H}_{\rm tan}$. We shall 
 eliminate  the dependence from the angle $\varphi_2$ up to order $\a^3$. To this end, define ${\rm h}_0$, ${\rm P}_0$ via ${\rm H}_{\rm tan}=:{\rm h}_0+\a{\rm P}_0$ and denote $\varpi:=\partial_{\Phi_2}{\rm h}=\frac{\Phi_2}{m_2{\rm r}_2^2}$. Since $\langle{\rm P}_0\rangle_{\varphi_2}=0$, a Hamiltonian vector field the time--one flow of which eliminates the dependence on $\varphi_2$ up to ${\rm O}(\a^2)$ has as  Hamiltonian the function $\psi_0$ defined as a primitive 
\beqano
\psi_0&=&\frac{1}{\varpi}\int^{\varphi_2}\a{\rm P}_0=
M_2 m_2^2\a\frac{{\rm r}_1}{\Phi_2}\sqrt{1-(\frac{\Theta}{\Phi_2})^2}\cos\varphi_2\eeqano
with $ \langle\phi_0\rangle_{\varphi_2}=0$.
It is a remarkable fact that ${\rm r}_2$ is cancelled. Since   $\phi_0$ is also independent of ${\rm R}_1$, ${\rm R}_2$ and $\vartheta$, this implies that its time--one flow, that we denote
\beqno
\phi_0:\quad (\tilde{\rm R}_1, \tilde{\rm R}_2, \tilde\Phi_2,\tilde\Theta, \tilde{\rm r}_2,\tilde{\rm r}_2, \tilde{\varphi}_2,\tilde\vartheta)\to ({\rm R}_1, {\rm R}_2, \Phi_2,\Theta, {\rm r}_2,{\rm r}_2, {\varphi}_2,\vartheta)\ ,\eeqno
 leaves $({\rm R}_2, {\rm r}_2, \Theta, {\rm r}_1)$  unvaried. Using again $\langle {\rm P}_0\rangle_{\varphi_2}=0$, we then have that ${\rm H}_0$ is conjugated to
$${\rm H}_1={\rm H}_{\rm tan}\circ\phi_0={\rm h}_0+\a^2{\rm P}_1+{\rm O}(\a^3)\ ,$$ 
where 
\beqano
{\rm P}_1=\frac{1}{2}\{\psi_0, {\rm P}_0\}=-\frac{M_2^2m_2^3}{2}\frac{\tilde{\rm r}_1^2}{\tilde{\rm r}_2^2\tilde\Phi_2^4}\Big(\tilde\Theta^2-\frac{1}{2}(\tilde\Phi_2^2-\tilde\Theta^2)(1+\cos2\tilde\varphi_2)\Big)\eeqano

\nl
A further step of averaging defined by the time--one flow 
\beqno
\phi_1:\quad (\hat{\rm R}_1, \hat{\rm R}_2, \hat\Phi_2,\hat\Theta, \hat{\rm r}_2,\hat{\rm r}_2, \hat{\varphi}_2,\hat\vartheta)\to (\tilde{\rm R}_1, \tilde{\rm R}_2, \tilde\Phi_2,\tilde\Theta, \tilde{\rm r}_2,\tilde{\rm r}_2, \tilde{\varphi}_2,\tilde\vartheta)
\eeqno
of 
\beqano
\psi_1&=&\frac{1}{\varpi}\int^{\varphi_2}\a^2({\rm P}_1-\langle{\rm P}_1\rangle)
\nonumber\\
&=&+\frac{\a^2}{\varpi}\int^{\varphi_2}\frac{M_2^2m_2^3}{4}\frac{{\rm r}_1^2}{{\rm r}_2^2\Phi_2^4}(\Phi_2^2-\Theta^2)\cos2\varphi\nonumber\\
&=&+\a^2\frac{M_2^2m_2^4}{8}\frac{{\rm r}_1^2}{\Phi_2^5}(\Phi_2^2-\Theta^2)\sin2\varphi_2
\eeqano
with $\langle\psi_1\rangle_{\varphi_2}=0$. As in the previous step, 
$\psi_1$ is  independent of $({\rm R}_1, {\rm R}_2, \vartheta)$ and, again ${\rm r}_2$, hence, $\phi_1$ leaves $({\rm R}_2, {\rm r}_2, \Theta, {\rm r}_1)$ unvaried. 
Then ${\rm H}_{1}$ is let into the form
$${\rm H}_2={\rm H}_1\circ\phi_1={\rm h}_0+\a^2{\rm P}_2+{\rm O}(\a^3)\ ,$$ 
where 
$${\rm P}_2=\langle{\rm P}_1\rangle_{\varphi_2}=-\frac{M_2^2m_2^3}{4}\frac{\hat{\rm r}_1^2}{\hat{\rm r}_2^2\hat\Phi_2^4}(3{\hat\Theta^2}-\hat\Phi_2^2)\ .$$
Including also the term ${\rm H}_{\rm rad}$ (left unvaried by this sequence of transformations) we finally have that the Hamiltonian ${\rm H}_{\rm red}$ in \equ{reduction} is transformed into 
\beqa{fin}
\hat{\rm H}&:=&{\rm H}_{\rm red}\circ\phi_0\circ\phi_1=\frac{\hat{\rm R}_2^2}{2m_2}-\frac{M_2m_2}{\hat{\rm r}_2}+\frac{1}{2m_2\hat{\rm r}_2^2}\Big({\hat\Phi_2^2}-\a^2\frac{M_2^2m_2^4}{2}\frac{\hat{\rm r}_1^2}{\hat\Phi_2^4}(3{\hat\Theta^2}-\hat\Phi_2^2)\Big)\nonumber\\
&&+{\rm O}(\a^3)\ .
\eeqa

\vskip.1in
\noi
Let now\beqno
(\hat y^\ppu, \hat y^\ppd, \hat x^\ppu, \hat x^\ppd)\in \real^3\times \real^3\times \real^3\times \real^3\eeqno be related to $({\rm C}_3, {\rm G}, \hat{\rm R}_1, \hat{\rm R}_2, \hat\Phi_2,\hat\Theta, \zeta, {\eufm g},\hat{\rm r}_1,\hat{\rm r}_2, \hat{\varphi}_2,\hat\vartheta)$ via  relations analogue to \equ{map}--\equ{incli}, \ie, $(\hat y, \hat x)=\phi^{-1}({\rm C}_3, {\rm G}, \hat{\rm R}_1, \hat{\rm R}_2, \hat\Phi_2,\hat\Theta, \zeta, {\eufm g},\hat{\rm r}_1,\hat{\rm r}_2, \hat{\varphi}_2,\hat\vartheta)$, with $\phi^{-1}$ as in \equ{def} and let $(\hat\L, \hat\ul, \hat \uz)$ be  defined via $(\hat y, \hat x)=\bar\phi(\hat\L,\hat\ul,\hat \uz)$, with $\bar\phi$ as in \equ{map**}. In the variables $(\hat\L, \hat\ul, \hat\uz)$, the he Hamiltonian \equ{fin} takes the form 

$$\tilde{\rm H}:=\hat{\rm H}\circ\phi^{-1}\circ\bar\phi=-\frac{M_2^2 m_2^3}{2\hat\L_2^2}-\frac{\a^2}{2m_2\hat{\rm r}_2^2}\frac{M_2^2m_2^4}{2}\frac{\hat{\rm r}_1^2}{\hat\Phi_2^4}(3{\hat\Theta^2}-\hat\Phi_2^2)$$
where $\hat\Theta=\hat{\rm C}^\ppd\cdot\hat x^\ppu$, $\hat{\rm r}_i=|\hat x^\ppi|$, $\hat\Phi_2=|\hat{\rm C}^\ppd|$, with $\hat{\rm C}^\ppd=\hat x^\ppd\times \hat y^\ppd$  have to be regarded as functions of $(\L,\hat\ul,\hat\uz)$.
A further $\hat\ul_2$--averaging, $\a^2$--close to the  identity
$$\widehat\phi:\qquad(\widehat\L, \widehat\ul, \widehat\uz)\to (\hat\L, \hat\ul,  \hat\uz)$$
 transforms $\tilde{\rm H}$  into 
\beq{averaged}\widehat{\rm H}:=\tilde{\rm H}\circ\widehat\phi=-\frac{M_2^2 m_2^3}{2\widehat\L_2^2}-\a^2\widehat{\rm P}(\widehat\L,\widehat\ul_1,\widehat \uz)\ ,\quad \widehat{\rm P}(\widehat\L,\widehat\ul_1,\widehat \uz):=\frac{M_2^2m_2^3}{4}\frac{\widehat{\rm r}_1^2}{\widehat\Phi_2^4}(3{\widehat\Theta^2}-\widehat\Phi_2^2)\frac{1}{2\p}\int_{\torus}\frac{d\widehat\l_2}{\widehat{\rm r}_2^2}\ ,\eeq
where $\widehat\Theta=\widehat{\rm C}^\ppd\cdot\widehat x^\ppu$, $\widehat{\rm r}_i=|\widehat x^\ppi|$, $\widehat\Phi_2=|\widehat{\rm C}^\ppd|$ have to be regarded as functions of $(\widehat\L, \widehat\l, \widehat z)$. Note that we have used that $\widehat\Theta=\widehat{\rm C}^\ppd\cdot\widehat x^\ppu$, $\widehat{\rm r}_1=|\widehat x^\ppu|$, $\widehat\Phi_2=|\widehat{\rm C}^\ppd|$ are independent of $\widehat\ul_2$.
By construction, the overall change
$$\bar\phi^{-1}\circ\phi\circ\phi_0\circ\phi_1\circ\phi^{-1}\circ\bar\phi
:\qquad (\widehat\L, \widehat\ul, \widehat \uz)\to (\L, \ul,  \uz)$$
is symplectic, $\a$--close to the identity and puts the Hamiltonian ${\rm H}$ in \equ{initial} into the form claimed in  Lemma \ref{NFT}. By the uniqueness claimed by this theorem, in comparison with the result of Lemma \ref{normal form lemma}, we  have that $\widehat{\rm P}(\widehat\L,\widehat\ul_1,\widehat\uz)$ in \equ{averaged} satisfies 
$$\a^2\widehat{\rm P}(\widehat\L,\widehat\ul_1,\widehat\uz)=-\a^2\frac{M_2^2m_2^3}{4}\frac{\widehat{\rm r}_1^2}{\widehat\Phi_2^4}(3{\widehat\Theta^2}-\widehat\Phi_2^2)\frac{1}{2\p}\int_{\torus}\frac{d\widehat\l_2}{\widehat{\rm r}_2^2}\equiv \a^2M_2 m_2{\rm P}^\ppu(\widehat\L,\widehat\ul_1,\widehat\uz)+{\rm O}(\a^3)\ ,$$
where ${\rm P}^\ppu$  is as in \equ{tildeP} (and, as above, $\widehat\Theta$, $\widehat\Phi_2$, $\widehat{\rm r}_1$ and $\widehat{\rm r}_1$ are regarded as functions of $(\widehat\L_1,\widehat\L_2,\widehat\ul_1,\widehat\ul_2, \widehat\uz)$). Taking the average with respect to $\widehat\ul_1$, we have the thesis.
\qed
\subsection{Proof of Proposition \ref{asymptotic formula}}\label{proof of 2.1}
 We shall need definitions and a result from \cite{chierchiaPi11c}, to which paper we refer for notations and details. 

\vskip.1in
\noi
Let, as in \cite{chierchiaPi11c}, $\cP_{\rm P}^{6n}$, $\cP_{\rm rps}^{6n}\subset \real^{3n}\times \real^{3n}$ denote the respective domains of the maps
$$\phi_{\rm P}:\ (y, x)\in\cP_{\rm P}^{6n}\to (\L,\ul,\uz)\in \real^n\times\torus^n\times \real^{4n}\ ,\quad \phi_{\rm rps}:\ (y, x)\in\cP_{\rm rps}^{6n}\to (\L,\l,z)\in \real^n\times\torus^n\times \real^{4n}$$
bbetween ``Cartesian'' and, respectively, Poincar\'e, {\sc rps} variables.
Consider the common domain   of $\phi_{\rm P}$ and $\phi_{\rm rps}$, \ie the set $\cP^{6n}_{\rm rps}\cap\cP^{6n}_{\rm P}$. 
On the $\phi_{\rm rps}$--image of such domain consider the symplectic map 
\beq{P/RPS map}
\phi_{\rm P}^{\rm rps}:\quad (\L,\l,z) \to (\L,\ul,\uz):=\phi_{\rm P}\circ\phi_{\rm rps}^{-1}
\eeq
which maps the {\scshape rps} variables onto the Poincar\'e variables. Such a map has a particularly simple structure:

 \begin{theorem}[\cite{chierchiaPi11c}]\label{prop: Poinc e Dep}
The symplectic map $\phi_{\rm P}^{\rm rps}$ in \equ{P/RPS map} has the form
\beqa{PoincDep}
\ul=\l+\varphi(\L,z)\qquad \uz={\cal Z}(\L,z)
\eeqa
where $\varphi(\L,0)=0$ and, for any fixed $\L$, the  map ${\cal Z}(\L,\cdot)$ is 1:1, symplectic\footnote{I.e., it preserves the two form $d\eta\wedge d\xi+dp\wedge dq$.} and its projections  verify
\beqno
\P_{\uh}{\cal Z}=\eta+{\rm O}(|z|^3)\ ,\ \P_{\ux}{\cal Z}=\xi+{\rm O}(|z|^3)\ ,\ \P_{\up}{\cal Z}={\cal V} p+{\rm O}(|z|^3)\ ,\ \P_{\uq}{\cal Z}={\cal V} q+{\rm O}(|z|^3)\eeqno
for some ${\cal V}={\cal V}(\L)\in {\rm SO}(n)$.
\end{theorem} 
Now we proceed to prove Proposition \ref{asymptotic formula}. Consider  the inverse maps 
\beqano
&&\phi_{\rm rps}^{-1}:\qquad (\L,\l,z)\in{\cal M}^{6n}_{\rm rps}\to  \big(y_{\rm rps}(\L,\l,z), x_{\rm rps}(\L,\l,z)\big)\nonumber\\\nonumber\\
&&\phi_{\rm P}^{-1}:\qquad(\L,\ul,\uz)\in{\cal M}^{6n}_{\rm P}\to  \big(y_{\rm P}(\L,\ul,\uz), x_{\rm P}(\L,\ul,\uz)\big)
\eeqano
with ${\cal M}^{6n}_{\rm rps}:=\phi_{\rm rps}({\cal P}^{6n}_{\rm rps})$, ${\cal M}^{6n}_{\rm P}:=\phi_{\rm rps}({\cal P}^{6n}_{\rm P})$. Let $y^\ppi_{\rm rps}\in \real^3$, $\cdots$ be the $i^{\rm th}$ projection of $y_{\rm rps}$, $\cdots$; \ie, to be defined by
$$y_{\rm rps}=\big(y^\ppu_{\rm rps}, \cdots y^\ppn_{\rm rps}\big)\ ,\quad \cdots$$
Let, finally, 
\beqano
\a_{ij}^2(f_{ij}^{(2)})_{\rm P}&&:=\frac{1}{(2\p)^2}\int_{\torus^2}d\ul_id\ul_j\nonumber\\
&&\quad\frac{3(x_{\rm P}^\ppi(\L,\ul_i, \uz)\cdot x_{\rm P}^\ppj(\L,\ul_j, \uz))^2-|x_{\rm P}^\ppi(\L,\ul_i, \uz)|^2|x_{\rm P}^\ppj(\L,\ul_j, \uz)|^2}{2|x_{\rm P}^\ppj(\L,\ul_j, \uz)|^5}\ .\nonumber\\
\eeqano
and\footnote{As observed in footnote \ref{foot24}, the map $\phi_{\rm rps}^{-1}$ depends explicitly on $(p_n,q_n)$, while  SO(3)--invariant expressions, such as  the right hand side of the formula below, do not. }
\beqano
\a_{ij}^2(f_{ij}^{(2)})_{\rm rps}&&:=\frac{1}{(2\p)^2}\int_{\torus^2}d\l_id\l_j\nonumber\\
&&\quad\frac{3(x^\ppi_{\rm rps}(\L,\l_i, z)\cdot x^\ppj_{\rm rps}(\L,\l_j, z))^2-|x^\ppi_{\rm rps}(\L,\l_i, z)|^2|x^\ppj_{\rm rps}(\L,\l_j, z)|^2}{2|x^\ppj_{\rm rps}(\L,\l_j, z)|^5}\ .\nonumber\\
\eeqano

\nl
 We shall use the following properties, easily deducible from \cite{chierchiaPi11c}:
\begin{itemize}
\item[{\rm (i)}] For $1\le i\le n$, $y^\ppi_{\rm rps}$, $x^\ppi_{\rm rps}$ depend on $\l$ only via $\l_i$. Analogously, $y^\ppi_{\rm P}$, $x^\ppi_{\rm P}$ depend on $\ul$ only via $\ul_i$. In particular $y^\ppi_{\rm P}$, $x^\ppi_{\rm P}$ depend on $\L$ only via $\L_i$ and 
depend on $\uz$ only via $\uz_i$,  
 but this will not be used.
\item[{\rm (ii)}] For any $1\le i<j\le n$, the map \beq{P**}(\L_i,\L_j,\ul_i,\ul_j,\uz_i,\uz_j)\to \big(y_{\rm P}^\ppi(\L, \ul_i. \uz), y_{\rm P}^\ppj(\L, \ul_j, \uz), x_{\rm P}^\ppi(\L, \ul_i. \uz), x_{\rm P}^\ppj(\L, \ul_j. \uz)\big)\eeq satisfies assumptions (A) and (B) of Proposition \ref{NF3BP}. Note that, unless we are in the case $n=2$, this is not true for the map \beq{rps**}(\L,\l_i,\l_j,z)\to \big(y_{\rm rps}^\ppi(\L,\l_i,z), y_{\rm rps}^\ppj(\L,\l_j,z), x_{\rm rps}^\ppi(\L,\l_i,z), x_{\rm rps}^\ppj(\L,\l_j,z)\big)\eeq
In particular,  both \equ{P**} and \equ{rps**} satisfy assumption (A) (for any $n$ and any $1\le i<j\le n$), 
but assumption (B) fails for \equ{rps**} (when $n>2$).
\item[{(iii)}] Letting ${\rm C}_{\rm rps}^\ppi:=x^\ppi_{\rm rps}\times y^\ppi_{\rm rps}$ and, analogously, ${\rm C}_{\rm P}^\ppi:=x^\ppi_{\rm P}\times y^\ppi_{\rm P}$, then, for any $1\le i\le n$, ${\rm C}_{\rm rps}^\ppi$ does not depend on $\l_i$ and, analogously, ${\rm C}_{\rm P}^\ppi$ does not depend on $\ul_i$. This is because, as remarked in (ii), both \equ{P**} and \equ{rps**} satisfy (A). 
\end{itemize}
By the previous items, may apply Proposition \ref{NF3BP} to the map \equ{P**}. We find 
\beqano
&&\a_{ij}^2(f_{ij}^{(2)})_{\rm P}(\L, \uz)=-\frac{M_jm_j^2}{4}\nonumber\\
&{\times}&\frac{\dst\frac{1}{2\p}\int_{\torus}\big(3({\rm C}_{\rm P}^\ppj(\L,\uz)\cdot x_{\rm P}^\ppi(\L,\ul_i,\uz))^2-|x_{\rm P}^\ppi(\L,\ul_i,\uz)|^2|{\rm C}_{\rm P}^\ppj(\L,\ul_j,\uz)|^2\big)d\ul_i}{\dst|{\rm C}_{\rm P}^\ppj(\L,\ul_j,\uz)|^4}\nonumber\\
&{\times}&\frac{1}{2\p}\int_{\torus}\frac{d\ul_j}{\dst |x_{\rm P}^\ppj(\L,\ul_j,\uz)|^2}\eeqano
Letting now $\uz={\cal Z}(\L,z)$ and changing the integration variables $\ul_i=\l_i+\varphi_i(\L,z)$ with ${\cal Z}$, $\varphi$ as in \equ{PoincDep}
we have the thesis. \qed

\section{Proof of Theorem {\rm A}}\label{ArnoldConj}
In this section, we aim to prove Theorem A.
\begin{remark}\rm
For definiteness, we prove Theorem A for the spatial three--body problem. In the case of the planar $(1+n)$--body problem with\footnote{For $n=2$ there is the stronger result of Theorem \ref{simpler planar}.} $n\ge 3$, assume the following asymptotic of semi--axes
\beq{asymp1} \underline a_j\le a_{j}\le \bar a_j\eeq
where
\beq{asymp2}\underline a_n:=\underline a\ ,\quad \bar a_n:=\bar a\ ,\quad \underline a_j:=c\underline\a^{(\frac{3}{2})^{n-j}} \underline a_{n}\ ,\quad  \bar a_j:=\underline\a^{(\frac{3}{2})^{n-j}} \underline a_{n}\eeq
where $0<\underline a<\bar a$ and $0<\underline\a<c<1$ are fixed. With this assumption, the rest of the proof of this case is similar to the one of the spatial three--body case presented below. Let us sketch it briefly.
An analogue splitting as in \equ{good split} below is available, with $N$, $\tilde N$ replaced by
$$N':=-\sum_{1\le i<j\le n}\bar m_i \bar m_j \sum_{k\in \{0,2\}}f^{(k)}_{ij}|_{\rm pl}\ ,\quad \widetilde N':=-\sum_{1\le i<j\le n}\bar m_i \bar m_j \sum_{k=3}^{\infty}f^{(k)}_{ij}|_{\rm pl}\ .$$
The integrability of $N'$ has been discussed in the Introduction. Moreover,  the functions in $f^{(2)}_{ij}|_{\rm pl}$ have all the same strength, $\frac{1}{a_2}(\underline\a)^{(\frac{3}{2})^{n-2}}$, while the ones with $k\ge 3$ (hence, the remainder $\tilde N$) are of order $\frac{\underline\a^{3/4}}{a_2}(\underline\a)^{(\frac{3}{2})^{n-2}}$. The remaining details are left to the reader. \qed
%
\end{remark}

\nl
Let us consider the spatial three--body Hamiltonian
\beq{sp3bp}\cH_{\rm 3b}:=h_{\rm Kep}(\L)+\m f_{\rm 3b}(\L,\l,\bar z)\eeq
namely, the Hamiltonian  $\cH_{\rm rps}$ in \equ{rps ham} for $n=2$.
Let $f^{(k)}_{ij}$ be as in \equ{alpha expansion}; define

\beq{contributions}N:=-\bar m_1 \bar m_2 \sum_{j\in \{0,2\}}f^{(j)}_{12}\ ,\quad \widetilde N:=-\bar m_1 \bar m_2 \sum_{j=3}^{\infty}f^{(j)}_{12}
\eeq 
so as to split
\beq{good split}(f_{\rm 3b}(\L,\l,\bar z))_{\rm av}=N+\widetilde N\ .\eeq
with $N$ integrable and $|\tilde N|\le \const\a^3$. Integrability of $N$ is known since \cite{lidovZ76} and will be discussed in this setting in Claim \ref{integration}. 

\subsection{Symmetries of the partially reduced system}\label{partial reduction}
We  recall some properties discussed in \cite{chierchiaPi11b} and \cite{chierchiaPi11c}, to which we refer for more details.

\nl
The Hamiltonian \equ{planetary} remains unvaried by   reflections with respect to coordinate planes $\{x_1=x_2\}$,  $\{x_3=0\}$ or rotations, for example, around the $k^\ppt$--axis. These transformations are, respectively,
\beqano
&&\begin{array}{llllllll}
\cR_{{}_{1\leftrightarrow 2}}:\quad &x^\ppi\to\big(x^\ppi_2,\ x^\ppi_1,\ x^\ppi_3\big)\ ,\quad &y^\ppi\to\big(-y^\ppi_2,\ -y^\ppi_1,\ -y^\ppi_3\big)\\  \\
\cR_3^-:\quad &x^\ppi\to\big(x^\ppi_1,\ x^\ppi_2,\ -x^\ppi_3\big)\ ,\quad &y^\ppi\to\big(y^\ppi_1,\ y^\ppi_2,\ -y^\ppi_3\big)\\ \\
\cR_g:\quad &x^\ppi\to {\rm R}_3(g)\,x^\ppi\ ,\quad &y^\ppi\to {\rm R}_3(g)\,y^\ppi
\end{array}\eeqano
where ${\rm R}_3(g)$ denotes the matrix
\beqno{\rm R}_3(g):=\left(
\begin{array}{lcc}
\cos{g}\ &-\sin{g}&0\\
\sin{g}\ &\cos{g}&0\\
0&0&1
\end{array}
\right)\ ,\qquad g\in \torus\ .\eeqno
Note, in particular, that $\cR_3^-$ and $\cR_g$ are symplectic transformations, while $\cR_{{}_{1\leftrightarrow 2}}$ is an involution.
The expressions of $\cR_{{}_{1\leftrightarrow 2}}$, $\cR_3^-$ and $\cR_g$ in terms of the variables \equ{z}  turn out  to be the same\footnote{See, for example \cite{herman09}.}  as in Poincar\'e variables. They are
\beqa{rot refl}
\begin{array}{llllllll}
\cR_{{}_{1\leftrightarrow 2}}\Big(\L,\ \l,\ z\Big):=\Big(\L,\ \frac{\p}{2}-\l,\ \cS_{{}_{1\leftrightarrow 2}} z\Big)\ ;\quad &\cR_3^-\Big(\L,\ \l,\ z\Big)=\Big(\L,\ \l,\ \cS_{34}^-z\Big)\\ \\
\cR_g\Big(\L,\ \l,\ z\Big)=\Big(\L,\ \l+g,\ \cS_gz\Big)& 
\end{array}
\eeqa
where 
\beqano
\arr{\cS_{{}_{1\leftrightarrow 2}}(\eta,\xi,p,q):=(\xi,\eta,q,p)\\ \\ \cS_{34}^-(\eta,\xi,p,q):=(\eta,\xi,-p,-q)\\ \\
\cS_{g}: \ \Big(\eta_j+{\rm i}\xi_j, p_j+{\rm i}q_j\Big)\to \Big(e^{-{\rm i}g}(\eta_j+{\rm i}\xi_j)\ ,\  e^{-{\rm i}g}(p_j+{\rm i}q_j)\Big)
}
\eeqano
with ${\rm i}:=\sqrt{-1}$.

\nl
Since also the Hamiltonian $\cH_{\rm rps}$ \equ{rps ham} is independent of $(p_n,q_n)$, in the above transformations, we may neglect this latter couple of variables and replace\footnote{Recall the definitions in \equ{z}--\equ{bar z}.} $z$ with $\bar z$ in \equ{rot refl}. In particular, the one--parameter group $\{\bar\cR_g\}_{g\in \torus}$ defined by
\beq{Rg}\bar\cR_g:\ (\L,\ \l,\ \bar z,\ p_n,\ q_n)\to (\L,\ \l+g,\ {\cal S}_g\bar z,\ p_n,\ q_n)\qquad g\in \torus\eeq
leaves $\cH_{\rm rps}$ unvaried. This group of transformations corresponds to be the time--$g$ flow of
\beq{G**}G=\sum_{i=1}^n\L_i-\sum_{i=1}^{n}\frac{\eta_i^2+\xi_i^2}{2}-\sum_{i=1}^{n-1}\frac{p_i^2+q_i^2}{2}\eeq
which is the Euclidean length of the angular momentum \equ{C}: $G=|{\rm C}|$, expressed in the variables \equ{z}. Therefore, $\bar\cR_g$ may be identified to be the group $g$--rotations  about the ${\rm C}$--axis.

\nl
In view of such relations, amusing symmetries (discussed\footnote{In \cite{chierchiaPi11b}, $\bar\cR_g$--invariance is called ``rotation invariance''. Here, to avoid confusions, we reserve this name only to the transformations \equ{rotation invariance}.}  in \cite{chierchiaPi11b}) appear among the Taylor coefficients of the expansion of the perturbation $f_{\rm rps}$ and hence also of its averaged value $(f_{\rm rps})_{\rm av}$. These symmetries are often referred to (for the classical Poincar\'e system \equ{prop deg}) as {\sl D' Alembert rules}.
 To describe such relations, we switch\footnote{$d\eta_i\wedge d\xi_i=d w_i\wedge dw_i^*$ and $dp_j\wedge dq_j=d w_{j+n}\wedge dw_{j+n}^*$} to ``Birkhoff coordinates''
\beq{Bir coord}w_i=\frac{\eta_i-{\rm i}\xi_i}{\sqrt 2}\ ,\quad w_{n+j}=\frac{p_j-{\rm i}q_j}{\sqrt 2}\ ,\quad w_i^\star=\frac{\eta_i+{\rm i}\xi_i}{{\rm i}\sqrt 2}\ ,\quad w^\star_{n+j}=\frac{p_j+{\rm i}q_j}{{\rm i}\sqrt 2}\eeq
with $1\le i\le n$ and $1\le j\le n-1$ and we regard  (abusively) $f_{\rm rps}$ and $(f_{\rm rps})_{\rm av}$ as  functions of $(\L,\l, w,w^\star)$.

\begin{claim}[\cite{chierchiaPi11b}, \cite{chierchiaPi11c}]\label{sym}
\item[{\rm (i)}] $\cR_3$--invariance implies that $f_{\rm rps}$ is even in $(w_{n+1}, \cdots w_{2n-1}, w^\star_{n+1}, \cdots w^\star_{2n-1})$ (equivalently, it is even in $(\bar p,\bar q)$); 
\item[{\rm (ii)}]
$\bar\cR_g$--invariance implies that, the only non--vanishing monomials  appearing  in the Taylor expansion of $(f_{\rm rps})_{\rm av}$ in powers $\{w_i$, $w_i^*\}_{1\le i\le 2n-1}$  are those with literal part $w^ {\a}
{w^*}^{\a^*}$ for which
\beq{a=a*}\sum_{i=1}^{2n-1}(\a_i-\a_i^*)=0\ .\eeq
\end{claim}
 Claim \ref{sym} and  the independence of $f^\ppd_{12}$ on the argument of $(\eta_2,\xi_2)$ (see the Introduction)  have the following corollary. Let $\cA=\cA(\a)$ denote a set of the form
\beq{cA}\cA:=(\L_1,\L_2):\ a_-<\frac{1}{M_1}(\frac{\L_1}{m_1})^2<\frac{\a }{M_2}(\frac{\L_2}{m_2})^2<a_+\eeq
(with $a_-<a_+$, $\a\in (0,1)$) and let $\cM^{10}_{\e_0}:=\cA\times\torus^2\times B_{\e_0}^6$.
\begin{claim} \label{integration}
$N$ (namely\footnote{Recall that $f_{12}^{(0)}$ is independent of $\bar z$.}, $f_{12}^{(2)}$) is integrable. More precisely: (i) it depends on $(\eta_2,\xi_2)$ only via $\frac{\eta_2^2+\xi_2^2}{2}$; (ii) one can find $\e_0>0$ and a symplectic change of variables
$$(\L,\breve\l,\breve z)\to(\L,\l,\bar z)$$
defined on the phase space $\cM^{10}_{\e_0}:=\cA\times\torus^2\times B_{\e_0}^6$  of the form
\beq{hat psi}\breve\phi:\quad \L=\L\ ,\quad \l=\breve\l+\varphi(\L,\breve z)\ ,\quad \bar z=\breve{\cal Z}(\L,\breve z)\eeq
defined for
$|\breve z|<\e_0$
which transforms $N$ into a new function $\breve N(\L,\breve z)$ depending only on $\dst\frac{\breve\eta_1^2+\breve\xi_1^2}{2}$, $\dst\frac{\breve\eta_2^2+\breve\xi_2^2}{2}$ and $\dst\frac{\breve p_1^2+\breve q_1^2}{2}$. In particular, $\breve\psi$ preserves $\dst\frac{\breve\eta_2^2+\breve\xi_2^2}{2}$ and $\dst\frac{\breve\eta_1^2+\breve\xi_1^2}{2}+\frac{\breve p_1^2+\breve q_1^2}{2}$.
\end{claim}
\proof Since $f^\ppd_{12}$ is even in $(p_1,q_1)$ and has only monomials with $\a_2=\a_2^\star$, 
Equation \equ{a=a*} with $n=2$ implies that $f^\ppd_{12}$ is even in $(\eta_1,\xi_1)$, $(\eta_2,\xi_2)$ and $(p,q)$ separately.
Moreover, 
$f^\ppd_{12}$  is integrable\footnote{To integrate $f^\ppd_{12}$, one can first reduce the integral $G_0:=\frac{\breve\eta_1^2+\breve\xi_1^2}{2}+\frac{p_1^2+q_1^2}{2}$
via the change of variables
$$\eta_1+{\rm i}\xi_1=(\breve\eta_1+\breve{\rm i}\xi_1)e^{{\rm i}g_0}\ ,\quad p_1+{\rm i}q_1=\sqrt{2(G_0-\frac{\breve\eta_1^2+\breve\xi_1^2}{2})}e^{{\rm i}g_0}$$
with $g_0$ cyclic in $f^\ppd_{12}$ (but not in $f_{\rm 3b}$). Note that this reduction does not cause singularities in $f_{\rm 3b}$, since $f_{\rm 3b}$ is even in $(p_1,q_1)$. Next, once $f^\ppd_{12}$  is reduced to one degree of freedom, its integration is trivial. }. 
Let $\bar z=\breve{\cal Z}(\L,\breve z)$ the transformation (parametrized by $\L$) verifying $$\sum_{i=1}^2d\eta_i\wedge d\xi_i+dp_1\wedge dq_1=\sum_{i=1}^2d\breve\eta_i\wedge d\breve\xi_i+d\breve p_1\wedge d\breve q_1$$ such that $\breve N(\L,\breve z):=\bar N\circ\breve{\cal Z}$ has the claimed properties. Then, it is standard to prove that $\breve z\to \breve{\cal Z}(\L,\breve z)$  may be lifted to a transformation as in \equ{hat psi} (compare, for example, \cite[Proposition 7.3]{chierchiaPi11b}).\qed
\subsection{{\sc kam} Theory}
In this section we complete the proof of Theorem {\rm A}.

\nl
Let $\e_0$, $\breve\phi$ is as in Claim \ref{integration}. For $(\L,\breve\l,\breve z)\in \cM^{10}:=\cA\times\torus^2\times B^6_{\e_0}$, define
\beqa{3b}\breve\cH_{\rm 3b}(\L,\breve\l,\breve z)&:=&\cH_{\rm 3b}\circ\breve\phi(\L,\breve\l,\breve z)\nonumber\\
&=&h_{\rm Kep}(\L)+\m \breve f_{\rm 3b}(\L,\breve\l,\breve z)\eeqa
where $\breve\phi$ is as in Claim \ref{integration}. By Claim \ref{integration}
\beq{splitting}(\breve f_{\rm 3b})_{\rm av}=\breve N+\tilde N\eeq
where $\breve N$ depends only on $\frac{\breve\eta_1^2+\breve\xi_1^2}{2}$, $\frac{\breve\eta_2^2+\breve\xi_2^2}{2}$, $\frac{\breve p_1^2+\breve q_1^2}{2}$ and
\beq{remainder}|\tilde N|\le \const\a^3\ .\eeq
To the system \equ{3b} we shall apply an abstract result (Theorem \ref{thm: ArnoldConjnew} below) that refines and generalizes  Theorem \ref{simplifiedFT}; see Remark \ref{generalization}. This is as follows.
\vskip.1in
\noi
Let $n_1$, $n_2\in\natural$, $B_\e^{2n_2}=\{y\in\real^{2n_2}: |y|<\e\}$ denote the $2n_2$--ball of radius $\e$ and let
\beq{defcP}
\cP_{\e_0}:=V\times \torus^{n_1}\times B^{2n_2}_{\e_0}\eeq
where $V$ is a open, connected set of $\real^{n_1}$. Let
\beq{pndham}
H(I,\f,p,q;\m):=H_0(I)+\m P(I,\f,p,q;\m)\ ,\eeq
be real--analytic on $\cP_{\e_0}$ and such that
\begin{itemize}
\item[{\rm (i)}] $\o_0:=\partial H_0$ is a real--analytic diffeomorphism of $V$;
\item[{\rm (ii)}] the average $\dst P_{\rm av}(I,p,q;\m)=\frac{1}{(2\p)^{n_1}}\int_{\torus^{n_1}}P(I,\f,p,q;\m)d\varphi$ has the form 
\item[] $\dst P_{\rm av}(I,p,q;\m, \a)=N(I,J;\m)+\tilde N(I,p,q;\m)$, where
\item[]  $J=(\frac{p_1^2+q_1^2}{2}, \cdots, \frac{p_{n_2}^2+q_{n_2}^2}{2})$ and  $\sup_{V\times B^{2n_2}}|\tilde N|\le \k$;
\item[{\rm (iii)}] the Hessians $\partial_I^2H_0$ $\partial^2_{I, J}N(I,J;\m)$ do not vanish, respectively, on $V$, $V\times B^{2n_2}_{\e_0}$.
\end{itemize}
\begin{theorem}\label{thm: ArnoldConjnew}
Under the previous assumptions, one can find positive numbers $C_*$, $\m_*$, $\k_\star$, $\e_1<\e_0$ depending only on $H$ and $\e_0$ and an integer $\b$ depending only on $n_1$, $n_2$, such that, for
\beq{simplifiedsmallness} |\m|<\m_*\ ,\quad |\k|<\k_*\ ,\quad |\m|<(\log\k^{-1})^{-2\b}\eeq
 a
set $\cK\subset  \cP_{\e_1}$ exists,  formed by the union of $H$--invariant $n$--dimensional tori, on which the $H$--motion is analytically conjugated to  linear Diophantine quasi--periodic motions. The  set $\cK$  is of positive Liouville--Lebesgue measure and satisfies
\beq{simplified measest}
\meas \cK> \Big(1- C_*({\sqrt[4]\m (\log \k^{-1})^{\b}}+\sqrt\k) \Big) \meas \cP_{\e_1}\ .
\eeq
\end{theorem}
\begin{remark}\label{generalization}\rm
Theorem \ref{thm: ArnoldConjnew} generalizes and refines Theorem \ref{simplifiedFT}: to obtain  Theorem \ref{simplifiedFT} from Theorem \ref{thm: ArnoldConjnew} it is sufficient to take $\k=\m$. In this case condition \equ{simplifiedsmallness} becomes just a smallness condition on $\m$ (as inTheorem \ref{simplifiedFT}) and, by \equ{simplified measest}, $\cK$ fills ${\cP}_{\e_1}$ up  to a set of density $(1-\tilde C\m^a)$ with any $0<a<\frac{1}{4}$. This should be compared with the measure estimate given in Theorem \ref{simplifiedFT}, where $a\sim \frac{1}{n}$.

\nl
The proof of Theorem \ref{thm: ArnoldConjnew}  is sketched in Appendix \ref{kam}.
\end{remark}

\noi
We are now ready to complete the
\proof {\bf of Theorem {\rm A}} Apply Theorem \ref{thm: ArnoldConjnew} to $\breve\cH_{\rm 3b}:=\cH_{\rm3b}\circ\breve\phi$ (where $\breve\phi$ is as in \equ{hat psi}), hence, with
$$n_1=2\ ,\quad n_2=3\ ,\quad V=\cA\ ,\quad \k=\const\a^3\ ,\quad N=\breve N$$
$\tilde N$ as in \equ{splitting} and $\e_0$ as in Claim \ref{integration}. \qed
\section{Proof of Theorem {\rm B}}\label{proof of stability}
In this section, we shall prove the following theorem, which is a more detailed statement of Theorem {\rm B}. 
Let
\beq{pl3bp}\cH_{\rm pl3b}=h_{\rm Kep}+\m f_{\rm pl3b}:=\cH_{\rm 3b}|_{p_1=q_1=0}\eeq
and denote as
\beq{planar domain}\cM^{8}_{\rm pl3b}:\quad\frac{\L_1^2}{M_1m_1^2}\ge \underline a_-\ ,\quad  \frac{M_2 m_2^2}{M_1m_1^2}\frac{\L_1^2}{\L_2^2}\le \a\ ,\quad \underline\e<|z_{\rm pl}|\le\e\le  \bar\e\ ,\quad \l_1,\ \l_2\in \torus\eeq its eight--dimensional phase space. Here, $\cH_{\rm 3b}$ is as in \equ{sp3bp} and $z_{\rm pl}:=(\eta_1,\eta_2,\xi_1,\xi_2)$.\begin{theorem}\label{stability Theorem detailed}
 There exists positive numbers $\bar\e$, $\bar\a$, $\bar\m$, $\bar\b$, $\t$, $\bar K_\star$, $\bar a$, $\bar b$, $\bar c$, $\bar d$  such that, if
$$
0<\a<\bar\a\ ,\quad 0<\m<\bar\m\ ,\quad 
\m<\bar c(\log\e^{-1})^{-\bar\b}
$$
one can find a an open  set $\bar\cM^{8}_{\rm pl3b}\subset \cM^{8}_{\rm pl3b}$ 
defined by the following inequalities for the Keplerian frequencies $\o_{\rm Kep}:=\partial_{\ovl\L}h_{\rm Kep}$ 
$$|\o_{\rm Kep}\cdot k|\ge \frac{\sqrt[4]\m}{\bar c\bar K}\quad \forall k:\ 0<|k|_1\le \bar K$$
with
\beq{Kappa}
\bar K=\bar K_\star\log(\e^{-1})
\eeq
such that for the  $\cH_{\rm pl3b}$--flow starting from $\bar\cM^{8}_{\rm pl3b}$ the following holds. This flow is symplectically  conjugated, via a
$\{\m^{1/12}, \e^2\}$
--close to the identity transformation $\phi$ to a flow
$$t\to (\tilde\L_1(t), \tilde\L_2(t), \tilde\eta_1(t), \tilde\eta_2(t), \tilde\xi_1(t), \tilde\xi_2(t))$$ 
such that, letting $\tilde t_i(t):=\frac{\tilde\eta_i^2+\tilde\xi_i^2}{2}$, then, for $i=1$, $2$,
$$|\tilde\L_i(t)-\tilde\L_i(0)|\le {\d^{\bar b}} \ ,\ |\tilde t_i(t)-\tilde t_i(0)|\le{\d^{\bar b}}\quad \forall\ 0\le t\le \frac{e^{\frac{1}{\d^{\bar a}}}}{\d}\ ,$$
with $\d:=\m^{\bar d}\e$.
\end{theorem}
For part of the proof, we shall deal with the system $\cH_{\rm 3b}$ in \equ{sp3bp}, which, reduces to the system $\cH_{\rm pl3b}$ in \equ{pl3bp} when $p_1=q_1=0$.

\vskip.1in
\noi
\proof {\bf Step 0.} Let us denote again as $\breve\phi$ a suitable symplectic transformation, whose existence is guaranteed by \cite{pinzari-th09}--\cite{chierchiaPi11b}, that conjugates $\cH_{\rm 3b}$ to a Hamiltonian $\breve\cH_{\rm 3b}$ having the same form as the one in \equ{3b}--\equ{remainder}, but with  $\breve N+\tilde N$ in Birkhoff normal form up to order $2m$, with possibly smaller $\cA$ of the form of  \equ{cA}, $\e_0$.  In the domain \equ{planar domain}, $\breve\phi$ is $\e^2$--close to the identity.
 
\subsection{Step 1: The Birkhoff normal form of order six}\label{Some dynamical application}
 In this section, we  aim to compute the Birkhoff normal form of order six if the three--body problem (planar and spatial).
   
\nl 
Let
 \beq{breve bir coord}\breve u_i:=\frac{\breve  \eta_i-{\rm i} \breve \xi_i}{\sqrt2}\ ,\quad  \breve u_i^\star:=\frac{ \breve \eta_i+{\rm i}\breve  \xi_i}{\sqrt2{\rm i}}\ ,\quad \breve v:=\frac{ \breve p_1-{\rm i} \breve q_1}{\sqrt2}\ ,\quad\breve v^\star:=\frac{ \breve p_1+\breve {\rm i} q_1}{\sqrt2{\rm i}}\ .\eeq
We shall show that, if $t_1:={\rm i}\breve u_1\breve { u}_1^\star$, $t_2:={\rm i}\breve u_2\breve {u}_2^\star$, $t_3:={\rm i}\breve v \breve v^\star$, 
\begin{claim}\label{BNF3} The Birkhoff normal form of order six of $(f_{\rm 3b})_{\rm av}$ is given by \equ{normal form of order 6}.
\end{claim}

\nl
Note that the $\big(1+{\rm O}(\frac{\L_1}{\L_2})\big)$--factor in \equ{normal form of order 6} has not been written for simplicity (it is available from below).

\proof 
By Claim \ref{integration}, the proof of \equ{normal form of order 6} amounts to compute the Birkhoff normal form of order six of $N$ in \equ{contributions}, up to an error of order $\frac{a_1^3}{a_2^4}$.  The constant term $f^{(0)}_{12}$ in \equ{contributions} contributes with $-\frac{\bar m_1\bar m_2}{a_2}$ to \equ{normal form of order 6}.
We  check that the Birkhoff normal form of $f^\ppd_{12}$  is corresponds to what remains in \equ{normal form of order 6}.
Recalling the definition of  $f_{12}^\ppd$ in \equ{a-exp neew} and the formulae in \equ{sum}, \equ{averaged planar} and \equ{averaged spatial3b}, we have that the explicit formula of  \equ{2aver**} in terms of {\sc rps} variables is
\beqa{Pav*}
f_{12}^\ppd
&=&\frac{a_1^2}{4a_2^3}\Big(1+3{\rm i}u_1 u_1^\star\bar e_1^2-3{\rm i}v v^\star\bar{\eufm s}^2-9({\rm i} u_1 u_1^\star)({\rm i} v v^\star)\bar{\eufm s}^2\bar e_1^2\nonumber\\
&&-\frac{15}{2}\big(( u_1^\star)^2 v^2+(v^\star)^2  u_1^2\big){ \bar e_1^2\bar{\eufm s}^2}\Big){\rm f}
\eeqa

\nl
where 
$\bar e_1$, $\bar{\eufm s}$, ${\rm f}$ are suitable functions of ${\rm i}u_1 u_1^\star$, ${\rm i}u_2 u_2^\star$ and ${\rm i}v v^\star$ (see Appendix  \ref{formula for Theta} for more details). 
Here we shall need only the first terms of their respective Taylor  expansions, which are 
\beqa{frequencies}
\bar e_1^2&=&\frac{1}{\L_1}-\frac{{\rm i}u_1 u_1^\star}{2\L_1^2}\nonumber\\
\bar{\eufm s}^2&=&\frac{1}{\L_1}+\frac{1}{\L_2}
+\frac{{\rm i}u_1u_1^\star}{\L_1^2} + \frac{{\rm i}u_2u_2^\star}{\L_2^2} - (\frac{1}{4 \L_1^2} + \frac{1}{4 \L_2^2} + \frac{1}{\L_1 \L_2}){\rm i}vv^\star
\nonumber\\
&+&\frac{1}{ \L_1^3} ({\rm i}u_1u_1^\star)^2 + \frac{1}{\L_2^3} ({\rm i}u_2u_2^\star)^2 
-(\frac{1}{\L_1^2\L_2}+\frac{1}{2\L_1^3}) ({\rm i}u_1u_1^\star)({\rm i}vv^\star) 
\nonumber\\
 &-&(\frac{1}{\L_1\L_2^2}+\frac{1}{2\L_2^3}) 
 ({\rm i}u_2u_2^\star)({\rm i}vv^\star)+ (\frac{1}{4\L_1\L_2^2}+\frac{1}{4\L_1^2\L_2}) ({\rm i}vv^\star)^2+\cdots\nonumber\\
{\rm f}&=&1+3\frac{{\rm i}u_2u_2^\star}{\L_2}+6(\frac{{\rm i}u_2u_2^\star}{\L_2})^2+10(\frac{{\rm i}u_2u_2^\star}{\L_2})^3+\cdots
\eeqa

\nl
Since  
$f^\ppd_{12}$ depends on $(u_2,u_2^\star)$ only via ${\rm i}u_2 u_2^\star$, this ``action'' (besides being preserved by the transformation $\breve\psi$ in \equ{hat psi}) is also preserved at any step of Birkhoff normalization. Since the factor ${\rm f}$ in \equ{Pav*} depends only on ${\rm i}u_2 u_2^\star$  (see Appendix \ref{formula for Theta}), we may leave such factor aside and look separately at the term inside parentheses 
\beqano
{\rm F}&:=&1+3{\rm i}u_1 u_1^\star\bar e_1^2-3{\rm i}v v^\star\bar{\eufm s}^2-9({\rm i} u_1 u_1^\star)({\rm i} v v^\star)\bar{\eufm s}^2\bar e_1^2-\frac{15}{2}\big(( u_1^\star)^2 v^2+(v^\star)^2  u_1^2\big){ \bar e_1^2\bar{\eufm s}^2}\ .
\eeqano
Using this expression and \equ{frequencies},  we see that the coefficients of ${\rm i}u_1u_1^\star$ and ${\rm i}v v^\star$ (``first order Birkhoff invariants''), are, respectively, given by\footnote{Note that we do not need to assume non--resonance of $(\O_{u_1}, \O_v)$ since $N$ in \equ{splitting} is integrable.}
$$\O_{u_1}=\frac{3}{\L_1}\ ,\quad \O_{v}=-3(\frac{1}{\L_1}+\frac{1}{\L_2})\ .$$
Letting 
$$f:=-\frac{15}{2}\big(( u_1^\star)^2 v^2+(v^\star)^2  u_1^2\big){ \bar e_1^2\bar{\eufm s}^2}\ ,\quad \phi:=-\frac{15}{2}\frac{1}{2{\rm i}(\O_{u_1}-\O_v)}\big(( u_1^\star)^2 v^2-(v^\star)^2  u_1^2\big){ \bar e_1^2\bar{\eufm s}^2}\ ,$$
one sees that the first step of Birkhoff normalization is obtained transforming ${\rm F}$ with the time--one flow of $\phi$. Then ${\rm F}$
is transformed into
\beqano
{\rm F}_1&:=&1+3{\rm i}u_1 u_1^\star\bar e_1^2-3{\rm i}v v^\star\bar{\eufm s}^2-9({\rm i} u_1 u_1^\star)({\rm i} v v^\star)\bar{\eufm s}^2\bar e_1^2+\frac{1}{2}\{\phi, f\}+{\rm o}(6)\ .
\eeqano
where ${\rm o}(6)$ stands for an expression starting with degree seven in $(u_1,v,u_1^\star,v^\star)$. The Birkhoff normal form of order six of ${\rm F}$,  obtained with a further step of Birkhoff normalization, is then
 \beqa{F2}
{\rm F}_2&:=&1+3{\rm i}u_1 u_1^\star\bar e_1^2-3{\rm i}v v^\star\bar{\eufm s}^2-9({\rm i} u_1 u_1^\star)({\rm i} v v^\star)\bar{\eufm s}^2\bar e_1^2+\frac{1}{2}\P\{\phi, f\}+{\rm o}(6)\ .
\eeqa
where $\frac{1}{2}\P\{\phi, f\}$ is obtained picking up normal  terms\footnote{Ie, monomials of the form $({\rm i}u_1{u_1^\star})^{\a} ({\rm i}vv^\star)^{\b}$.}
 of $\frac{1}{2}\{\phi, f\}$. But, 
\beqa{cubic terms}
\frac{1}{2}\P\{\phi, f\}&=&\frac{225}{2}\frac{1}{(\O_{u_1}-\O_v)}
(({\rm i}u_1 u_1^\star)({\rm i}v v)^2-({\rm i}u_1 u_1^\star)^2({\rm i}v v))
\bar{\eufm s}^4\bar e_1^4
\eeqa
where it is enough to replace $\bar{\eufm s}$, $\bar e_1$ with their respective lowest order terms in \equ{frequencies}.

\nl
In view of \equ{Pav*}, \equ{frequencies}, \equ{F2} and \equ{cubic terms}, we have that \equ{normal form of order 6} follows.\qed
\subsection{Step 2: Full reduction of the SO(3)--symmetry}\label{sec: reduction}
The next step is to reduce completely the SO(3)--symmetry from the system $\breve\cH_{\rm 3b}$. Recall the definition of $\cA$ in \equ{cA}, $\e_0$ as in Claim \ref{integration}.

\nl
Since the procedure we follow is analogue\footnote{The formulae in  \cite[\S 9]{chierchiaPi11b} are a bit different (but obviously, equivalent) from \equ{reduction**}, since in  \cite[\S 9]{chierchiaPi11b} we reduce the last couple of variables, denoted as \cite[$(\breve p_{n-1}, \breve q_{n-1})$]{chierchiaPi11b}  (corresponding to  $(\breve p_1,\breve q_1)$ in our case), while in \equ{reduction**}, we reduce the first couple. This different choice has two reasons: (i) it provides simultaneously reduction in the planar and the spatial problem and (ii) formulae are a bit simpler, since the term $t_1^3$ does not appear in \equ{normal form of order 6}.} to the one  in \cite[\S 9]{chierchiaPi11b}, we shall skip some detail and refer to \cite[\S 9]{chierchiaPi11b} for complete information. We switch to a new set of symplectic variables $(\L_1,\L_2,G, \hat u_2,\hat u_2,\hat\l_1,\hat\l_2, \hat g,\hat u_2^\star, \hat u_3^\star)$ defined via\footnote{Analogue transformations were considered in \cite{maligeRL02}.} 
\beqa{reduction**}
\hat\phi:\ \ &&\arr{\L_i=\L_i\\ \\
\breve \l_i=\hat\l_i+\hat g
}\ \ \arr{\dst {\breve u}_2=\hat u_2 e^{{\rm i}\hat g}\\ \\
\dst {\breve  u}^\star_2=\hat u_2^\star e^{-{\rm i}\hat g}
}\ \ \arr{\dst {\breve  v}_2=\hat v_2 e^{{\rm i}\hat g}\\ \\
\dst {\breve v}^\star_2=\hat v_2^\star e^{-{\rm i}\hat g}
}\nonumber\\  &&\arr{
\breve u_1=\sqrt{\varrho^2/2-\hat t_2-\hat t_3}e^{{\rm i}\hat g}\\ \\
\breve u_1^\star=-{\rm i}\sqrt{\varrho^2/2-\hat t_2-\hat t_3}e^{-{\rm i}\hat g}
}
\eeqa
with $\breve u_1$, $\breve u_2$, $\breve v$, $\breve u_1^\star$, $\breve u_2^\star$, $\breve v^\star$ defined as in \equ{breve bir coord}, $\varrho^2/2:=\L_1+\L_2-G$, $\hat t_2:={\rm i}\hat u_2\hat u_2^\star$, $\hat t_3:={\rm i}\hat v\hat v^\star$. From the last couple of definitions, one sees that $G$ is just the function in\footnote{As discussed in \cite[Proposition 7.3]{chierchiaPi11b} any step of Birkhoff normalization commutes with $\bar\cR_g$ in \equ{Rg}, the the time--$g$ flow of $G$ in \equ{G**}; equivalently, it preserves $G$.} \equ{G**} (with $n=2$) and hence its conjugated angle, $\hat g$, is cyclic in the system. 
Letting $(\hat\eta_2,\hat\xi_2)$, $(\hat p_1,\hat q_1)$ the real variables associated, respectively,  to $(u_2,u_2^\star)$, $(v,v^\star)$ via \equ{Bir coord} and $\hat z:=(\hat\eta_2,\hat p_1,\hat\xi_2,\hat q_1)$. Fix $\varrho_\star<\e_0$.
There follows from \cite[Remark 9.1-(iv)]{chierchiaPi11b} that $\hat\phi$ is well defined and symplectic in the domain defined by $(\l_1,\l_2,\hat g)\in \torus^3$ and
$$G\in \real\ ,\ (\L_1,\L_2)\in \cA_G:=\{(\L_1,\L_2)\in \cA:\ 0<\varrho_\star\le\varrho(\L,G)<\e_0\}\ ,\quad |\hat z|<\varrho_\star\ .$$
 As usual, being $\hat g$ cyclic, we regard
$G$  as an external fixed parameter so as to 
have a reduced (four--dimensional) phase space for the variables $(\L,\hat\l,\hat z)$.

\vskip.1in
\noi
Let 
\beq{reduced system1}\hat\cH_G:=\breve\cH_{\rm 3b}\circ\hat\phi=h_{\rm Kep}+\m\hat f_{G}(\L,\hat\l,\hat z)\eeq
denote the fully reduced system (where $\breve\cH_{\rm 3b}$ is as in Claim \ref{integration}) on the phase space
\beq{reduced system2}\hat\cM_G^8:=\cA_G\times\torus^2\times B^4_{\varrho_\star}\eeq
We may assume that the function $\hat{N}+\check{N}$, where $\hat{N}:=\breve {N}\circ\hat\phi$ and $\check N:=\tilde N\circ\hat\phi$, is again in Birkhoff normal form of order $2m$. If not, proceeding as in  \cite[Proof of Proposition 5.1]{chierchiaPi11c}, one can find a $\e^{2m+1}$--close to the identity symplectic transformation $\check\phi$ such that $\hat{N}'+\check{N}":=(\hat{N}+\check{N})\circ\check\phi$ is so. In the following statement, replace  eventually $\hat\phi$, $\hat N$ and $\check N$ with, respectively, $\hat\phi':=\hat\phi\circ\check\phi$, $\hat N'$, $\check N'$.
\begin{proposition}\label{prop: reduced normalized}
The system \equ{reduced system1}--\equ{reduced system2} verifies 
\beqno
(\hat f_G)_{\rm av}=\hat{N}+\check N\eeqno
 where $\hat{N}+\check N$ is in Birkhoff normal form of order $2m$, $|\check N|\le \const\a^3$. Moreover, the first three orders of $\hat N$ are given by
\beqa{normal form of order 6 reduced} 
\hat{N}&:=&-\frac{\bar m_1\bar m_2}{a_2}-\bar m_1\bar m_2\frac{a_1^2}{4a_2^3}\Big(\big(1
-3(\frac{1}{\L_1}-\frac{1}{\L_2}){\hat t_2}-3(\frac{2}{\L_1}+\frac{1}{\L_2})\hat t_3\Big)\nonumber\\
&-&\bar m_1\bar m_2\frac{a_1^2}{4a_2^3}\Big(-\frac{3}{2}\frac{\hat t_2^2}{\L_1^2}+9\frac{\hat t_2\hat t_3}{\L_1^2}
+12\frac{\hat t_3^2}{\L_1^2}\nonumber\\
&-&\frac{9}{2}\frac{\hat t_2^3}{\L_1^2\L_2}-\frac{105}{4}\frac{\hat t_2^3\hat t_3}{\L_1^3}-\frac{315}{4}\frac{\hat t_2\hat t_3^2}{\L_1^3}-\frac{105}{2}\frac{\hat t_3^3}{\L_1^3}
\big)\big(1+{\rm O}(\frac{\L_1}{\L_2})+{\rm O}(\varrho^2)\big)\nonumber\\
&+&{\rm O}(|t|^{7/2})\Big)\ .
\eeqa
 \end{proposition}
\proof
The term $\hat N$ is easily computed from \equ{normal form of order 6} and \equ{reduction**}, which amounts to replace, in \equ{normal form of order 6}
$$t_1:=\frac{\varrho^2}{2}-\hat t_2-\hat t_3\ ,\quad t_2=\hat t_2\ ,\quad t_3=\hat t_3\ .$$
We then find \equ{normal form of order 6 reduced}. \qed
\subsection{Step 3: Averaging fast angles}
In the next step we introduce, on a suitable phase space
\beq{reduced averaged phsp}\ovl\cM^8_G:=\bar D\times\torus^2\times{B}^4_{\e_1/4}\subset \hat\cM_G^8\ ,\eeq
(where $\hat\cM_G^8$ is as in \equ{reduced system2}; $\e_1\le \varrho_\star$ will be arbitrary) a new system
\beq{averaged reduced}\ovl\cH_{G}:=h_{\rm Kep}(\ovl\L)+\m (\hat{ N}(\ovl\L,\ovl z)+
\hat N_\star(\ovl\L,\ovl z))+\m \ovl f_{G}(\ovl\L,\ovl\l,\ovl z)\eeq
where $\hat N$ 
is
as in the previous sections, $\hat N_\star$ (as well as $\hat N$) depends only on $\ovl t_1={\rm i}\ovl u_1 \ovl u_1^\star$, $\ovl t_2={\rm i}\ovl u_2 \ovl u_2^\star$, $\ovl t_3={\rm i}\ovl v \ovl v^\star$ and is suitably small and $\ovl f_G$ is suitably small. 
\begin{lemma}\label{aver3b}
There exist positive numbers $\bar M$ $\r_0$, $s_0$,    depending only of $h_{\rm Kep}$ and $f_{\rm 3b}$ in \equ{rps ham} such that, for any given $m\in \natural$, one can find  $\g_\star$, $\a_*$, $\m_\star$, $C$ (depending only on $m$, $\e_0$, $s_0$) such that for any $\m$, $\a$, $\bar\g>0$, $\t>2$, $\bar K>\frac{6}{s_0}$, 
verifying $0<\a<\a_*$, $0<\m<\m_\star$,
\beq{gamma new}\bar\g\ge\g_\star\max\{\sqrt\m\bar K^{\t+1},\ \sqrt[3]{\m\e_1}\bar K^{\t+1}\}\ ,\quad \bar\r:=\frac{\bar\g}{2\bar M\bar K^{\t+1}}\le \r_0\ ,\eeq
an open set $\bar D\subset \cA_G$ 
with 
$$\meas\Big(\cA_G\setminus{\bar D}\Big)\leq C\bar\g\meas\cA_G$$
defined by the following inequalities for the Keplerian frequencies $\o_{\rm Kep}:=\partial_{\ovl\L}h_{\rm Kep}$ 
$$|\o_{\rm Kep}\cdot k|\ge \frac{\bar\g}{\bar M\bar K^\t}\quad \forall k:\ 0<|k|_1\le \bar K$$
such that for any positive number $\e_1\le \varrho_\star$ a real--analytic transformation\footnote{We refer to \cite{poschel93} for (now, standard)  notations of the kind $\cA_\r$, or $\torus^n_s$, where $\cA$ is a subset of the reals and $\r$, $s$ are positive numbers.}
$$\ovl\phi:\quad (\ovl\L,\ovl\l,\ovl z)\in \bar D_{\bar\r/16}\times\torus^2_{s_0/48}\times{B}^4_{\e_1/4}\to (\L,\hat\l,\hat z)\in(\cA_G)_{\r_0}\times\torus^2_{s_0}\times B^4_{\varrho_\star}$$
exists, which is 
$\{\frac{\m\bar K^{2(\t+1)}}{\bar\g^2},\ \frac{\m\e_1\bar K^{3(\t+1)}}{\bar\g^3}\}$
-close to the identity and lets the Hamiltonian \equ{reduced system1}--\equ{reduced system2}  into $\ovl\cH_{G}:=\hat\cH_G\circ\ovl\phi$ as in \equ{averaged reduced} with $\hat N$ as in Proposition \ref{prop: reduced normalized}, $\hat N_\star$ in Birkhoff normal form of order $m$, with Birkhoff invariants $\frac{\m\bar K^{2\t+1}}{\bar\g^2}$--close to $0$ and
\beq{averaged reduced1}| \ovl f_{G}|\le C\m \max\{e^{-\bar Ks_0/6},\ \e_1^{2m+1}\}\ .\eeq
\end{lemma}
The proof (sketched below) of Lemma \ref{aver3b} relies on Normal Form (Averaging\footnote{Sometimes distinction between ``Normal Form'' and ``Averaging'' Theory is made, depending on the strength of the remainder. For an exponentially small remainder, as in \cite{nehorosev77}, \cite{poschel93}, \cite{biascoCV03}, ``Normal Form'' Theory is often used (after \cite{poschel93}); for a quadratically--small remainder, ``Averaging'' Theory is used, after \cite{arnold63}. Normal form Theory is obtained with suitably many steps of averaging. 
}) Theory for properly--degenerate systems and the classical Birkhoff theory (see, \eg, \cite{hoferZ94}). As for Normal form theory, we refer to the theory developed in \cite{biascoCV03} (see also \cite{chierchiaPi10}), which, in turn, generalizes ideas and techniques of \cite{poschel93} to the degenerate case. For information on Normal Form theory, see \cite{arnold63}, \cite{nehorosev77}, \cite{poschel93}, \cite{biascoCV03}, \cite{chierchiaPi10} and references therein.

\vskip.1in
\noi
{\bf Sketch of proof  of Lemma \ref{aver3b}} 
We use analogue techniques as the ones  in  \cite{chierchiaPi10}, therefore, we shall limit to describe the necessary changes. We refer, in particular, to  \cite[{ Steps 1--4} in the proof of Theorem 1.4]{chierchiaPi10}.
First of all, choice, in \cite[{ Steps 1--4} in the proof of Theorem 1.4]{chierchiaPi10},
\beqano
&&n_1=2\ ,\quad n_2=2\ ,\quad V=\cA_G\ ,\quad \k=\a^3\ ,\quad \e_0:=\varrho_\star\ ,\quad H:=\hat\cH_G\nonumber\\
&&h=h_{\rm Kep}\ ,\quad P_{\rm av}=\hat N+\check N\ ,\quad P:=\m \hat f_G\ ,\quad\nonumber\\
&& I=(\L_1,\L_2)\ ,\quad \varphi:=(\hat\l_1,\hat\l_2)\ ,\quad p:=(\hat\eta_2, \hat p_1)\ ,\quad q:=(\hat\xi_2,\hat q_1)\nonumber\\
&& \O:=\frac{3}{4}\bar m_1\bar m_2\frac{a_1^2}{a_2^3\L_1}(\frac{1}{\L_1}-\frac{1}{\L_2},\frac{2}{\L_1}+\frac{1}{\L_2}\ )+{\rm O}(\frac{a_1^3}{a_2^4})\ .
\eeqano
Next, modify \cite[{ Steps 1--4} in the proof of Theorem 1.4]{chierchiaPi10} as follows.

\nl
In \cite[Step 1]{chierchiaPi10}, neglect \cite[Eq. (36)]{chierchiaPi10}, so as to ``leave $\bar K$ free'' and hence replace $\log\e^{-1}$ with $\frac{s_0}{30}\bar K$ wherever it appears (\ie, \cite[Eqs. (41), (42), (43)]{chierchiaPi10}). Neglect the second line in \cite[Eq. (40)]{chierchiaPi10}. At the end of \cite[Step 1, 2, 3, 4]{chierchiaPi10}, in the definition of $\bar H$, $\tilde H$, $\check H$, $\breve H$, respectively, replace $\e^5$ with $e^{-{\bar K}s_0/6}$.  At the beginning of \cite[Step 2, 3, 4]{chierchiaPi10}, in the definition of, respectively, $\tilde v$, $\hat v$, $\check v$, replace $\e$ with $\e_1\le \e_0$. 
In \cite[Step 2]{chierchiaPi10} replace  ``$\bar N$ also has a $\dst {\m(\log{\e^{-1}})^{2\t+1}}{\bar\g^{-2}}$--close--to--$0$ elliptic equilibrium point'' with ``$\bar N$ also has a $\dst{\m\bar K^{2\t+1}}{\bar\g^{-2}}$--close--to--$0$ elliptic equilibrium point''.
Replace\footnote{In \cite[Eq (45)]{chierchiaPi10} $\m\e$ should be replaced by $\m$. This does not affect  the thesis of \cite[Theorem 1.4]{chierchiaPi10}}  \cite[Eqs. (43), (44), (45), (46)]{chierchiaPi10} with, respectively: (43)': $|\bar p-\tilde p|$, $|\bar q-\tilde q|$ $\leq C\frac{\m\bar K^{2\t+1}}{\bar\g^2}$, $ |\bar\varphi-\tilde\varphi|\leq C\max\Big\{\frac{\e_1^2\bar K^{\t+1}}{\bar\g},\frac{\m \e_1\bar K^{3\t+2}}{\bar\g^3},\Big\}$; (44)': $ |\tilde p-\hat p|,\ |\tilde q-\hat q|\leq C\max\{\frac{\m \e_1\bar K^{2\t+1}}{\bar\g^{2}}\}\ ,\quad |\tilde\varphi-\hat\varphi|\leq 
C\max\{\frac{\m \e_1^2\bar K^{3\t+2}}{\bar\g^3}\}
$;  (45)': $ |\hat\O-\O|,\quad |\hat R|\leq C\frac{\m \bar K^{2\t+1}}{\bar\g^{2}}$ and (46)': $ |\hat p-\check p|$, $|\hat q-\check q|\leq C\frac{\m \e_1^2\bar K^{2\t+1}}{\bar\g^{2}}$, $|\hat\varphi-\check\varphi|\leq C\frac{\m \e_1^3\bar K^{3\t+2}}{\bar\g^3}$, by suitably modifying the proofs below. Moreover replace Equation just before \cite[Eq. 45]{chierchiaPi10} with\footnote{The symbol $\hat N$ used in \cite{chierchiaPi10} is here replaced with $\hat{\cal N}$, to avoid confusions with \equ{normal form of order 6 reduced}. } $\hat{\cal N}(I,p,q):=\tilde N\circ\hat\phi=\hat N+\hat R$ (where $\hat N$ is as in \equ{normal form of order 6 reduced}) and replace the last line in  \cite[Step 4]{chierchiaPi10} with ``where $\breve N(\breve I,\breve r)$ is a polynomial of degree $m$ in $(I_{n_1+1}, \cdots, I_{n_2})$''. Lemma \ref{aver3b} follows, with $\hat N_\star:=\breve N-\hat N$ and $\hat f_G:=\m (e^{-\bar Ks_0/6}\breve P+{\rm O}(\e_1^{2m+1}))$.\qed
We then apply Lemma \ref{aver3b}  to the system \equ{reduced system1}--\equ{reduced system2} with $\bar K$ as in \equ{Kappa}, with $\e$, $\a$ replaced by $\e_1$, $\a_*$ and
\beqa{gamma}
\bar\g=\g_\star\sqrt[4]\m\bar K^{\t+1}\ .\eeqa
where $\g_\star$ is as in \equ{gamma new}.  By the thesis of Lemma \ref{aver3b}, we conjugate $\hat\cH_G$ in \equ{reduced system1}--\equ{reduced system2} to $\ovl\cH_G$ in \equ{reduced averaged phsp}--\equ{averaged reduced}, with $\ovl f_G$ satisfying \equ{averaged reduced1}, via a symplectic transformation which, by the choice of $\bar\g$  in  \equ{gamma}, is $\m^{1/12}$--close to the identity.
\subsection{Step 4: Nehoro{\v{s}}ev Theory}\label{steepness}
We apply Nehoro{\v{s}}ev Theory (\ie, Theorem \ref{Nek thm}) to the system $\ovl\cH_G$ in \equ{reduced averaged phsp}--\equ{averaged reduced}, in the {\sl planar} case, \ie, with $\hat t_4=0$.

\nl
 For information on the tools that are used, compare \cite{nehorosev73}, \cite{nehorosev77}, \cite{nehorosev79} and Appendix \ref{Nekhorossev}.
 
 \nl
In applying Theorem \ref{Nek thm}, we shall take
\beqa{system}
&&n_1=3\ ,\ n_2=0\ ,\ V=\bar D\ ,\ B^4:=B^4_{\e_1/8}\ ,\ \r:=\min\{\bar\r/16, s_0/48, \e_1/8\}
\nonumber\\
&&H_0(\ovl\L_1,\ovl\L_2,\ovl t_1):=h_{\rm Kep}(\ovl\L_1,\ovl\L_2)+\m(\hat N+\hat N_\star)(\ovl\L_1,\ovl\L_2,\ovl t_1)\ ,\  P:=\m\ovl f_G\eeqa
where $\bar\r$, $s_0$ and $\e_1$ are as in Lemma \ref{aver3b}
\vskip.1in
\noi
 We have to check\footnote{Recall that,  for $n_2=0$, as it is in our case, condition \equ{inclusion} is void; see Appendix \ref{Nekhorossev}.} steepness of $H_0(\ovl\L_1,\ovl\L_2,\ovl t_1)$ and the smallness condition \equ{Nek smallness} of $P$.  The first check is provided by the following claim.

\begin{claim}\label{clm: steepness}
The  function $H_0$ in \equ{system}
is  $( g, m,  C_1, C_2, {\eufm a}_1, {\eufm a}_2, \d_1, \d_2)$--{\sl steep}, with 
\beq{steepness parameters}g=\hat g\ ,\quad m=\hat m\ ,\quad  {\eufm a}_i=\hat{\eufm a}_i\ ,\quad \d_i=\min\{\sqrt{\a_*},\e_1^2\}\hat\d_i\ ,\quad C_i=\m\a_*^2\hat C_i\eeq
where $( \hat g, \hat m,  \hat C_1, \hat C_2, \hat {\eufm a}_1, \hat {\eufm a}_2, \hat \d_1, \hat \d_2)$ suitable numbers independent of $\a_*$, $\m$, $\e_1$.
\end{claim}

\nl

\proof We take, in \equ{normal form of order 6 reduced}, $\hat t_3=0$. The system has three degrees of freedom. 
We firstly prove {\sl steepness} for a suitable ``rescaled'' system associated to ${\rm F}$. That is, if $\hat N_0:=-\frac{\bar m_1\bar m_2}{a_2}$  is as  in \equ{normal form of order 6 reduced} and $\hat N_1:=\hat N-\hat N_0$ we consider the system
\beqa{rescaled}{\rm F}_{\rm resc}(\hat\L_1,\hat\L_2,\hat t_2)&:=&\bar m_1^2{\bar m_0\a_*}\Big(h_{\rm Kep}^\ppu(\bar m_1\sqrt{\bar m_0\a_*}\hat\L_1)+\b_2h_{\rm Kep}^\ppd(\hat\L_2)+\m \hat N_0(\L_2)\nonumber\\
&+&\m\b_3(\hat N_1+ \hat N_\star)(\bar m_1\sqrt{\bar m_0\a_*}\hat\L_1,\hat\L_2,\e_1^2\hat t_2)\Big)\eeqa
with $\a_*$, $\e_1$ as in Lemma \ref{aver3b}
\beq{choices of parameters}\b_2:=\a_*^{-3/2}\ ,\quad \b_3:=\m^{-1}\a_*^{-3}\e_1^{-2}\ .\eeq
We check that ${\rm F}_{\rm resc}$ is {\sl steep} by verifying the three--jet condition: See Appendix \ref{Steepness conditions}.
The three--jet condition \equ{three jet} for the system \equ{rescaled} is

\beqa{three jet1}
\arr{\dst\eta_1+\b_2\a_*^{3/2}(\frac{\hat a_1}{\hat a_2})^{3/2}\eta_2+ \b_3 \a_*^3\e_1^2\m\frac{3}{4}\frac{\bar m_2}{\bar m_0}(\frac{\hat a_1}{\hat a_2})^3\eta_3=0\\
\dst\dst\eta_1^2+\frac{\bar m_1}{\bar m_2}\b_2\a_*^2(\frac{\hat a_1}{\hat a_2})^{2}\eta_2^2- \b_3 \a_*^3\m\e_1^4{\frac{1}{4}}\frac{\bar m_2}{\bar m_0}(\frac{\hat a_1}{\hat a_2})^3\eta_3^2=0\\
\dst\dst\eta_1^3+(\frac{\bar m_1}{\bar m_2})^2\b_2\a_*^{5/2}(\frac{\hat a_1}{\hat a_2})^{5/2}\eta^3_2+ \b_3 \a_*^{7/2}\m\e_1^6{\frac{9}{16}}\frac{\bar m_1}{\bar m_0}(\frac{\hat a_1}{\hat a_2})^{7/2}\eta_3^3=0
}
\eeqa
where we have used $m_i=\bar m_i+{\rm O}(\m)$, $M_i=\bar m_0+{\rm O}(\m)$ and neglected higher order terms going to zero with $\m$, $\e_1$, $\a_*$. If we
eliminate $\eta_1$ from the first and the second equation and from the first and the third equation,  we obtain a homogeneous system of two equations in $(\eta_2,\eta_3)$ that, in view of \equ{choices of parameters}, generically, the has only solution $\eta_2=\eta_3=0$, implying that also $\eta_1=0$. This implies that the function ${\rm F}_{\rm resc}$ \equ{rescaled} is $(2\hat g,\hat m/2, \hat C_1,\hat C_2,\hat {\eufm a}_1,\hat {\eufm a}_2,\hat \d_1,\hat \d_2)$--{\sl steep} with suitable values  of $(\hat g, \hat m, \hat C_1,\hat C_2,\hat {\eufm a}_1,\hat {\eufm a}_2,\hat \d_1,\hat \d_2)$ which are of order $1$ in $\m$, $\a_*$, $\e_1$. This readily implies that ${\rm F}$ in \equ{system} is $( g, m,  C_1, C_2, {\eufm a}_1, {\eufm a}_2, \d_1, \d_2)$--{\sl steep}, with $( g, m,  C_1, C_2, {\eufm a}_1, {\eufm a}_2, \d_1, \d_2)$ as in \equ{steepness parameters}.\qed
\begin{remark}\label{steepness in space}\rm
In the case of the {\sl spatial} three--body  problem, instead of \equ{three jet1},
we would have 
\beqano
\arr{\dst\eta_1+\b_2\a_*^{3/2}(\frac{\hat a_1}{\hat a_2})^{3/2}\eta_2+ \b_3\e_1^2 \a_*^3\m\frac{3}{4}\frac{\bar m_2}{\bar m_0}(\frac{\hat a_1}{\hat a_2})^3\eta_3+ \b_3\e_1^2 \a_*^3\m\frac{3}{2}\frac{\bar m_2}{\bar m_0}(\frac{\hat a_1}{\hat a_2})^3\eta_4=0\\
\dst\dst\eta_1^2+\frac{\bar m_1}{\bar m_2}\b_2\a_*^2(\frac{\hat a_1}{\hat a_2})^{2}\eta_2^2- \b_3\e_1^4 \a_*^3\m{\frac{1}{4}}\frac{\bar m_2}{\bar m_0}(\frac{\hat a_1}{\hat a_2})^3\eta_3^2+\b_3\e_1^4 \a_*^3\m\frac{3}{2}\frac{\bar m_2}{\bar m_0}(\frac{\hat a_1}{\hat a_2})^3\eta_3\eta_4\\
\dst\qquad\qquad+{2}\b_3\e_1^4 \a_*^3\m\frac{\bar m_2}{\bar m_0}(\frac{\hat a_1}{\hat a_2})^3\eta_4^2=0\\
\dst\dst\eta_1^3+(\frac{\bar m_1}{\bar m_2})^2\b_2\a_*^{5/2}(\frac{\hat a_1}{\hat a_2})^{5/2}\eta^3_2+ \b_3\e_1^6 \a_*^{7/2}\m\e_1^6{\frac{9}{16}}\frac{\bar m_1}{\bar m_0}(\frac{\hat a_1}{\hat a_2})^{7/2}\eta_3^3+\b_3\e_1^6 \a_*^3\m{\frac{70}{64}}\frac{\bar m_2}{\bar m_0}(\frac{\hat a_1}{\hat a_2})^3\eta_3^2\eta_4\\
\dst\qquad\qquad +\b_3\e_1^6 \a_*^3\m{\frac{105}{32}}\frac{\bar m_2}{\bar m_0}(\frac{\hat a_1}{\hat a_2})^3\eta_3\eta_4^2+\b_3\e_1^6 \a_*^3\m{\frac{105}{16}}\frac{\bar m_2}{\bar m_0}(\frac{\hat a_1}{\hat a_2})^3\eta_4^3=0
}
\eeqano
It is not clear to the author if this system exhibits non--trivial solutions, so the analysis of this case is deferred to a subsequent paper.
\end{remark}
We can now complete the
\proof {\bf of Theorem \ref{stability Theorem detailed}}
It remains only to check condition \equ{Nek smallness}, with $P$, $\r$ as in \equ{system}. In view of \equ{M*}, \equ{M0},  \equ{steepness parameters} and the choice of $\bar\g$ in \equ{gamma}, we have
$\r\ge \r_\star\min\{\e_1,\ \sqrt[4]\m\}$ and hence
\beqano
M_\star&\ge&\frac{\tilde c}{\r}\min\{(\m\a_*^2)^q,\r^q\}\ge\frac{\tilde c}{\r}\min\{(\m\a_*^2)^q,\e_1^q,\ \bar\r^q\} \nonumber\\
&\ge&\frac{\tilde c}{\r}\min\Big\{(\m\a_*^2)^q,\e_1^q,\ (\frac{\bar\g_\star}{2\bar M})^q\m^{q/4}\Big\}\ge \frac{c_\star}{\r}\min\Big\{(\m\a_*^2)^q,\e_1^q\Big\}
\eeqano
for some $q>1>c_\star$ depending only on $n_1$, $n_2$, ${\eufm a}_1$, $\eufm a_2$. 
Noticing that \equ{averaged reduced1} and Cauchy inequality imply
$$M:=\sup|\partial P|=\m\sup|\partial\ovl f_G|\le \tilde C\frac{\m}{\r}\max\{e^{-\bar Ks_0/6},\ \e_1^{2m+1}\}$$
one sees that condition \equ{Nek smallness} is met, provided one previously fixes, in Lemma \ref{aver3b}, $2m+1\ge q$, $\bar K$ as in \equ{Kappa}, with a suitable $\bar K_\star$ and takes $\e_1< (\frac{c_\star}{\tilde C})^{1/(2m+1)}(\m\a_*^2)^{q/(2m+1)}$.
The thesis then follows, with $\a$, $\e$ replaced by $\a_*$, $\e_1$ and $\phi:=\breve\phi\circ\hat\phi\circ\ovl\phi\circ\hat\phi^{-1}$. \qed

 \appendix
\section{The Fundamental Theorem and another result in Arnold's 1963 paper}
\label{ArnoldKAM}
\renewcommand{\theequation}{
\arabic{equation}}
Here we recall two theorems  in \cite{arnold63}. The former is named ``The Fundamental Theorem'' in  \cite{arnold63} and is as follows.

\nl
Recall the definition of $\cP_{\e_0}$ in \equ{defcP}.
\begin{theorem}[{V. I. Arnold}, {\cite[p. 143]{arnold63}}]\label{fundamental theorem}
Consider a  Hamiltonian of the form
$$H(I,\varphi,p,q)=H_0(I)+\m P(I,\varphi,p,q)$$
which is real--analytic on $\cP_{\e_0}$
where $V\subset \real^{n_1}$ is open and connected, $B^{2n_2}_{\e_0}\subset \real^{2n_2}$ is a ball of radius $\e_0$ around the origin and $\torus:=\real/(2\p\integer)$. Assume that
\begin{itemize}
\item[{\rm (i)}] $I\in V\to \partial_I H_0$ is a diffeomorphism;

\item[{\rm (ii)}]  $P_{\!\rm av}$ is in Birkhoff normal form\footnote{We refer to \cite{hoferZ94} for information on Birkhoff Theory.} of order $6$;

\item[{\rm (iii)}]  the matrix $\b$ of the ``second order Birkhoff invariants'': is not singular: $|\det \b|\ne 0$ on $V$.

\end{itemize}
  Then, there exists $\e_0>0$ such that, for \beq{A5} 0<\e<\e_0\ ,
  \qquad 0<\m<\e^8\ , \eeq one can find a set $\cK_{\m,\e}\subset
  \cP_\e\subset \cP_{\e_0}$, with $$ {\meas{\cK_{\m,\e}}}\ge(1-\e^{16(n_1+n_2)}){\meas \cP_\e}$$ formed by the union of $H$--invariant
  $(n_1+n_2)$--dimensional 
  tori on  which the $H$--motion is analytically conjugated to linear
  Diophantine\footnote{{\it I.e.}, the flow is conjugated to the
    Kronecker flow $\theta\in\torus^{n_1+n_2}\to \theta+\o\,t\in \torus^{n_1+n_2}$,
    with $\o\in \real^{n_1+n_2}$ satisfying $|\o\cdot k|\ge \g|k|_1^{-\t}$ for all $k\ne 0$, for suitable $\g$, $\t>0$.}
  quasi--periodic motions.
\end{theorem}
The latter is less general, but used in \cite{arnold63} to prove Theorem \ref{simpler planar}.
\begin{theorem}[V. I. Arnold, {\cite[p. 128]{arnold63}}]\label{simplifiedFT}
Under the same assumptions as in Theorem \ref{fundamental theorem}, but replacing {\rm (ii)}, {\rm (iii)} and \equ{A5} with
\begin{itemize}
\item[{\rm (ii)$'$}] $P_{\rm av}$ has the form
\item[] $\dst P_{\rm av}(I,p,q)=N(I,J)+\tilde N(I,p,q)$ where $J=(\frac{p_1^2+q_1^2}{2},\cdots, \frac{p_{n_2}^2+q_{n_2}^2}{2})$ and $\tilde N={\rm o}(\m)$;
\item[{\rm (iii)$'$}] the Hessian $\partial^2_{(I,J)}N$ is non--singular: $\det \partial^2_{(I,J)}N\ne 0$ on $V\times B^{2n_2}_{\e_0}$
 \end{itemize}
 and condition \equ{A5} with condition
 $$|\m|<\m_*$$
 one can find a set $\cK_{\m}\subset
  \cP_{\e_0}$, with $$ {\meas{\cK_\m}}\ge(1-\m^{a}){\meas \cP_{\e_0}}$$
(where $a$ decreases with $n_1+n_2$) having the same properties as the set $\cK_\e$ of Theorem \ref{fundamental theorem}.
\end{theorem}
\section{Proof of \equ{averaged planar}, \equ{averaged spatial3b} and \equ{Pav*}.}
\label{formula for Theta}
\renewcommand{\theequation}{
\arabic{equation}}
The formulae in \equ{averaged planar} and \equ{averaged spatial3b} are a consequence of Proposition \ref{asymptotic formula} and the formulae developed in \cite{pinzari-th09}--\cite{chierchiaPi11b} (see, \eg, \cite[Appendix A]{chierchiaPi11b}). Indeed, from such papers there results that, if
\beqano
&&a_i:=\frac{1}{M_i}(\frac{\L_i}{m_i})^2\ ,\quad e_i^2=\frac{\eta_i^2+\xi_i^2}{\L_i}-(\frac{\eta_i^2+\xi_i^2}{2\L_i})^2=:2{\rm i}u_i u_i^\star\bar e_i^2\nonumber\\
&& \zeta_i:\ \zeta_i- e_i\sin\zeta_i=\l_i+\arg(\eta_i,\xi_i)\nonumber\\
&&\bar{\eufm c}:=\frac{2\L_1+2\L_2-2{\rm i}u_1 u_1^\star-2{\rm i}u_2 u_2^\star-{\rm i}v v^\star}{4(\L_1-{\rm i}u_1 u_1^\star)(\L_2-{\rm i}u_2 u_2^\star)}\ ,\quad\bar{\eufm s}:=\sqrt{2\bar{\eufm c}(1-{\rm i} v v^\star\bar{\eufm c})}
\eeqano
then, the expressions of ${\rm C}^\ppd\cdot x^\ppu
$, $|{\rm C}^\ppd|$, ${\rm r}_1=|x^\ppu|$ and ${\rm r}_2=|x^\ppd|$ in terms of the  {\sc rps} variables    are 
\beqano
&&{\rm C}^\ppd\cdot x^\ppu
=\Big((\hat u_1 v^\star-\hat u_1^\star v)x^\ppu+{\rm i}(\hat u_1 v^\star+\hat u_1^\star v)x^\ppd\Big)
\Big)\bar{\eufm s}|{\rm C}^\ppd|\\
&&|{\rm C}^\ppd|=\L_2-{\rm i}u_2u_2^\star\ ,\quad 
{\rm r}_i=a_i(1- e_i\cos\zeta_i)\\
&&x^\ppu_1:=\frac{1}{M_1}(\frac{\L_1}{m_1})^2(\cos\zeta_1- e_1)\ ,\quad x^\ppu_2:=\frac{1}{M_1}(\frac{\L_1}{m_1})^2\sqrt{1- e_1^2}\sin\zeta_1
\eeqano
with
$$\hat u_i:=\frac{u_i}{\sqrt{{\rm i} u_i u_i^\star}}=\frac{ \eta_i-{\rm i} \xi_i}{\sqrt2\sqrt{ \eta_i^2+ \xi_i^2}}\ ,\quad \hat u_i^\star:=\frac{u_i^\star}{\sqrt{{\rm i} u_i u_i^\star}}=\frac{ \eta_i+{\rm i} \xi_i}{\sqrt2{\rm i}\sqrt{ \eta_i^2+ \xi_i^2}}\ .$$
Then we have
\beqano
({\rm C}^\ppd\cdot x^\ppu)^2
&=&\Big({\rm i}v v^\star \,
\big((x^\ppu_1)^2+(x^\ppu_2)^2\big)
+\big((\hat u_1 ^\star)^2 v^2+(v^\star)^2 \hat u_1 ^2\big)((x^\ppu_1)^2-(x^\ppu_2)^2)\nonumber\\
&+&2{\rm i}\big((\hat u_1 ^\star)^2 v^2-(v^\star)^2 \hat u_1 ^2\big)x^\ppu_1x^\ppu_2
\Big)\bar{\eufm s}^2|{\rm C}^\ppd|^2\ .
\eeqano
and hence, taking the $\ul_1$--average (recall the relation $d\ul_2=(1-e_2 \cos\zeta_2)d\zeta_2$)
\beq{last}\frac{1}{2\p}\int_\torus ({\rm C}^\ppd\cdot x^\ppu)^2d\ul_1=\Big({\rm i} v v^\star a_1^2(1+\frac{3}{2} e_1^2)+{\frac{5}{2}}\big(( u_1^\star)^2 v^2+(v^\star)^2  u_1^2\big)\frac{a_1^2 e_1^2}{2{\rm i}u_1 u_1^\star}\Big)\bar{\eufm s}^2|{\rm C}^\ppd|^2\eeq
Here, we have used 
\beqano
&&\frac{1}{2\p}\int_{\torus}\big((x^\ppu_1)^2+(x^\ppu_2)^2\big)d\l_1=\frac{1}{2\p}\int_{\torus}{\rm r}_1^2=\frac{1}{2\p}\int_{\torus}a_1^2(1- e_1\cos\zeta_1)^3d\zeta_1\nonumber\\
&&\qquad=a_1^2(1+\frac{3}{2} e_1^2)\nonumber\\
&&\frac{1}{2\p}\int_{\torus}\big((x^\ppu_1)^2-(x^\ppu_2)^2\big)d\l_1=\frac{1}{2\p}\int_{\torus}d\zeta_1\big(a_1^2(\cos2\zeta_1+ e_1^2\nonumber\\
&&\qquad+ e_1^2\sin^2\zeta_1-2e_1\cos\zeta_1)(1- e_1\cos\zeta_1)\big)=\frac{5}{2}{a_1^2 e_1^2}\nonumber\\
&&\frac{1}{2\p}\int_{\torus}x^\ppu_1x^\ppu_2d\l_1=\frac{1}{2\p}\int_{\torus}a_1^2\sqrt{1- e_1^2}(\cos\zeta_1- e_1)(1- e_1\cos\zeta_1)\sin\zeta_1 d\zeta_1=0
\eeqano
Note that \equ{last} implies \equ{averaged spatial3b}.
In turn, \equ{averaged planar} for $n=2$ and \equ{averaged spatial3b} give \equ{Pav*}, with $\dst {\rm f}:=\frac{\frac{1}{2\p}\int_{\torus}\frac{d\zeta}{(1-e_2\cos\zeta)}}{(1-\frac{\eta_2^2+\xi_2^2}{2\L_2})^2}$.
\nl

\section{Proof of Theorem \ref{thm: ArnoldConjnew}}
\label{kam}
\renewcommand{\theequation}{
\arabic{equation}}
Theorem \ref{thm: ArnoldConjnew} is an easy consequence\footnote{To obtain Theorem \ref{thm: ArnoldConjnew} from Theorem
\ref{thm:**}, it is sufficient to choose
$$\bar\g=\g_*\sqrt[4]\m\log(\a^{-1})^{\t+1}\ ,\quad \g_1= \g_2=\g_*^2\max\{\a^2,\ \sqrt[4]\m\}<\g_*\e_0^2\ ,\quad C_*:=\frac{C}{\e_0^2}\ .$$} of the following more technical statement.

\thm{thm:**} Under the same notations  and assumptions as in Theorem~\ref{thm: ArnoldConjnew}, one can find
  $\g_*$, $C_*$ such that, for any $\e_0$, one can find positive numbers $\e_1<\e_0$, $\m_*$ and  $\a_*$ such that,  for any $\a$, $\m$, $\g_1$, $\bar\g_2$, $\bar\g$ verifying 
  $$|\a|<\a_*\ ,\quad |\m|<\m_*\ ,\quad \m\sg\le \fg$$
  and 
\beqa{how small gammas}
\arr{\g_*{\sqrt[4]{\m}(\log{\a^{-1}})^{\t+1}}\le \bar\g\le \g_*\\
\g_{*} \max\big\{\a^{2},\ \frac{\sqrt\m(\log\a^{-1})^{\t+1}}{\bar\g}\big\}<\fg<\g_*\\
\g_{*}\max\big\{\a^{2}(\log{(\fg^2/\a^3)})^{\t_*+1},\\  \ 
\sqrt\m(\log\a^{-1})^{\t+1}\bar\g^{-1}\ \Big(\log\Big(\frac{\fg^2}{\m(\log\a^{-1})^{2\t+1}\bar\g^{-2}}\Big)\Big)^{\t+1}
\big\}<\sg<\g_*\e_0^2\ ,
}
\eeqa
where $\t>n:=n_1+n_2$, then, one can find a set $\cK\subset \cP$  formed by the union of $H$--invariant $n$--dimensional tori, on which the $H$--motion is analytically conjugated to  linear Diophantine quasi--periodic motions. The  set $\cK$  is of positive  measure and satisfies
\beqno
\meas \cK> \left[1-C(\bar\g+\g_1+\frac{\bar\g_2}{\e_0^2}+\a^{n_2})\right]\meas \cP_{\e_1}\ .
\eeqno
Furthermore, the flow on each $H$--invariant torus in $\cK$ is analytically conjugated to a translation $\psi\in\torus^n\to \psi + \o t\in\torus^n$ with Diophantine frequencies.
\ethm

\vskip.1in
\noi
This result  is a slight modification of \cite[Theorem 1.4]{chierchiaPi10} (which, in turn, had been obtained in \cite{pinzari-th09}).  Then here we briefly sketch its proof, describing only the necessary changes with respect to \cite[Proof of Theorem 1.4]{chierchiaPi10}  and referring the reader to that paper for more details.

\nl
To proceed, we need to recall
\begin{itemize}
\item[--] the definition of ``two velocities'' Diophantine vector\footnote{This is a suitable  generalization of the standard definition of Diophantine numbers,  introduced in \cite{arnold63}.} in \cite[Eq. (19)]{chierchiaPi10};
\item[--] the functional setting and notations described at the beginning of  \cite[\S 2]{chierchiaPi10}; 
\item[--] the ``averaging (iterative) Theorem'' \cite[Lemma A.1]{chierchiaPi10};
\item[--] the ``two--scale {\sc kam} Theorem'' \cite[Proposition 3]{chierchiaPi10}.
\end{itemize}

\nl
 {\bf Sketch  of proof of Theorem \ref{thm:**}} 
 Let $\r_0$, $s_0$, $\e_0$ (possibly with a smaller value of $\e_0$) be positive numbers such that $H$ in \equ{pndham} has analytic extension on the complex set
 $${\cal P}_{\r_0,s_0, \e_0}=V_{\r_0}\times \torus^{n_1}_{s_0}\times B^{2n_2}_{\e_0}\ .$$
Take three numbers $\bar\g$, $\g_1$, $\g_2=\m\bar\g_2$ verifying \equ{how small gammas} and $\m\bar\g_2<\g_1$,
where $\g_*$ is some large number, depending only on $n_1$, $n_2$, to be chosen below.

\nl
As in \cite[Proof of Theorem 1.4, Step 1]{chierchiaPi10}, start with removing, in $H$, the dependence on $\varphi$ up to high orders. But, at difference with \cite[Proof of Theorem 1.4, Step 1]{chierchiaPi10}, apply  \cite[Lemma A.1]{chierchiaPi10} (instead of \cite[Proposition 1]{chierchiaPi10}) , with  $\ell_1=n_1$, $\ell_2=0$, $m=n_2$ $h=H_0$, $g\equiv 0$, $f=\m P$, ${B}={B'}=\{0\}$, $r_p=r_q=\e_0$, $s=s_0$, $\r_p=\r_q=\e_0/3$, $\s=s_0/3$, 
$\L=\{0\}$,  
\beq{bar K}e^{-\bar Ks_0/3}:=\k\qquad\textrm{\ie,}\qquad \bar K=\frac{3}{s_0}\log{\k^{-1}}\ ,\eeq
$A=\bar D$,  $r=\bar\r$, $\r=\bar\r/3$, where 
$\bar D$, $\bar\r$ are defined  
as in \cite[(37)]{chierchiaPi10} 
By \cite[(38)]{chierchiaPi10}, and the choice of $\bar\g$, the following standard measure estimate holds
\beqno
\meas\Big({V}\setminus{\bar D}\Big)
\le C\g_*\sqrt{\m}(\log{\k^{-1}})^{\t+1}\meas V\eeqno
where $C$ depends on the $C^1$--norm of $H_0$. 
Proceeding as \cite[(39)]{chierchiaPi10} and the immediately following formula, one sees that the ``non--resonance'' condition \cite[(64)]{chierchiaPi10} on $\bar D_{\bar\rho}$ and the ``smallness'' condition   \cite[(65)]{chierchiaPi10} are  then verified , provided $\m$ is chosen small enough, because of the choice of $\bar\g$ and $\g_*$.
By  the thesis of \cite[Lemma A.1]{chierchiaPi10}, we find  a real--analytic symplectomorphism
$$
\bar\phi:\ (\bar I,\bar \varphi,\bar p,\bar q)\in W_{(\bar\r, \e_0)/3,s_0/3}\to (I,\varphi,p,q)\in W_{v_0,s_0}
$$
where  $W_{v_0,s_0}:={\bar D}_{\r_0}\times \torus^{\fd}_{s_0}\times {B}_{\e_0}$ ($v_0=(\r_0,\e_0)$), and,  by the choice of $\bar K$  in (\ref{bar K}),  $H$ is transformed into\footnote{
$\P_{0} T_{\bar K}  P=\dst  P_{\rm av}=\frac{1}{(2\p)^n}\int_{\torus^n}  Pd\varphi$.
}
\beqa{averaged**}
\bar H:=H\circ \bar\phi&=& h+\m P_{\rm av}+\m \bar P\nonumber\\
&=& h+\m N+\m\tilde N+\m \bar P
\eeqa
where $\dst  P_{\rm av}=N+\tilde N$ corresponds to $g_+$ of \cite[Lemma A.1]{chierchiaPi10}, $\bar P$ corresponds to $f_+$
and
hence, by the choice of $\bar K$ in \equ{bar K}, the assumption on $\tilde N$ and the thesis \cite[(68)]{chierchiaPi10} of \cite[Lemma A.1]{chierchiaPi10}, one has that the new perturbation $\m\tilde N+\m\bar P$ verifies
\beqa{smallest}\|\m\tilde N+\m\bar P\|_{v_0/3, s_0/3}&\leq& C\m\max\{\frac{\bar K^{2\t+1}}{\bar\g^2}\m,\  e^{-\bar K s_0/3},\ \k\}\nonumber\\
&\le& C\m\max\{\frac{\bar K^{2\t+1}}{\bar\g^2}\m,\  \k\}
\eeqa
In view of \cite[(69)]{chierchiaPi10}, the transformation $\bar\phi$  verifies
\beqano
|I-\bar I| ,\ |p-\bar p|,\ |q-\bar q|\leq C\frac{\m(\log{\k^{-1}})^{\t}}{\bar\g}\ ,\quad |\varphi-\bar\varphi|\leq C\frac{\m(\log{\k^{-1}})^{2\t+1}}{\bar\g^2}\ .
\eeqano

\vskip.1in
\noi
Continue as in \cite[Proof of Theorem 1.4, Step 5]{chierchiaPi10}, but replacing the set in \cite[(47)]{chierchiaPi10}. 
with the set
\beq{annulus}{\mathcal A}:=\Big\{J\in \real^{n_2}:\r_1< J_i<{\e_0^2}/9\ ,\quad 1\le i\leq\sd\Big\}\eeq
 where  $\r_1<{\e_0^2}/9$ will be fixed in the next step, on so as to maximize the measure of  preserved tori. Next define $\cD$ as in  \cite[(48)]{chierchiaPi10} (but with $\cA$ as in \equ{annulus}) and
\beq{rho}
\r:=\min\{\r_1,
 \ \bar\r/3\}\ ,\quad s:=s_0/3\eeq 
 Introduce the change of variables
$$(J,\psi)=\big((J_1,J_2),(\psi_1,\psi_2)\big)\in {\cal D}_\r\times \torus^{n_1+n_2}_s\to (\bar I,\bar\varphi,\bar p,\bar q)$$
defined as in \cite[(49)]{chierchiaPi10}, but replacing ``checks'' with ``bars'', 
This lets the Hamiltonian \equ{averaged**} into
\beqano
&&  H(J,\psi)=H_0(J_1)+\m N(J)+\m(\bar P+\tilde N)\ ,\quad (J,\psi)\in {\cal D}_\r\times \torus^{n_1+n_2}_s
\eeqano
\vskip.1in
\noindent
Next, analogously to \cite[Proof of Theorem 1.4, Step 6]{chierchiaPi10}, construct the Kolmogorov set and estimate  its measure via \cite[Proposition 3]{chierchiaPi10}.

\nl
To this end, fix $\g_1$ and $\g_2=\m\bar\g_2$, with $\g_1$, $\bar\g_2$ satisfying $\m\bar\g_2\le \g_1$ and \equ{how small gammas}. Let $\r_1$ in\equ{annulus}--\equ{rho} be chosen so that
$$\r_1=\check c_1\max\big\{\sqrt\k,\ \frac{\sqrt\m(\log\k^{-1})^{\t+1/2}}{\bar\g}\big\}\ .$$
with $\check c_1$ some large number depending only on $n_1$, $n_2$ to be fixed below. Note that the needed condition $\r_1<\e_0^2/9$ (compare the previous step; Eq. \equ{annulus}) is satisfied for
$\k<(\e_0/(3\sqrt{\check c_1}))^4$ and\footnote{Use the definition of $\bar\g$ in \equ{how small gammas}.} $\m<\g_\star^4(\e_0/(3\sqrt{\check c_1}))^8$. 
The assumption that the frequency map $\o:=\partial(H_0(J_1)+\m N(J))$ is a diffeomorphism of  $\cD_{\r}$ is trivially satisfied. Moreover,  the numbers $M$, $\hat M$, $\cdots$, $\bar M_2$ involved in \cite[Proposition 3]{chierchiaPi10} may be chosen as in \cite[Proof of Theorem 1.4, Step 6]{chierchiaPi10}, apart for $E$, which is chosen as\footnote{Compare, in particular, \equ{smallest} for the choice of $E$ and recall Equation  
 \equ{rho} and the definition of $\bar\r$ and of $\bar K$  in \equ{bar K}.}
\beqano
\pertnorm =C\max\{\m \k\ ,\ {\bar K^{2\t+1}\m^2}{\bar\g^{-2}}\}\ .
\eeqano
Then, we can take $L$ as in \cite[Proof of Theorem 1.4, Step 6]{chierchiaPi10}, while

$$
K=C\log{(E/(\m\fg^2))^{-1}}$$
and
\beqno
\begin{split}
&\hat\r\\
= & c\,\min\Big\{\frac{\fg}{(\log{(E/(\m\fg^2))^{-1}})^{\t+1}},\frac{\sg}{(\log{(E/(\m\fg^2))^{-1}})^{\t+1}} ,  \frac{\bar\g}{(\log{\k^{-1}})^{\bar\t+1}} , \r_1, 
\ \r_0\Big\}.
\end{split}
\eeqno
To check the ``{\sc kam}--smallness condition'' \cite[(32)]{chierchiaPi10}, we divide the two cases $E=C\m \k$ or $E=C\bar K^{2\t+1}\m^2\bar\g^{-2}$. If $E=\m \k$,
\beqno
\hat c\KAM\le
C \max\Big\{
\k \Big(\log\Big(\frac{\fg^2}{\k}\Big)\Big)^{2(\t+1)} \max\{ \su{\g_1^2},\su{\bar\g_2^2}\}, \frac{\k(\log{\k^{-1}})^{2(\t+1)}}{\bar\g^2}
,\  
\frac{\k}{\r_1^2},\
\frac{\k}{\r_0^2}\Big\},
\eeqno
with a constant $C$  not involving $\check c_1$. Then, from \equ{how small gammas} and $\r_1\ge \check c_1\sqrt\k$ there follows 
\beq{nuova}
\hat c\hat E<C\max\Big\{\frac{1}{\g_*},\frac{1}{\check c_1^2},\frac{\k}{\r_0^2}\Big\}<1
\eeq
provided $\g_*$, $\check c_1^2>C$ and $\k<C^{-1}\r_0^2$.
On the other hand, in the case $E=C\m^2\bar K^{2\t+1}\bar\g^{-2}$
\beqano
\hat c\KAM&\le&
C \max\Big\{
\m(\log\k^{-1})^{2\t+1}\bar\g^{-2}\ \Big(\log\Big(\frac{\fg^2}{\m(\log\k^{-1})^{2\t+1}\bar\g^{-2}}\Big)\Big)^{2(\t+1)} \max\{ \su{\g_1^2},\su{\bar\g_2^2}\}, \nonumber\\
&& \frac{\m(\log{\k^{-1}})^{4(\t+1)}}{\bar\g^4}
,\  
\frac{\m(\log\k^{-1})^{2\t+1}\bar\g^{-2}}{\r_1^2},\
\frac{\m(\log\k^{-1})^{2\t+1}\bar\g^{-2}}{\r_0^2}\Big\},
\eeqano
Using now that $\r_1\ge \check c_1\frac{\sqrt\m(\log\k^{-1})^{\t+1/2}}{\bar\g}$ and again the definition of $\bar\g$ in \equ{how small gammas}, we again find an inequality like in \equ{nuova}, but with $\frac{\k}{\r_0^2}$ replaced by $\frac{\sqrt\m}{\r_0^2\g_*^2}$

\nl
 Finally, since the KAM condition $\hat c\hat E<1$ is met,\cite[Proposition 3]{chierchiaPi10} holds in this case. Then, we
can find  a set of invariant tori
$${\mathcal K}_*\subset\bar D_r\times\torus^{n_1}\times\big\{2\r_1< p_i^2+q_i^2<2(\e_0/3)^2\ ,\ \forall\ i\big\}_r\subset( \cP_{\sqrt{2}\e_0/3})_r$$
(with $r<C\bar\g_2$) satisfying the measure estimate
\beqa{meas est}
\meas\big( \cP_{\sqrt{2\check c_2}\e_0}\setminus \cK_*\big)&\le& \meas\big( \cP_{\sqrt{2\check c_2}\e_0})_r\setminus \cK_*\big)\nonumber\\
&\le& C(\bar\g+\g_1+\frac{\bar\g_2}{\e_0^2}+\k^{n_2/4})\meas \cP_{\sqrt{2}\e_0/3}).
\eeqa
We omit  to detail how \equ{meas est} follows from \cite[(34)]{chierchiaPi10}. For example, the reader may easily modify the end of \cite[Proof of Theorem 1.4, Step 6]{chierchiaPi10}. 

\nl
The theorem is so proved with $\cK:=\cK_*\cap \cP_{\e_0/3}$,  $\e_1=\sqrt{2}\e_0/3$, $\k_*$ $:=$ $\min$ $\{C^{-1/4}$ $\sqrt{\r_0}$, $\e_0/(3\sqrt{\check c_1})\}$, $\m_*$ $:=$ $\min$ $\{
C^{-2}\r_0^4\g_*^4
,\ \g_\star^4(\e_0/(3\sqrt{\check c_1}))^8$\}. \qed
\section{The Theorem by N. N. Nehoro{\v{s}}ev}\label{Nekhorossev}
\renewcommand{\theequation}{
\arabic{equation}}
Below is a more technical statement of Theorem \ref{Nek thm simplified}, as it follows from \cite{nehorosev77} and, especially, \cite{nehorosev79}.

\vskip.1in
\noi
The statement in \cite{nehorosev77}--\cite{nehorosev79}  is based on  the notion of  ``steepness'' for a given smooth  function $H_0(I)=H_0(I_1,\cdots, I_{n_1})$ of $n_1$ arguments. We shall adopt the definition given in \cite{nehorosev77}.
 This definition involves a number of parameters, denoted, in \cite{nehorosev77}, as $(g$, $m$, $C_1$, $\cdots$, $C_{n_1-1}$, $\d_1$, $\cdots$, $\d_{n_1-1}$, ${\eufm a}_1$, $\cdots$, ${\eufm a}_{n_1-1})$.
Accordingly, we shall call a given function $(g$, $m$, $C_1$, $\cdots$, $C_{n_1-1}$, $\d_1$, $\cdots$, $\d_{n_1-1}$, ${\eufm a}_1$, $\cdots$, ${\eufm a}_{n_1-1})$--steep, if it is steep with such parameters. See \cite[p. 28 and p. 36]{nehorosev77} for details.

\begin{theorem}[{\cite{nehorosev77}, p. 30; \cite{nehorosev79}}]\label{Nek thm} Let $H=H_0(I)+ P(I,\varphi,p,q)$ be real--analytic on $\cP_\r:=V_\r\times \torus_\r^{n_1}\times B_\r^{2n_2}$ and assume that $I\in V\to H_0(I)$ is $(g$, $m$, $C_1$, $\cdots$, $C_{n_1-1}$, $\d_1$, $\cdots$, $\d_{n_1-1}$, ${\eufm a}_1$, $\cdots$, ${\eufm a}_{n_1-1})$--steep, with $\r<1<m$. Then, one can find $a$, $b\in (0,1)$ and\footnote{We changed a bit notations of \cite{nehorosev77}. 
Let us call $\bar\cP$, $\bar\r$ the quantities that 
in the statement of  \cite[The main theorem, p. 30]{nehorosev77} are called $F$, $\r$ (clearly, $s$, $n$, $H_1$, $G$, $D$ of \cite{nehorosev77} correspond to our $n_1$, $n_2$, $P$, $V$, $B^{2n_2}$). In the statement of  \cite[The main theorem, p. 30]{nehorosev77}, condition \equ{inclusion} is required, with $\cP$ replaced by
 $\bar\cP_{-2r}$, where  $\bar\cP_{-2r}$ is a real set defined   as the   biggest subset  $\cA\subset\bar\cP$ for which $\cA_{2r}\subset \bar\cP$.  Plainly $(\bar\cP_{-2r})_{2r+\bar\r}=\bar\cP_{\bar\r}$. Letting $\cP:=\bar\cP_{-2r}$ and $\r:=2r+\bar\r$ we have our statement. Our condition $M_\star<\r^{1/b}$ corresponds to \cite{nehorosev77} 's assumption $\bar\r>0$.
 }  $0<M_\star<\r^{1/b}$ such that, if \beq{Nek smallness}M:= \sup_{\cP_\r}|\partial P|\in (0,M_\star)\eeq
any trajectory $t\to \g(t)=(I(t), \varphi(t), p(t), q(t))$ solution of $H$ such that \beq{inclusion}(p(t), q(t))\in B^{2n_2}\ ,\quad \forall\ 0\le t\le T:= \frac{1}{M}e^{\frac{1}{M^a}}\eeq
verifies
$$|I(t)-I(0)|\le r:=\frac{1}{2}M^b\qquad \forall\ 0\le t\le T\ .$$
The number $M_\star$ can be taken to be\footnote{See \cite[p. 53]{nehorosev79}. By the previous note, we have to replace $\r$ in \cite[p. 53]{nehorosev79} with $\bar\r:=\r-2r$. Note that condition $M_\star<(\frac{\r}{2})^{1/b}$ implies $\r\ge\r-2r=\r-M^b\ge \r-M_\star^b\ge \frac{\r}{2}$. With this observation, we are allowed to identify  $\r$ of \cite[p. 53]{nehorosev79} with our $\r$. 
Letting then $M_0$, $M_1$ and $M_2$ as in \cite[p. 53]{nehorosev79}, one sees, using the formulae in \cite[pp. 48--57]{nehorosev79}, that $M_1=\frac{c_1}{\r^8 m^4}(\frac{C_{n_1-1}}{g})^p$, while, since $\r<1<m$, $\dst M_2=c_2\min\{\frac{1}{m\r^2}(\frac{\r}{m})^p, \frac{1}{m\r^2}(\frac{C_{n_1-1}}{g})^p, \frac{1}{\r^2 m}(\frac{m}{C_r})^p,\ {\frac{1}{m\r^2}(\frac{1}{\max_r\d_r})}^p,\ \frac{1}{m\r^2}\}$. Therefore, $M_0:=\min\{M_1,\ M_2\}$ verifies the inequality in \equ{M0}.
}
\beq{M*}M_\star=\min\{(\frac{\r}{2})^{1/b}\ ,\quad M_0\}\eeq
where $M_0$ verifies
\beq{M0}M_0\ge \frac{c_0}{\r}\min\{(\frac{C_{n_1-1}}{g})^p,\ (\frac{\r}{m})^p,\ (\frac{m}{C_r})^p,\ (\frac{1}{\max_r\d_r})^p,\ 1\}\eeq
for some $c_0<1<p$ depending only on $n_1$, $n_2$ and ${\eufm a}_1$, $\cdots$, ${\eufm a}_{n_1-1}$.
\end{theorem}

\subsection{Steepness conditions}\label{Steepness conditions}
In \cite{nehorosev77}, a function $H_0=H_0(I)$ of $n_1$ variables $(I_1,\cdots, I_{n_1})$ is called ``quasi--convex'' in $I$ if the system $$\arr{
\dst \sum_{j=1}^{n_1}\partial_{I_j}H_0(I)\eta_j=0\\
\dst \sum_{j, k=1}^{n_1}\partial^2_{I_jI_k}H_0(I)\eta_j\eta_k=0
}$$ has the only trivial solution. Concave or convex functions, having  definite in sign Hessian $\partial^2_{I_jI_k}H_0(I)$, are in particular quasi--convex. Moreover, $H_0$ is said to satisfy the three--jet conditions if, again, the system
 \beq{three jet}\arr{
\dst \sum_{j=1}^{n_1}\partial_{I_j}H_0(I)\eta_j=0\\
\dst \sum_{j, k=1}^{n_1}\partial^2_{I_jI_k}H_0(I)\eta_j\eta_k=0\\
\dst \sum_{j, k, h=1}^{n_1}\partial^3_{I_jI_k, I_h}H_0(I)\eta_j\eta_k\eta_h=0
}\eeq
has the only trivial solution.

\nl
In \cite{nehorosev73} it is proved that quasi--convex functions and functions satisfying the three--jet condition are steep.

\def\cprime{$'$} \def\cprime{$'$}

\end{document}